\documentclass[a4paper,11pt]{amsart}
\input{Macros}
\usepackage[frenchb, english]{babel}
\usepackage{tabu}

\newcommand{\fm}{\mathfrak{m}}

\DeclareMathOperator{\gr}{gr}
\DeclareMathOperator{\Fitt}{Fitt}
\DeclareMathOperator{\NMSing}{Sing^{\text{nm}}}
\DeclareMathOperator{\D}{D}
\DeclareMathOperator{\DM}{DM}

\DeclareMathOperator{\VDir}{VDir}
\DeclareMathOperator{\cl}{cl}

\title{Resolution of 1-foliation singularities on surfaces and threefolds}
\author{Quentin Posva}
\date{}

\address{Mathematisches Institut der Heinriche-Heine-Universität Düsseldorf, Universitätsstr.1, 24.13.03.33}
\email{quentin.posva@hhu.de}

\begin{document}

\maketitle

\begin{quote}
\textsc{Abstract.} We consider resolution of singularities for $1$-foliations on varieties of dimension at most three in positive characteristic. We prove that such singularities can be completely resolved if we allow tame regular Deligne--Mumford stacks as underlying spaces. If one restricts to underlying varieties, we show that $1$-foliations singularities can be simplified into multiplicative ones.
\end{quote}

\tableofcontents

\section{Introduction}
On a normal variety $X$ over a field of positive characteristic $p>0$, a $1$-foliation is a saturated sub-sheaf of $T_X$ that is closed under Lie brackets and $p$-th powers. It follows from their definition that $1$-foliations are uniquely determined by their generic stalks, and so it is relatively easy to produce examples. However they generally acquire complicated singularities at special points, and a desingularization statement was previously known only in the special case $p=2=\dim$ \cite{Hirokado_Singularities_of_mult_derivations_and_Zariski_surfaces}. If $p>2$, then straightforward computations show that already on surfaces we cannot hope to resolve the singularities of the $1$-foliation germ $x\partial_x+\lambda y\partial_y$ with $\lambda\in\bF_p\setminus\{0,1\}$ on $(\bold{0}\in\bA^2_{x,y})$ while preserving the regularity of the underlying surface (see for example \cite[Remark 4.8]{Tziolas_Quotient_by_alpha_p_and_mu_p}). But these $1$-foliations have rather simple singularities (see \cite[\S\S 4--5]{Posva_Singularities_of_quotients_by_foliations}), and one can wonder if they are the only obstruction to a full desingularization statement. One purpose of this article is to show that it is indeed the case, at least for $1$-foliations of corank $1$ on surfaces and threefolds:

\begin{theorem_intro}[\autoref{thm:surface_foliations_resolution} and \autoref{thm:resolution_3fold_rank_2}]\label{thm_intro:schematic_resolution}
Let $X$ be a normal variety of dimension $d\leq 3$ over an algebraically closed field $k$ of positive characteristic, and $\sF$ be a $1$-foliation on $X$. Then there exists a projective birational morphism $f\colon X'\to X$ such that if $s\in \Sing(f^*\sF)$, then there exists \'{e}tale coordinates $x_1,\dots,x_d\in \sO_{X',s}^\text{h}$ such that $f^*\sF\otimes \sO_{X',s}^\text{h}$ is generated by derivations of the form
		$$\frac{\partial}{\partial x_i}\quad \text{or}
		\quad \sum_{i=1}^d\lambda_i x_i\frac{\partial}{\partial x_i} \ \text{  with }\lambda_i\in\bF_p.$$
\end{theorem_intro}

In the case $\dim X=2$, the starting point of the proof is a result due to Giraud \cite{Giraud_Forme_normale_sur_une_surface_en_car_pos} which asserts that a given section $u\in \Gamma(X^{(1)},\sO_{X^{(1)}})$, after a sequence of blow-ups, can be given everywhere \'{e}tale-locally the normal form
		$u=v^p+\bold{x}^\bold{a}$
where $\bold{x}$ is a $p$-basis and $\bold{a}\neq \bold{0}$ modulo $p$. This allows us to describe the normalization of the $p$-cyclic cover $X^{(1)}[\sqrt[p]{u}]$ as a ring of constants for derivations of the desired type (\autoref{prop:inv_subring_normaliz_hypersurface}). Using Jacobson's correspondence \cite[Theorem 2.5.12]{Posva_Singularities_of_quotients_by_foliations}, it follows that the $1$-foliation on $X$ defining the sub-extension $K(X^{(1)})\left[\sqrt[p]{u}\right]\subset K(X)$ is \'{e}tale-locally generated by derivations of the desired form. This settles the surface case. 

The corank $1$ case on threefolds is proved similarly using Cossart's generalization of Giraud's result \cite{Cossart_Forme_normal_de_fonctions_en_dim_3}. We actually need a stronger statement than the main result of \cite{Cossart_Polyhedre_caracteristique}, but which follows quite easily from the same arguments. In the \autoref{appendix} we outline Cossart's methods and explain why they imply the result that we need.

The Jacobson correspondence allows one to deduce some results about the corank $2$ case on threefolds, but if we want to ensure regularity of the underlying scheme we are faced with the problem of lifting derivations along resolutions of linearly reductive quotient singularities. We work out the surface case in \autoref{prop:T-equivariance_toric_surfaces}, relying on earlier work of Wahl and Hirokado, but we leave the threefold case open to an interested reader. See \autoref{rmk:schematic_rk_1_result} for more details. Nevertheless, as we will see below, there is a nicer approach to the corank $2$ case.

\medskip
We can go further than \autoref{thm_intro:schematic_resolution}. It has been realized during recent years that stacks and weighted blow-ups provide powerful tools for resolution of singularities, simplifying and streamlining the resolution algorithms in characteristic $0$ (see \cite{Abramovich_Birational_geometry_using_weighted_blowups} for a survey). Following this trend, one may ask whether weighted blow-ups can help with the resolution of $1$-foliation singularities. Now, local calculations show that weighted blow-ups are exactly the missing piece to obtain a full desingularization for $1$-foliations on surfaces! More precisely, to resolve $x\partial_x +\lambda y\partial_y$ one blow-ups $(x,y)$ with weights $(1,\tilde{\lambda})$, where $\tilde{\lambda}\in \{1,\dots,p-1\}$ reduces to $\lambda$ modulo $p$. This settles at once the surface case. 

The threefold case is more delicate, since the singular locus of the $1$-foliation is in general one-dimensional and the \'{e}tale coordinates given by \autoref{thm_intro:schematic_resolution} are not natural in any helpful way. Nonetheless, if we work with a given set of coordinates, then there is a simple algorithm that resolves the singularities of $1$-foliations of corank $1$ with two steps of weighted blow-ups, and one only has to prove that the Rees algebras defining these blow-ups are independent of the coordinates: this is shown in \autoref{claim:weighted_ideal_indep_of_choices}. 

Once this is done, the corank $2$ should in principle follow at once through the Jacobson correspondence. However at this stage we are dealing with $1$-foliations on Deligne--Mumford stacks, so some foundational work has to be taken care off. We do so in \autoref{section:foliations_on_DM_stacks}: everything follows from the schematic counterparts and from the properties of \'{e}tale morphisms. 

Putting everything together, we obtain:
\begin{theorem_intro}[\autoref{thm:stacky_resolution_surfaces}, \autoref{thm:stacky_resolution_rk_2} and \autoref{thm:stacky_resolution_rk_1}]\label{thm_intro:stacky_resolution}
Let $X$ be a normal variety of dimension $\leq 3$ over a algebraically closed field $k$ of positive characteristic $p>2$, and $\sF$ be a $1$-foliation on $X$. Then there exists a proper birational morphism $f\colon \sX'\to X$ such that 
	\begin{enumerate}
		\item $\sX'$ is a connected separated regular tame Deligne--Mumford stack of finite type over $k$ with diagonal representable by schemes, and
		%\item $f$ is an isomorphism over $\Reg(\sF)$, and
		\item $f^*\sF$ is regular.
	\end{enumerate}
\end{theorem_intro}

For $p=2$ we are able to show that we can take $\sX'$ to be a scheme, hence recovering and generalizing Hirokado's result \cite{Hirokado_Singularities_of_mult_derivations_and_Zariski_surfaces}:
\begin{theorem_intro}[\autoref{thm:resolution_in_char_2}]\label{thm_intro:resolution_char_2}
Let $X$ be a normal variety of dimension $\leq 3$ over a algebraically closed field $k$ of positive characteristic $p=2$, and $\sF$ be a $1$-foliation on $X$. Then there exists a projective birational morphism $f\colon X'\to X$ such that 
	\begin{enumerate}
		\item $X'$ is a regular variety, and
		%\item $f$ is an isomorphism above $\Reg(\sF)$, and
		\item $f^*\sF$ is regular.
	\end{enumerate}
\end{theorem_intro}

\begin{remark}[Blowing-up over $\Reg(\sF)$]
In \autoref{thm_intro:stacky_resolution} and \autoref{thm_intro:resolution_char_2}, if $\dim X=2$ we can ensure that $f$ is an isomorphism over $\Reg(\sF)$ (see \autoref{thm:stacky_resolution_rk_2}). If $\dim X=3$ our method does not guarantee this feature: see \autoref{rmk:u_vs_F} and \autoref{rmk:blow_up_over_reg_locus}.
\end{remark}

In the two-dimensional case, the arguments leading to \autoref{thm_intro:stacky_resolution} can be organized functorially with respect to \'{e}tale morphisms. We obtain:
\begin{theorem_intro}[\autoref{thm:functorial_resolution}]
There exists a function $R(\bullet,\bullet)$, which takes as input pairs $(\sS,\sF)$ where $\sS$ is a separated normal two-dimensional Deligne--Mumford stack of finite type over $k$ and $\sF\subset T_{\sS/k}$ is a $1$-foliation of rank one, and whose output is a proper birational morphism
		$$R(\sS,\sF)\colon \sS'\longrightarrow \sS$$
such that:
	\begin{enumerate}
		\item $\sS'$ is a regular Deligne--Mumford stack,
		\item $R(\sS,\sF)^*\sF$ is regular,
		\item $R(\sS,\sF)$ is an isomorphism over $\Reg(\sF)$.
	\end{enumerate}
Moreover, if $f\colon \sT\to \sS$ is an \'{e}tale morphism, then we have a cartesian diagram
		$$\begin{tikzcd}
		\sT' \arrow[r]\arrow[d, "R(\sT{,}f^*\sF)" left] & \sS'\arrow[d, "R(\sS{,}\sF)"] \\
		\sT\arrow[r, "f"] & \sS.
		\end{tikzcd}$$
\end{theorem_intro}

We refer to \autoref{rmk:several_algorithms} for further remarks on the algorithm $R(\bullet,\bullet)$, and to \autoref{rmk:functorial_resolution_dim_3?} for some issues regarding a similar functorial result in dimension three.

\begin{remark}[Resolution of foliations in characteristic $0$ in dimension $\leq 3$]
Similar resolution statements for foliations on complex varieties are known (\footnote{
	See \cite{ABTW_Principalization_of_log_foliated_orbifolds} and the references therein for the most recent results about resolution of complex foliations in arbitrary dimensions.}), with the caveat that a larger range of final singularities must be allowed. Thus in the light of \autoref{thm_intro:schematic_resolution} and of \autoref{thm_intro:stacky_resolution}, we see that the $p$-closed condition in positive characteristic puts strong restrictions on the possible foliation singularities.

In dimension $2$, resolution of complex foliation singularities is a well-known theorem of Seidenberg. It guarantees that we can simplify the singularities of a line foliation $\sF$ on a surface $S$ until at each $s\in \Sing(\sF)$ with local coordinates $x,y\in\sO_{S,s}$ we can write
		$$\sF\otimes\sO_{S,s}=\sO_{S,s}\cdot (a(x,y)\partial_x+b(x,y)\partial_y)$$
with the property that the eigenvalues $\lambda_1,\lambda_2$ of the Jacobian matrix $\Jac(a,b)(s)$ are not both zero, and that their ratio does not belong to $\bQ_{>0}$. Furthermore, the Poincar\'{e}--Dulac theorem asserts that if the ratio does not belong to $\bQ$, then $\sF\otimes\widehat{\sO}_{S,s}$ is generated by $\lambda_1 \bar{x}\partial_{\bar{x}}+\lambda_2\bar{y}\partial_{\bar{y}}$ for some formal coordinates $\bar{x},\bar{y}\in\widehat{\sO}_{S,s}$. Not every final foliation singularities is formally linearizable (e.g.\ the saddle-node); but as noticed by McQuillan \cite[Proposition II.1.3]{McQuillan_Canonical_models_foliations}, and in accordance with our \autoref{thm_intro:schematic_resolution}, if a complex derivation has $p$-closed reduction for all-but-finitely-many primes $p$ then it is formally linearizable.

In dimension $3$, the complex corank $1$ case is due to Cano \cite{Cano_Reduction_sing_corank_1_foliations_dim_3}. His methods are similar in spirit to those of Cossart \cite{Cossart_Forme_normal_de_fonctions_en_dim_3}, and so our \autoref{thm_intro:schematic_resolution} can be seen as a positive characteristic analogue. 

The corank $2$ case in dimension $3$ is more delicate. By \cite{Cano_Roche_Spivakovsky_Reduction_sing_rank_1_foliations_dim_3}, one can always obtain only ``log elementary" singularities using smooth blow-ups. But as explained in \cite[Fact III.iii.6]{McQuillan_Panazzolo_Almost_etale_resolution_of_foliations} we cannot always obtain elementary singularities (equivalently log canonical singularities, see \cite{McQuillan_Canonical_models_foliations} or \cite[\S 4]{Posva_Singularities_of_quotients_by_foliations}) while preserving the regularity of the underlying variety, as demonstrated by an example of Sanz. It was realized in \cite{McQuillan_Panazzolo_Almost_etale_resolution_of_foliations} that the situation is largely clarified if one works with regular Deligne--Mumford stacks as underlying spaces and weighted blow-ups as modifications: in this set-up one can simplify foliations until they only have canonical singularities. See also \cite{Rebelo_Reis_Resolution_of_rk_1_foliations_on_3folds} for an alternative proof and an extended discussion on that topic.
\end{remark}

\begin{remark}[Resolution of vector fields in positive characteristic]
Our methods rely heavily on Jacobson's correspondence, and thus give no insight into the resolution of non-$p$-closed foliations in positive characteristic. Van den Essen's proof of Seidenberg theorem  for foliations on surfaces works equally well in positive characteristic (see \cite[Theorem 1]{Brunella_Birational_geometry_of_foliations}), but on threefolds essentially nothing is known.
\end{remark}

\subsection{Digression: $1$-foliations and formal groupoids}
A desingularization statement is usually the first step for a Minimal Model Program of some sort. We shall comment briefly on the topic of the MMP for ($1$-)foliations.

From the point of view of birational geometry, a general $1$-foliation in positive characteristic is not be well-behaved \cite{Bernasconi_Counterexample_MMP_for_foliations_in_pos_char}: the cone theorem, the base-point freeness theorem and the existence of Mori fiber spaces might all fail. Thus there is no hope to obtain an MMP for $1$-foliations, while in characteristic $0$ the MMP for foliations has been developed on surfaces \cite{Brunella_Birational_geometry_of_foliations, McQuillan_Canonical_models_foliations, Mendes_Kodaira_dim_sing_foliations} and (in most cases) on threefolds 
\cite{McQuillan_Canonical_models_foliations, Cascini_Spicer_MMP_corank_1_foliations_on_3folds, Cascini_Spicer_MMP_rank_1_foliations_on_3folds, Spicer_Svaldi_Applications_MMP_cork_1_foliations}, and even in greater generality for foliations induced by fibrations \cite{Cascini_Spicer_MMP_for_integrable_foliations, Chen_Han_Liu_Xie_MMP_for_integrable_foliations_and_gen_pairs}.

In characteristic $0$, there are several equivalent points of view on foliations: a first order datum $\sF\subset T_{X/k}$ determines a sub-algebra $\sF^\bullet\subset \Diff_{X/k}$, which in turn can be seen as an infinitesimal equivalence relation on $X$. In positive characteristic, all these notions are different and $1$-foliations (that is, first order data) are the less structured ones. So one can expect that more constrained objects behave better from the perspective of birational geometry.

In an effort to develop a characteristic-free study of foliations, McQuillan introduces in \cite{McQuillan_Formal_groupoids} the notion of \emph{formal groupoids}, which encapsulates the idea of infinitesimal equivalence relation. A formal groupoid on $X$ is an infinitesimal groupoid $g\colon \mathfrak{F}\to X\times_k X$ in the category of formal schemes such that, generically over $X$, both projections $p_i\circ g$ are formally smooth and $g$ factors through an embedding into the completion $\widehat{\Delta}$ of $X\times X$ along its diagonal. This notion leads to a satisfactory MMP for (suitable mildly singular) formal groupoids on surfaces \cite[\S\S 5-6]{McQuillan_Formal_groupoids}. So formal groupoids appear to be more suited to birational geometry than $1$-foliations.

However, apart from smooth dominant morphisms, we currently lack a good source of examples of midly singular formal groupoids in positive characteristic. Any $\mathfrak{F}\rightrightarrows X$ determines a $1$-foliation by truncation to the first order (see \cite[\S 4.5]{Bongiorno_Infinitesimal_groupoids_and_quotient_spaces} for a clear explanation of this procedure), but not every $1$-foliation can be integrated into a formal groupoid. The obstruction already appears for regular $1$-foliations $\sF$: one has to, at the very least, extend $\sF$ into a sub-algebra $\sF^\bullet\subset \Diff_{X/k}$ of the sheaf of Grothendieck differential operators with the property that $\gr(\sF^\bullet)$ is isomorphic to the symmetric divided power algebra $\Gamma(\sF)$ (cf.\ \cite[1.6]{McQuillan_Formal_groupoids}). This is possible locally at any point $s\in X\setminus\Sing(\sF)$: for then by \cite[Lemma 2.4.5]{Posva_Singularities_of_quotients_by_foliations} there are local coordinates $x_1,\dots,x_d\in\sO_{X,s}$ such that $\sF$ is generated by the $\partial_{x_1},\dots,\partial_{x_r}$, and one can integrate $\sF\otimes\sO_{X,s}$ into the sub-algebra of $\Diff(\sO_{X,s})$ generated by the operators $\frac{1}{a!}\frac{\partial^a}{\partial x_i^a}$ $(i\leq r, a\geq 1)$. The problem is that these extensions are not unique: see \cite[Example 5.7]{Grabowski_Foliations_and_Galois_theorey_in_pos_char} for an extension of $\partial_y\in \Der_k(k(x,y))$ into an algebra $\sF^\bullet$ of differential operators that is not generated by $\{\frac{1}{a!}\frac{\partial^a}{\partial y^a}\mid a\geq 1\}$, but which nonetheless satisfies $\gr(\sF^\bullet)\cong \Gamma(\sF^1)$.

A natural question is whether we can generalize our resolution theorems to the setting of formal groupoids. We will not investigate it here, but we make the following remark.

\begin{remark}[Resolution of formal groupoids]\label{rmk:resolution_formal_foliation}
What is meant by resolution of formal groupoids depends on the desired structure at the singular points. Some caution is necessary: as demonstrated in \cite{Staccone_Non_resolution_and_algebraization_of_foliations}, over-optimistic statements do not hold. Indeed, one finds in \emph{op. cit.}, \S 7, an $\bF_2$-surface $S$ together with a formal groupoid $\mathfrak{F}\rightrightarrows S$, induced by a generically smooth fibration $S\to C$ onto a curve, that does not admit any proper birational model $\mathfrak{F}'\rightrightarrows S'$ with the properties that:
	\begin{enumerate}
		\item $S'$ is regular,
		\item the $1$-truncation of $\mathfrak{F}'$ is everywhere log canonical (see \cite[Proposition 3.0.6, Corollary 3.0.7]{Posva_Singularities_of_quotients_by_foliations} for equivalent descriptions), in particular if $s\in \Sing(\mathfrak{F}'^1)$ then there are formal coordinates $x,y$ for $\widehat{\sO}_{S',s}$ such that
				$$\widehat{\sO}_{S',s'}/\mathfrak{F}'^1=k(s)\llbracket x^ay^b\mid a+\lambda_1 b=0 \ (p)\rrbracket$$
		for some $\lambda_1\in\bF_p^\times$, and
		\item such local description extends to higher truncations: if $s\in S'$ is a singular point of $\mathfrak{F}'$ then there are formal coordinates $x,y$ at $s$ and a $p$-adic integer $\lambda$ such that
				$$\text{for every }r\geq 1: \quad \sO_{S',s}/\mathfrak{F}'^r= k(s)\llbracket x^ay^b\mid a+\lambda b=0\ (p^r)\rrbracket.$$
		Here $\mathfrak{F'}^r$ is the $r$-truncation of $\mathfrak{F'}$ affording a factorization $S'\to S'/\mathfrak{F}'^r\to S'^{(r)}$, see \cite[3.1]{Staccone_Non_resolution_and_algebraization_of_foliations} or \cite[4.4.2 and 4.5.6]{Bongiorno_Infinitesimal_groupoids_and_quotient_spaces}.
	\end{enumerate}
%Notice that the second condition (log canonicity of $\mathfrak{F'}^1$) is implied by the $r=1$ instance of the third one (use \autoref{example:ring_of_csts_toric_der} and \autoref{prop:lc_foliation_and_linear_alg}).
\end{remark}

The article is organized as follows. After specifying our conventions and spelling out a useful lemma related to normalization, we consider lift of derivations to resolutions of surfaces in \autoref{section:lift_of_derivations} and we recall the construction of weighted blow-ups in \autoref{section:weighted_blow_ups}. We develop the basic theory of $1$-foliations on normal Deligne--Mumford stacks in \autoref{section:foliations_on_DM_stacks}. Then we proceed to the resolution of $1$-foliations on surfaces in \autoref{section:resolution_surfaces}, and on threefolds in \autoref{section:resolution_3folds}.

\subsection{Acknowledgements}
I am grateful to Andr\'{e} Belotto for several useful discussions and comments, and to Przemys{\l}aw Grabowski for putting in perspective the content of \cite{Grabowski_Foliations_and_Galois_theorey_in_pos_char}. I thank the referees for many useful remarks and suggestions. The author is supported by the grant Postdoc Mobility $\sharp$P500PT/210980 of the Swiss National Science Foundation, and is thankful to Stefan Schröer's group at the Heinrich Heine Universität of Düsseldorf for hosting him.

\section{Preliminaries}
\subsection{Conventions and notations}
Unless stated otherwise, we work over an algebraically closed field $k$ of positive characteristic $p>0$ (the main exception will be \autoref{section:lift_of_derivations}).
	\begin{enumerate}
		\item A \emph{variety} (over $k$) is an integral quasi-projective $k$-scheme of finite type. A \emph{curve} (resp. \emph{surface}, \emph{threefold}) is $k$-variety of dimension one (resp. two, three).
		\item Normalizations of integral Noetherian schemes and algebras are denoted by $(\bullet)^\nu$.
		\item Normalization commutes with completion for excellent rings \cite[7.8.3.(vii)]{EGA_IV.2}: we will use this fact several times without further comment.
		\item Let $X\to S$ be a morphism of $k$-schemes. We let $F_S$ be the absolute Frobenius of $S$. Then we can form the \emph{relative Frobenius} $F_{X/S}\colon X\to X^{(1)}=X\otimes_{S,F_S}S$: it is an $S$-linear morphism. If $F_S$ is invertible, for example in the case $S=\Spec(k)$, we can construct the sequence of $S$-linear morphisms (all denoted by $F_{X/S}$)
				$$\dots\to X^{(-1)}\to X\to X^{(1)}\to \dots$$
		Notice that in this case, the schemes $X^{(n)}$ ($n\in \bZ$) are abstractly (but not $S$-linearly) isomorphic. (The notation for this sequence of morphisms is not uniform in the literature.)
		%Our notation is used in, eg, \cite{Ekedahl_Canonical_models_surfaces_in_pos_char}. 
		%One also finds the notation $X^{(1)}=X^{(p/S)}$.)
		\item The conditions $S_i$ are \emph{Serre's conditions}, see \cite[0341]{Stacks_Project}.
		%\item We use at several places standard MMP terminology for singularities, as defined in \cite[\S 2.3]{Kollar_Mori_Birational_geometry_of_algebraic_varieties}.
		\item Stacks appear at a few places in this article. They will always be separated (which by definition means that its diagonal is proper \cite[D\'{e}finition 7.6]{Laumon_Moret_Bailly_Champs_algebriques}) and of finite type over the base field, normal, with diagonal representable by schemes. In particular, if their inertia is quasi-finite then it is finite. We do not require familiarity with stacks theory beyond, say, \cite[Sections 3--4]{Alper_Stacks_and_moduli}. Most of the time, we will present our stacks as quotients of \'{e}tale groupoids of schemes $[V \rightrightarrows U]$, see \cite[4.3]{Laumon_Moret_Bailly_Champs_algebriques} or \cite[044O]{Stacks_Project}.
		\item If such a stack $\sX$ over $k$ has finite inertia, we say that $\sX$ is \emph{tame} if the stabilizer of any geometric point is a finite linearly reductive group scheme over $k$ \cite{Abramovich_Olsson_Vistoli_Tame_stacks}.
		\item Given a normal variety (or Deligne--Mumford stack) $X$ over $k$, we let 
			$$T_{X/k}=\sHom_{\sO_X}(\Omega_{X/k}^1,\sO_X)$$ 
		be the \emph{tangent sheaf} of $X$. In the literature, this sheaf is sometimes denotes $\Theta_{X/k}$. We usually regard $T_{X/k}$ as the sheaf of $k$-linear derivations of $\sO_X$, denoted $\Der_k(\sO_X)$.
		\item We refer to \cite{Posva_Singularities_of_quotients_by_foliations} for definitions and basic properties of $1$-foliations on schemes. We also refer to \emph{op.cit.}\ for an extended bibliography of the topic of $1$-foliations: many constructions that we use or generalize appear in earlier work of Ekedahl, Miyaoka, Shafarevich--Rudakov, etc.
	\end{enumerate}
	
\subsection{Some lemmas about normalization}
We record several useful results related to normalization.

\begin{proposition}\label{prop:Euclid_algorithm}
Let $B$ be a Noetherian reduced ring with normalization $B^\nu$, and assume that $x,y\in B$ are such that $x^a=y^b$ for coprime positive integers $a$ and $b$. Then there exists an element $u\in B^\nu$ such that $x=u^b$ and $y=u^a$.
\end{proposition}
\begin{proof}
This can be proved using properties of semi-normalization (see \cite[\S 2]{Vitulli_Wk_normality_and_seminormality}), but let us give an elementary proof that is an application of Euclid's algorithm. We may assume that $a>b$. For clarity we change notations and write $a=a_0, b=a_1$ and $x=x_1,y=x_0$. The relation $x^a=y^b$ then becomes 
        $x_1^{a_0}=x_0^{a_1}$.
If we run Euclid's algorithm on $a_0$ and $a_1$, we find a descending sequence
        $a_0>a_1>a_2>\dots a_N>a_{N+1}=(a_0,a_1)=1$
and positive integers $m_1,\dots,m_N$ such that
        $a_{i-1}=m_ia_i+a_{i+1}$ for all $i=1,\dots,N$.
Then in the total field of fractions $Q(B)$ the original equation yields
        $$x_1^{a_2}=\frac{x_0^{a_1}}{x_1^{m_1a_1}}=\left(\frac{x_0}{x_1^{m_1}}\right)^{a_1},$$
so if we let $x_2=x_0/x_1^{m_1}$ we obtain the equation
        $x_1^{a_2}=x_2^{a_1}$.
Continuing this way, we obtain a sequence of elements $x_0,\dots,x_{N+1}$ of $Q(B)$ that satisfies the following relations:
        $$x_{i+1}=\frac{x_{i-1}}{x_i^{m_i}},\quad x_i^{a_{i+1}}=x_{i+1}^{a_i},\quad \forall i=1,\dots, N.$$
In particular, $x_N=x_{N+1}^{a_N}$. We also compute:
        $$x_{N-1}=x_N^{m_N}x_{N+1}=x_{N+1}^{a_Nm_N+1}=x_{N+1}^{a_{N-1}}.$$
More generally, we find by descending induction that
        $x_i=x_{N+1}^{a_i}$ for all $i=0,\dots, N$.
This implies that $x=x_1=x_{N+1}^{a_1}=x_{N+1}^b$ and similarly that $y=x_{N+1}^a$. So $x_{N+1}\in B^\nu$ (\footnote{
	This is clear if $B$ is integral. This also holds for $B$ reduced, by \cite[035F, 035P]{Stacks_Project}.
}) and by taking $u=x_{N+1}$ we obtain the statement.
\end{proof}

\begin{proposition}\label{prop:inv_subring_normaliz_hypersurface}
Let $k$ be an algebraically closed field of characteristic $p>0$. Let $A=k\llbracket x_1,\dots,x_n\rrbracket\supset k\llbracket u_1=x_1^p,\dots,u_n=x_n^p\rrbracket$ and let $J\subset \{1,\dots,n\}$ be a non-empty set of indices. Consider the intermediate extension
	$$k\llbracket \bold{x}\rrbracket \supset
	B=\left(
	\frac{k\llbracket \bold{u}\rrbracket [T]}{
	\left(T^p-\prod_{j\in J}u_j^{a_j}\right)}
	\right)^\nu
	\supset
	k\llbracket \bold{u}\rrbracket$$
where the $a_j$ are positive integers that are all coprime with $p$. Fix an index $j\in J$: then the sub-module
	$$\sF= \sum_{r\notin J} A\cdot \partial_{x_r}
	+\sum_{j'\in J\setminus\{j\}} 
	A\cdot (a_{j'}x_j\partial_{x_j}-a_jx_{j'}\partial_{x_{j'}})
	\subset \Der_k^\text{cont}(A)$$
is a $1$-foliation, and $B=A^\sF$.
\end{proposition}
\begin{proof}
The sub-module $\sF$ is generated by $n-1$ commuting $p$-closed continuous derivations: it is therefore a $1$-foliation. These derivations are linearly independent over the fraction field of $A$, and thus the generic rank of $\sF$ is $n-1$. By the Jacobson correspondence it follows that $[\Frac(A)^\sF:\Frac(k\llbracket \bold{u}\rrbracket)]=p$. 

The ring $B$ can be described as the normalization of the sub-ring $k\llbracket \bold{u}, \prod_{j\in J}x_j^{a_j}\rrbracket$ in its fraction field. Since $\prod_{j\in J}x_j^{a_j}$ belongs to $A^\sF$, and as the latter is normal, we deduce that $B\subseteq A^\sF$. As both $k$-algebras are normal, finite and generically of degree $p$ over $k\llbracket u\rrbracket$, we have $B=A^\sF$ as claimed.
\end{proof}

\begin{lemma}\label{lemma:quotient_is_normalization}
Let $k$ be an algebraically closed field of characteristic $p>0$, let $X=\Spec(A)$ be an affine normal variety over $k$, and let $\sF$ a $1$-foliation of corank one on $X$. Then:
	\begin{enumerate}
		\item $K(X/\sF)$ is purely inseparable of degree $p$ over $K(X^{(1)})$.
		\item If $u\in A^p$ is such that $K(X/\sF)=K(X^{(1)})[u^{1/p}]$, then $X/\sF=X^{(1)}[u^{1/p}]^\nu$.
	\end{enumerate}
\end{lemma}
\begin{proof}
By the Jacobson correspondence \cite[Theorem 2.5.12]{Posva_Singularities_of_quotients_by_foliations} we have
		$$[K(X):K(X/\sF)]=p^{\rk \sF}.$$
As $[K(X):K(X^{(1)})]=p^{\dim X}$ and $\sF$ has corank one it follows that $[K(X/\sF):K(X^{(1)})]=p$. The second part follows from the fact that $X/\sF$ is the normalization of $X^{(1)}$ in $K(X/\sF)$.
\end{proof}
	
\subsection{Lifts of derivations along resolutions}\label{section:lift_of_derivations}
There is no difficulty in pulling back foliations along birational morphisms (see \cite[Construction 2.4.11]{Posva_Singularities_of_quotients_by_foliations}). But whether or not a given derivation lifts to a regular derivation is much more complicated. We shall briefly review some results about this question. \emph{In this section we work over an algebraically closed field $k$ of characteristic $p\geq 0$}.

In case of blow-ups of regular varieties at regular centers, the transformation rules for derivations are easy to understand:
\begin{example}[Pullback along smooth blow-ups]\label{example:smth_blowup_of_foliations}
We explain how foliations pullback along smooth blow-ups on the affine space. %To simplify the notations, if $\xi$ is a coordinate let us write $\frac{\partial}{\partial \xi}=\partial_\xi$.
	\begin{enumerate}
		\item Let $S\to \bA^2_{x,y}$ be the blow-up of the origin. One chart of the blow-up is given by
				$$k[x,y]\longrightarrow k[x',y'],\quad (x,y)\mapsto (x',x'y').$$
		Then $\partial_x(x')=1$ and $\partial_x(y')=-y/x^2=-y'/x'$. 
		This implies that $\partial_x=\partial_{x'}-(y'/x')\partial_{y'}$. 
		Similar computations show that
			$\partial_y=x'\partial_{y'}$.
		A simple way to summarize these identities is the following:
				\begin{equation}\label{eqn:pullback_blowup_surface}
				\begin{cases}
				x\partial_{x} &=x'\partial_{x'}-y'\partial_{y'},\\
				y\partial_{y} &= y'\partial_{y'}.
				\end{cases}
				\end{equation}									
		As the differentials $\partial_x$ and $\partial_y$ generate $T_{\bA^2}$, this shows how vector fields on $\bA^2$ pullback to $S$ on this chart. To describe pullbacks of foliations, we just have to saturate the sub-module generated by the pullbacks of the vector fields in $T_S$. The transformation rules on the other chart are symmetric.
		\item Let $X\to \bA^3_{x,y,z}$ be the blow-up of the origin. One chart of the blow-up is given by
				$$k[x,y,z]\longrightarrow k[x',y',z'],\quad
				(x,y,z)\mapsto (x',x'y',x'z').$$
		After calculations similar to the ones above, we find the transformation rules
				\begin{equation}\label{eqn:pullback_blowup_threefold}
				\begin{cases}
				x\partial_x &= x'\partial_{x'}-y'\partial_{y'}-z'\partial_{z'},\\
				y\partial_y &= y'\partial_{y'},\\
				z\partial_z &= z'\partial_{z'}.
				\end{cases}
				\end{equation}
		The situation on the two other charts are obtained by permuting the roles of the coordinates.
		\item More generally, the reader will check that if $\pi\colon \bA^n_{\bold{y}}\to \bA^n_{\bold{x}}$ is the $y_1$-chart of the blow-up of $(x_1,\dots,x_r)$ for $r\leq n$, which means that we have
				$$(x_1,\dots,x_n)\mapsto 
				(y_1,y_1y_2,\dots,y_1y_r,y_{r+1},\dots,y_n),$$
		then the transformation rules are
				$$\pi^*\partial_{y_1}=\partial_{y_1}-\sum_{i=2}^r\frac{y_i}{y_1}\partial_{y_i}, \quad
				\pi^*\partial_{x_i}=\frac{1}{y_1}\partial_{y_i} \ (1<i\leq r), \quad
				\pi^*\partial_{x_j}=\partial_{y_j} \ (j>r).$$
	\end{enumerate}
\end{example}

It is easy to produce examples of non-liftable derivations: the calculations of \autoref{example:smth_blowup_of_foliations} show that $\partial_x$ does not lift to a regular derivation along $\Bl_{\bold{0}}\bA^2_{x,y}$. The question becomes more interesting if we restrict our attention to those morphisms which resolve singularities.

Classically, resolutions of singularities are obtained by sequences of blow-ups and normalizations. Derivations need not lift along any of these two operations: we record well-known facts in the following remarks.

\begin{remarks}
Let $A$ be an integral ring and $D\in\Der_\bZ(A)$.
	\begin{enumerate}
		\item If $I\subset A$ is an ideal and $D(I)\subseteq I$, then $D$ lifts to a global section of $T_{\Bl_I\Spec(A)/\bZ}$. In fact, let $D_+(a)=\Spec(A[I/a])$ be any standard affine patch of $\Bl_I\Spec(A)$ \cite[052P]{Stacks_Project}. An element of $A[I/a]$ is a fraction $x/a^n$ where $x\in I^n$ for some $n\geq 0$. We compute 
			$$D\left( \frac{x}{a^n}\right)=
			\frac{D(x)}{a^n}+ \frac{(-n)xD(a)}{a^{n+1}}.$$
Since $D(I^n)\subset I^n$ for any $n$, we see that $D(x)\in I^n$ and $xD(a)\in I^{n+1}$, and thus $D(x/a^n)\in A[I/a]$. Hence $D$ lifts to a section of $T_{\Bl_I\Spec(A)}$ over $D_+(a)$, and as $a$ is arbitrary we get our claim.
		\item If $A$ is a finitely generated $\bF_p$-algebra, then $D$ need not preserve the radical ideal defining the singular locus of $A$. Nonetheless, if $D$ defines an $\alpha_p$ or a $\mu_p$-action, then $D$ preserves the Fitting ideals $\Fitt_i(\Omega_{A/\bF_p})$ \cite[Proposition 3.1]{Brion_Schroeer_Inverse_Galois_problem_for_alg_groups}.
		\item If $A$ is a non-normal $\bF_p$-algebra, then $D$ need not lift to $\Der(A^\nu)$: there are examples where $\alpha_p$ or $\mu_p$-actions on $A$ fail to lift to $A^\nu$ (see e.g.\ \cite[Example 5.9]{Brion_Action_fiite_grps_on_curves}).
		\item By contrast, if $A$ is a finitely generated $\bQ$-algebra, then $D$ always preserves the ideal of the singular locus (because it can be integrated into a $\bQ\llbracket t \rrbracket$-automorphism of $A\llbracket t\rrbracket$ using exponentiation) and lifts to the normalization (essentially by the same token, see \cite{Seidenberg_Derivation_and_integral_closure}).
	\end{enumerate}
\end{remarks}

This leads us to the following definition. To state it, observe that given a morphism of $k$-schemes $\pi\colon X\to \Spec(A)$ such that $\pi_*\sO_X=A$, we have a natural map $\pi_*\colon H^0(X,T_{X/k})\to \Der_k(A)$. Note that if $\pi$ is birational then $\pi_*D=0$ implies that $D=0$ at the generic points of $X$: so $\pi_*$ is injective as soon as $\pi$ is birational and $T_{X/k}$ torsion-free.

\begin{definition}
Let $A$ be a normal essentially of finite type $k$-algebra admitting a birational projective morphism $\pi\colon X\to \Spec(A)$ such that $X$ is regular. We say that $\pi$ is a $\mathbf{T}$\textbf{-equivariant resolution} if the natural injective map $\pi_*\colon H^0(X,T_{X/k})\to \Der_k(A)$ is bijective.
\end{definition}

%\begin{remark}\label{rmk:derivation_lifts_after_scaling}
%If $X\to \Spec(A)$ is an arbitrary birational morphism of finite type, then $D\in \Der_k(A)$ only lifts to a \emph{rational derivation} of $X$, meaning that $D$ is a global section of $\Der_k(\sO_X,\Frac(A))$. If $\Spec(B)\subset X$ is an affine open subset, then as $B$ is a finite type $A$-algebra with the same fraction field we can write $B=A[s_1,\dots,s_n]$ with $s_i=a_i/b_i\in \Frac(A)$. Then clearing denominators we see that $\left(\prod_i b_i^2\right)\cdot D$ is a $k$-derivation of $B$.
%\end{remark}

\begin{remark}[On terminology]
In some sources such as 
\cite{Burns_Wahl_Local_contrinutions_to_global_deform_surfaces, Wahl_Vanishing_theorems_for_resolutions_of_surfaces, Wahl_Equisingular_deformations_surfaces, Hirokado_Further_evaluation_of_Wahl_vanishing}, what I have baptised a $T$-equivariant resolution is simply called an equivariant resolution. This is unfortunate, because confusion arises when a group $G$ is acting on the singularity $A$ and we look for a $G$-equivariant resolution. 

For example, if $S$ is a surface singularity with minimal resolution $\pi^m\colon S^m\to S$, then the action of any abstract group on $S$ lifts to $S^m$: this is because $S^m$ is obtained by repeatedly normalizing and blowing-up the reduced singular locus, and the abstract group action lifts along both operations. But in positive characteristic $\pi^m$ might not be $T$-equivariant, as we will see shortly.
\end{remark}

In characteristic $0$, post-Hironaka algorithms for resolution of singularities put an emphasis on functoriality. Equivariance for algebraic group actions is a consequence of functoriality with respect to smooth morphisms \cite[3.9.1]{Kollar_Lectures_on_resolution_of_sing}, which is indeed a feature of some well-known algorithms \cite{Bierstone_Milman_Functoriality_in_resolution}. $T$-equivariance is more subtle\footnote{Any derivation can be integrating into a formal action of the formal group $\widehat{\bG}_a$, but stronger functoriality properties are needed to handle this case (the action morphism is formally smooth but not smooth). Equivalently, one need to ensure that the steps of the algorithm blow-up ideals that are preserved by a given set of derivations. No proof is yet written down, but to my understanding it is in reach of the current techniques, see \cite{Temkin_Relative_and_log_resolution_of_singularities}
and \cite{ABTW_Principalization_of_log_foliated_orbifolds}.}.

\medskip
Let us concentrate on surface singularities. For normal surfaces in characteristic $0$, the minimal good resolution always provides a $T$-equivariant resolution \cite[Proposition 1.2]{Burns_Wahl_Local_contrinutions_to_global_deform_surfaces}. However, this is not true in positive characteristic, even for RDP's; we follow the analysis of \cite{Wahl_Vanishing_theorems_for_resolutions_of_surfaces}. Let $A$ be a germ of normal surface over an algebraically closed field of characteristic $p>0$, and let $\pi\colon X\to \Spec(A)$ be its minimal resolution with exceptional curve $E=\bigcup_i E_i$. We have an exact sequence of $\sO_X$-modules
		$$0\to S\to T_{X/k}\to \bigoplus_i N_{E_i/X}\to 0$$
where $S$ is a rank $2$ vector bundle, see \cite[Proposition 1.7]{Wahl_Vanishing_theorems_for_resolutions_of_surfaces}. Using the long exact sequence with (local) cohomology for $S$ associated to $U=X\setminus E$, one finds \cite[Proposition 1.12]{Wahl_Vanishing_theorems_for_resolutions_of_surfaces} an exact sequence
		\begin{equation}\label{eqn:T-equivariance}
		0\to H^0(X,T_{X/k})\overset{\pi_*}{\longrightarrow} \Der_k(A)\to H^1_E(X,S).
		\end{equation}
So $H^1_E(X,S)=0$ is a sufficient condition for $T$-equivariance of $\pi$. If $A$ is an RDP (\footnote{I.e.\ a rational double point, equivalently a canonical singularity in the sense of the MMP.}), then vanishing of $H^1_E(X,S)$ is also a necessary condition for $T$-equivariance, by \cite[Theorems 4.1 and 5.1]{Hirokado_Further_evaluation_of_Wahl_vanishing}. 

More generally, for toric surface singularities we have the following result:
\begin{proposition}\label{prop:T-equivariance_toric_surfaces}
Let $k$ be an algebraically closed field of characteristic $p>0$, and $A$ be a germ of normal surface singularity over $k$ with minimal good resolution $\pi\colon X\to \Spec(A)$. Assume that the exceptional divisor $E$ is a chain of rational curves\footnote{This includes $\mu_p$-quotient of regular surface germs, see \cite[Theorem 2.3]{Hirokado_Singularities_of_mult_derivations_and_Zariski_surfaces}.}. Then:
	\begin{enumerate}
		\item If $\widehat{A}\cong k\llbracket x,y,z\rrbracket /(z^{n+1}-xy)$ is an $A_n$-singularity with $p|n+1$, then $\dim_k H^1_E(X,S)=1$ with generator the image of $\partial_z$.
		\item In every other case $H^1_E(X,S)=0$ and so $\pi$ is $T$-equivariant.
	\end{enumerate}
\end{proposition}
\begin{proof}
First assume that $A$ is a RDP: then it is an $A_n$-singularity. If $p$ does not divide $n+1$ then $H^1_E(X,S)=0$ by \cite[Theorem 5.17]{Wahl_Vanishing_theorems_for_resolutions_of_surfaces}. If $p|n+1$ then, by \cite[Theorems 4.1 and 5.1]{Hirokado_Further_evaluation_of_Wahl_vanishing}, the sequence \autoref{eqn:T-equivariance} is a short exact sequence with $\dim_k H^1_E(S)=1$. To identify the derivation generating this group, let us work with the model $\widehat{A}\cong k\llbracket x,y,z\rrbracket /(z^{n+1}-xy)$. Notice that $\partial_z\in \Der_k^\text{cont}(k\llbracket x,y,z\rrbracket)$ sends the ideal $(z^{n+1}-xy)$ to itself, and so it descends to a (continuous) derivation of $\widehat{A}$. I claim that $\partial_z$ does not lift along $\pi$. Indeed, we can factor $\pi$ as
		$$X\overset{\varphi}{\longrightarrow} Y \overset{b}{\longrightarrow} \Spec(\widehat{A})$$
where $b$ is the blow-up of the singular point. If $\partial_z$ lifts to an element $D$ of $H^0(X,T_{X/k})$, then $\varphi_*D$ is a lift of $\partial_z$ to $Y$. So it is sufficient to prove that $\partial_z$ does not lift to $Y$. Suppose it does: then it must restrict to a derivation of the affine blow-up algebra $\widehat{A}[\fm/z]$, where $\fm=(x,y,z)$ is the ideal of the singular point. But we compute
		$$\partial_z\left( \frac{x}{z}\right) =\frac{-x}{z^2}
		\notin \widehat{A}[\fm/z],$$
and therefore $\partial_z$ does not lift to $Y$ and neither to $X$. It follows that $H^1_E(S)$ is generated by the image of $\partial_z$ through $\Der_k(A)\to H^1_E(X,S)$. Notice that this implies that $z\partial_z$ lifts to $H^0(X,T_{X/k})$.

Next assume that $A$ is not an RDP. We will use \cite[Proposition 2.14]{Wahl_Vanishing_theorems_for_resolutions_of_surfaces} to deduce that $H^1_E(S)=0$. For this, we need to show that there does not exist any \emph{bad cycle} $Z$ in the sense of \cite[(2.13.1)]{Wahl_Vanishing_theorems_for_resolutions_of_surfaces}, i.e.\ any cycle $Z>0$ supported on $E$ such that
		$$Z\cdot E_i<0\quad\Rightarrow\quad
		Z\cdot E_i=-1\text{ and }  t_{E_i}=1$$
where $t_{E_i}=E_i\cdot (E-E_i)$ is the number of intersections of $E_i$ with the rest of $E$ (so in our case, $t_{E_i}\in \{1,2\}$). Suppose for the sake of contradiction that there exists such a bad cycle $Z$: then by \cite[Theorem 4.3]{Wahl_Vanishing_theorems_for_resolutions_of_surfaces}, as $A$ is not an RDP, there exists an irreducible component $G$ of $E$ (\footnote{
		The component $G$ is given in \cite[Theorem 4.3.(3)]{Wahl_Vanishing_theorems_for_resolutions_of_surfaces}. In the notation of \emph{op.cit.}\ it is of the form $s(L)$ which implies that it is irreducible, cf.\ \cite[(3.2)]{Wahl_Vanishing_theorems_for_resolutions_of_surfaces}.
}) such that 
	\begin{equation}\label{eqn:def_prop_of_G}
	Z\cdot G<t_G-2.
	\end{equation} 
Now, since $E$ is a string we have $t_G\leq 2$ and thus $Z\cdot G<0$. Therefore $Z\cdot G=-1$ and $t_G=1$ as $Z$ is a bad cycle. But using \autoref{eqn:def_prop_of_G} we get $-1=Z\cdot G<t_G-2=-1$ which is a contradiction.
\end{proof}

Functorial and combinatorial resolutions of toroidal singularities are known in every dimension over any perfect field, see \cite[Theorem 1.0.1]{Wlodarczyk_Functorial_resolution_except_for_toroidal_locus}. So in principle one could play a similar game in dimension three and higher; but we will not attempt to do so here.

\subsection{Weighted blow-ups}\label{section:weighted_blow_ups}
Stacks will enter the picture via weighted blow-ups, which we quickly review now. We refer to \cite{Quek_Rydh_Weighted_blowups} for details and proofs.

Let $A$ be a ring. A \textbf{Rees algebra} over $A$ is a finitely generated graded sub-$A$-algebra $I_\bullet$ of $A[t]=\bigoplus_{i=0}^\infty A\cdot t^i$ such that $I_0=A$ and $I_{n}\cdot t\supset I_{n+1}$ for every $n$. (Here $t$ is a dummy variable used to keep track of degrees.) The grading gives a $\bG_m$-action, and one can form the \textbf{weighted blow-up of $A$ along $I_\bullet$}, defined by
		$$\Bl_{I_\bullet}\Spec(A)=\left[
		\Spec_A (I_\bullet)\setminus V(I_+)/\bG_m\right].$$
It is a tame algebraic stack. The natural morphism $\Bl_{I_\bullet}\Spec(A)\to \Spec(A)$ is proper \cite[Proposition 1.6.1]{Quek_Rydh_Weighted_blowups}, but usually not projective (see \cite[Corollary 1.6.2]{Quek_Rydh_Weighted_blowups}).

An important instance of this construction is the following. For $x\in A$ and $d\geq 1$ let $(x,d)$ be the graded ideal of $A[t]$ defined by 
		$$(x,d)=A\oplus (x)t\oplus \dots 
		\oplus (x) t^{d}\oplus (x^2)t^{d+1}\oplus \dots $$
or more succinctly by $(x,d)_m=(x^{\lceil m/d\rceil})t^m$. Given a collection $\bold{x}=(x_1,\dots,x_n)$ of elements of $A$ and positive integers $\bold{d}=(d_1,\dots,d_n)$, we let $(\bold{x},\bold{d})$ be the Rees algebra defined by
		\begin{equation}\label{eqn:weighted_ideals}
		(\bold{x},\bold{d})_m=\sum_{m_1+\dots+m_n=m}\left( \prod_{i=1}^n(x_i,d_i)_{m_i}\right).
		\end{equation}
Notice that if $d_i=1$ for every $i$, then $(\bold{x},\bold{d})_m=(x_1,\dots,x_n)^m$ and so $\Bl_{I_\bullet}\Spec(A)=\Bl_I\Spec(A)$ is the usual blow-up of $A$ along $I$.
		
Weighted blow-ups commute with flat extensions \cite[Corollary 3.2.14.(iii)]{Quek_Rydh_Weighted_blowups}. More precisely, suppose that $A\to A'$ is a flat ring extension. Then for any Rees algebra $I_\bullet$ over $A$, there is a natural isomorphism
		$$\Bl_{I_\bullet}\Spec(A)\times_{\Spec(A)}\Spec(A')\cong 
		\Bl_{I_\bullet\otimes_A A'}\Spec(A').$$
Hence these constructions sheafify in the \'{e}tale topology, allowing us to define weighted blow-ups of schemes and of Deligne--Mumford stacks along Rees algebra sheaves in the expected way.
		
%If $A$ is normal and excellent and $I_\bullet$ is a Rees algebra, then the integral closure $\IC(I_\bullet)$ of $I_\bullet$ is again a Rees algebra \textbf{reference}. Then $\Bl_{\IC(I_\bullet)}\Spec(A)$ is called the \textbf{normalized weighted blow-up of $A$ along $I_\bullet$}.

\begin{remark}\label{rmk:separated_stacks}
Assume that $X$ is a separated scheme or algebraic stack, and let $\sI_\bullet$ be a Rees algebra on $X$. Then $\Bl_{\sI_\bullet}X$ is separated. Indeed, the morphism $\Bl_{\sI_\bullet}X\to X$ is proper, hence separated, and as separated morphisms of stacks are stable under composition \cite[Remarque 7.8.1.(3)]{Laumon_Moret_Bailly_Champs_algebriques}, we deduce that $\Bl_{\sI_\bullet}X$ is separated. 

The stacks $\sX'$ appearing in \autoref{thm_intro:stacky_resolution} will be obtained as successive weighted blow-ups of varieties. In particular they will automatically be separated.
\end{remark}

\begin{proposition}\label{prop:properties_weighted_blow_ups}
Suppose that $A$ is regular and that $I_\bullet=(\bold{x},\bold{d})$ where $\bold{x}$ is a regular sequence. Then:
	\begin{enumerate}
		\item $\Bl_{I_\bullet}\Spec(A)$ is a regular tame Artin stack;
		\item if every $d_i$ is invertible in $A$, then $\Bl_{I_\bullet}\Spec(A)$ is Deligne--Mumford;
		\item the natural morphism $\Bl_{I_\bullet}\Spec(A)\to \Spec(A)$ is an isomorphism above the complement of $V(\bold{x})$;
		\item the coarse moduli space of $\Bl_{I_\bullet}\Spec(A)$ is the scheme $\Proj_A(I_\bullet)$.
	\end{enumerate}
\end{proposition}
\begin{proof}
The statement about the coarse moduli space is \cite[Proposition 1.6.1.(iii)]{Quek_Rydh_Weighted_blowups}, and given the compatibility with flat base-change, the third statement is an immediate consequence. The first two statements follow from the explicit local description \cite[Corollary 5.2.5]{Quek_Rydh_Weighted_blowups}, which in the case of invertible weights is recalled in \autoref{example:weighted_blow_ups} below.
\end{proof}

The coarse moduli space $\Proj_A(I_\bullet)$ is usually singular, and one of the advantages of the stacky point of view is that we may continue to assume that the ambient space retains the regularity properties of $A$, if we are willing to work equivariantly.

Our use of weighted blow-ups will be limited to cases of the following nature.

\begin{example}[Cf.\ {\cite[Corollary 5.2.5]{Quek_Rydh_Weighted_blowups}}]\label{example:weighted_blow_ups}
Consider $A=k[x_1,\dots,x_r]$ with $k$ algebraically closed, and let $d_1,\dots,d_n$ $(n\leq r)$ be positive integers that are invertible in $k$. The weighted blow-up $\Bl_{\sum_{i=1}^n (x_i,d_i)}\Spec(A)$ is covered by affine quotient stacks $D_+(x_i)$ for $i=1,\dots,n$, which can be described as follows:
		$$D_+(x_i)=
		[ \Spec k[u,v_j\mid j=1,\dots, \hat{i},\dots,r]/\mu_{d_i}]$$
with action 
		$$u\mapsto \xi u,\quad
		v_j\mapsto \xi^{-d_j} v_j \ (j\leq n),\quad
		v_j\mapsto v_j \ (j>n)$$
where $\xi\in k$ is a primitive $d_i$-th root of unity. The $k[x_1,\dots,x_r]$-algebra structure is given by
		$$x_i=u^{d_i},\quad x_j=u^{d_j}v_j \ (1\leq j\leq n, \ j\neq i), \quad
		x_j=v_j \ (j>n).$$
We also record how the derivations $\partial_{x_j}$ pullback. One checks that
		$$\partial_{x_i}=\frac{1}{d_i u^{d_i}}\left[
		u\partial_u -\sum_{j\in \{1,\dots, \hat{i},\dots,n\}} d_jv_j\partial_{v_j}\right],\quad
		\partial_{x_j}=\frac{1}{u^{d_j}}\partial_{v_j} \ (1\leq j\leq n, \ j\neq i),\quad
		\partial_{x_j}=\partial_{v_j} \ (j>n).$$
All these derivations are invariant by the $\mu_{d_i}$-action on $\Der_k(k[u,v_j\mid j\neq i])$ (which is spelled out in \autoref{example:foliation_on_qt_stack} below), and thus they descend to the quotient stack $D_+(x_i)$. 

Notice that if, say, $d_i$ was not invertible in $k$, then no non-zero scaling of $\partial_{x_i}$ would lift to the blow-up.
\end{example}

The following result on Rees algebras will be very useful:

\begin{proposition}\label{prop:Rees_algebra_is_functorial}
Let $k$ be an algebraically closed field of characteristic $p>0$, and let $(\sO,\fm,\kappa )$ be a complete regular local $k$-algebra of dimension two. Let $\Lambda\in \{1,\dots,p-1\}$ and let $\lambda\in\bF_p^\times$ be its reduction modulo $p$. Let $x,y$ and $u,v$ be two systems of parameters of $\sO$. Assume that $x\partial_x+\lambda y\partial_y$ and $u\partial_u+\lambda v\partial_v$ 
%generate the same sub-module of $\Der^\text{cont}_k(\sO)$ and 
define the same $\kappa$-linear endomorphism of $\fm/\fm^2$. Then $(x,1)+(y,\Lambda)=(u,1)+(v,\Lambda)$.
\end{proposition}
\begin{proof}
For $a,b\in \fm$ we write $I(a,b)=(a,1)+(b,\Lambda)$. Then we have to show the equality of Rees algebra $I(x,y)=I(u,v)$. By symmetry, it suffices to prove one containment.

If $\Lambda=1$ then $I(u,v)_m=\fm^m=I(x,y)_m$ for every $m\geq 0$, so we may assume that $\lambda\in\bF_p\setminus \{0,1\}$ for the rest of the proof.

We look for a power series expansion of $u,v$ in terms of $x,y$. Write $\xi=x\partial_x+\lambda y\partial_y$ and $\xi'=u\partial_u+\lambda v\partial_v$. Then $\xi$ and $\xi'$ both define $\mu_p$-actions on $\sO$ that fix the closed point \cite[Proposition 2.3.7, Lemma 2.3.10]{Posva_Singularities_of_quotients_by_foliations}. As $\mu_p$ is diagonalizable, these actions are linearizable: they are uniquely determined by the induced actions on the cotangent space $\fm/\fm^2$ via the canonical isomorphism
		$$\Sym^\bullet_{\kappa}(\fm/\fm^2)\cong \sO.$$
These actions on the cotangent space are precisely given by the images of $\xi$ and $\xi'$ in $\End_{\kappa}(\fm/\fm^2)$. By assumption these two endomorphisms are equal; thus the two $\mu_p$-actions are the same, which means that $\xi=\xi'$.

From $u=\xi(u)$ and $\lambda v=\xi(v)$ it follows that we have the power series expansion
		$$u=\sum_{i\geq 0}\sum_{n\geq 0} \alpha_{i,n} x^iy^{\varphi(i)+np},\quad
		v=\sum_{j\geq 0}\sum_{n\geq 0} \beta_{n,j}x^{\psi(j)+np}y^j$$
in $\sO\cong \kappa \llbracket x,y\rrbracket$, with the functions $\varphi,\psi\colon \bN\to \{0,\dots,p-1\}$ defined by
		$$i+\varphi(i)\lambda =1,\quad \psi(j)+j\lambda=\lambda\quad \text{in }\bF_p.$$
%Notice that $i$ and $\varphi(i)$ cannot be simultaneously zero; and similarly for $j$ and $\psi(j)$.
Since $(u,v)=\fm_\eta=(x,y)$, we must have a linear term in the expansions of $u$ and $v$. Because $\lambda\neq 1$, we cannot have $i=0,\varphi(i)=1$ nor $j=0,\psi(j)=1$. So we obtain that
		$$u=\alpha_{1,0}x+\sum_{i\geq 1}\sum_{n\geq 0} \alpha_{i,n} x^iy^{\varphi(i)+np}, \quad
	v=\beta_{0,1}y+\sum_{j\geq 1}\sum_{n\geq 0} \beta_{n,j}x^{\psi(j)+np}y^j$$
where $\alpha_{1,0},\beta_{0,1}\in \kappa^\times$. We can rewrite this as
		$$u=\alpha x+ y^{d}f,\quad v=\beta y+x^{\Lambda}g$$
where $\alpha,\beta\in \sO^\times$ and $f,g\in \sO$, and where $d\in \{1,\dots,p-1\}$ reduces to $\lambda^{-1}$ modulo $p$. %and where $\Lambda\in\{1,\dots,p-1\}$ reduces to $\lambda$ modulo $p$.

Since $I(u,v)=I(\alpha^{-1}u,\beta^{-1}v)$, we may assume that $\alpha=\beta=1$. So we reduce to show that $I(x+y^df, y+x^\Lambda g)\subseteq I(x,y)$.

Consider the parts of degree $m$: the ideal $I(x+y^df, y+x^\Lambda g)_m$ of $\sO$ is generated by the elements
				$$(x+y^df)^{m-i}\cdot (y+x^\Lambda g)^{\lceil i/\Lambda \rceil}=
				\left(\sum_{s=0}^{m-i}\binom{m-i}{s}x^{m-i-s}(y^df)^s\right)
				\left(\sum_{t=0}^{\lceil i/\Lambda\rceil} \binom{\lceil i/\Lambda\rceil}{t}y^t (x^\Lambda g)^{\lceil i/\Lambda \rceil -t}\right),$$
	for $i=0,\dots,m$. We are going to show that each one of these elements belong to $I(x,y)_m$. It suffices to show that each possible product on the right-hand side belongs to $I(x,y)_m$. So fix $s\in \{0,\dots,m-i\}$ and $t\in \{0,\dots,\lceil i/\Lambda\rceil\}$, and observe that
			\begin{equation}\label{eqn:combinatorics}
			x^{m-i-s}(y^df)^sy^t (x^\Lambda g)^{\lceil i/\Lambda \rceil -t}\in \left( x^{m-s-\Lambda t+\epsilon_\Lambda(i)} y^{t+ds}\right)
			\end{equation}
	where the function $\epsilon_\Lambda(\bullet)\colon \bN\to \{0,\dots,\Lambda-1\}$ is defined by
			$$\epsilon_\Lambda(a)=\Lambda\cdot \lceil a/\Lambda\rceil -a.$$
	The ideal on the right-hand side of \autoref{eqn:combinatorics} is included in $I(x,y)_m$ if we can exhibit an $j\in\{0,\dots,m\}$ such that
			\begin{equation}\label{eqn:combinatorics_2}
			\begin{cases}
			m-s-\Lambda t+\epsilon_\Lambda(i) &\geq  m-j \\
			t+ds & \geq \lceil j/\Lambda\rceil.
			\end{cases}
			\end{equation}
	The first inequality always has a solution $0\leq j\leq m$, for
			$$s+\Lambda t-\epsilon_\Lambda (i)\leq m-i+\Lambda\cdot \lceil i/\Lambda\rceil -\epsilon_\Lambda(i)=m.$$
	So let us choose $j=\max\{0,s+\Lambda t-\epsilon_\Lambda(i)\}$. We claim that the second inequality of \autoref{eqn:combinatorics_2} is satisfied. This is clear if $j=0$. If $j>0$, we argue by contradiction: if the inequality is not satisfied then
			\begin{eqnarray*}
			t+ds &<& \left\lceil \frac{s+\Lambda t-\epsilon_\Lambda(i)}{\Lambda}\right\rceil \\
			&=& t+\left\lceil \frac{s-\epsilon_\Lambda(i)}{\Lambda}\right\rceil\\
			&\leq & t+s,
			\end{eqnarray*}
	which is a blatant contradiction.
	
	Therefore, a generating set of $I(x+y^df, y+x^\Lambda g)_m$ is contained in $I(x,y)_m$. As $m$ is arbitrary, we obtain that $I(x+y^df, y+x^\Lambda g)\subseteq I(x,y)$. This proves the proposition.
\end{proof}

\section{1-Foliations on Deligne--Mumford stacks}\label{section:foliations_on_DM_stacks}
In this section, we generalize the notions of $1$-foliations and of induced quotients (as in \cite[\S 2.4--5]{Posva_Singularities_of_quotients_by_foliations}) to Deligne--Mumford stacks. Throughout the section, $k$ denotes an algebraically closed field of characteristic $p>0$.

\subsection{Preliminaries and definitions}
Let $f\colon V\to U$ be an \'{e}tale morphism of normal $k$-schemes of finite type. Then the canonical map $df\colon f^*\Omega_{U/k}\to \Omega_{V/k}$ is an isomorphism. Dualizing, we obtain a canonical isomorphism 
	\begin{equation}\label{eqn:tangent_sheaf_and_etale}
	\partial f\colon T_{V/k}\cong \sHom_{\sO_U}(f^*\Omega_{U/k},\sO_U)\cong f^*T_{U/k}
	\end{equation}
where the second isomorphism follows from the flatness of $f$ \cite[0C6I]{Stacks_Project}. Now suppose that $\iota\colon G\subset T_{U/k}$ is a coherent sub-module. Again by flatness of $f$, the natural map $f^*\iota\colon f^*G\to f^*T_{U/k}$ is injective, and so we shall identify $f^*G$ with its image in $T_{V/k}$.

\begin{lemma}\label{lemma:pullback_under_etale}
$G$ is a $1$-foliation if and only if $f^*G$ is a $1$-foliation.
\end{lemma}
\begin{proof}
Since $f^*(T_{U/k}/G)=T_{V/k}/f^*G$ and $f$ is \'{e}tale, we see that $T_{U/k}/G$ is $S_1$ is and only if $T_{V/k}/f^*G$ is $S_1$. Assume these sheaves are $S_1$. Then whether $G$ and $f^*G$ are closed under Lie brackets and $p$-th powers can be seen at the generic points \cite[Remark 2.4.2]{Posva_Singularities_of_quotients_by_foliations}: hence we may assume that $T_{U/k},G,T_{V/k}$ and $T_{U/k}/F$ are locally free. Consider the map
		$$L_G\colon \bigwedge^2 G\to T_{U/k}/G, \quad v\wedge w\mapsto [v,w] \mod G,$$
and define similarly $L_{f^*G}$. Since we have a natural isomorphism 
	$$f^*\left(\bigwedge^2G\right)\cong \bigwedge^2f^*G,\quad f^*(v\wedge w)\mapsto f^*v\wedge f^*w,$$ 
we see that $L_G=0$ if and only if  $L_{f^*G}=0$. So $G$ is closed under Lie brackets if and only if $f^*G$ is. Assuming they both are, we get that the map
			$$P_G\colon F^*_{U/k}G\to T_{U/k}/G,\quad v\mapsto v^{[p]} \mod G$$
is $\sO_U$-linear, by Hochschild's formula and \cite[(15) p.209]{Jacobson_Abstract_derivations_and_Lie_algebras}. We define similarly the $\sO_V$-linear map $P_{f^*F}$. By naturality of the Frobenius we see that $f^*P_G=P_{f^*G}$, and it follows that $G$ is $p$-closed if and only if $f^*G$ is so.
\end{proof}

Thus we can make the following definition:
\begin{definition}
Let $f\colon V\to U$ be an \'{e}tale morphism of normal $k$-schemes of finite type, and $\sF\subset T_{U/k}$ be a $1$-foliation. Then we let $f^*\sF\subset T_{V/k}$ be the pullback of $\sF$ as an $\sO_U$-module, and we regard it as an $1$-foliation on $V$.
\end{definition}

Now let $\sX$ be a separated normal Deligne--Mumford stack of finite type over $k$ with diagonal representable by schemes (cf.\ our conventions). Choose an \'{e}tale presentation $U\to \sX$ with associated groupoid $s,t\colon U\times_\sX U\rightrightarrows U$ (here $s,t$ denote the two projections). In general (\cite[\S 4.1]{FGA_explained} or \cite[Exerc. 4.1.18]{Alper_Stacks_and_moduli}), given a coherent $\sO_U$-module $M$, a \emph{descent datum} for $M$ with respect to $s,t$ is an isomorphism $\alpha\colon s^*M\cong t^*M$ that satisfies the equality
		$$p_{23}^*\alpha\circ p_{12}^*\alpha=p_{13}^*\alpha
		\quad \text{(the cocycle condition)},$$
where the $p_{ij}\colon U\times_\sX U\times_\sX U \to U\times_\sX U$ are the projection morphisms. Since all descent data are effective for $U\times_\sX U\rightrightarrows U$, such data defines uniquely a coherent $\sO_\sX$-module $\mathcal{M}$ satisfying $(U\to\sX)^*\mathcal{M}=M$. Conversely, any coherent $\sO_\sX$-module arises in this way.

The tangent sheaf $T_{U/k}$ has a natural descent datum, given by
		$$\begin{tikzcd}
		(\partial t)\circ (\partial s)^{-1}\colon 	
		s^*T_{U/k} \arrow[r, "\sim"] &
		T_{U\times_\sX U/k} \arrow[r, "\sim"]
		& t^*T_{U/k}
		\end{tikzcd}$$
where $\partial t$ and $\partial s$ are defined as in \autoref{eqn:tangent_sheaf_and_etale}. We denote the corresponding $\sO_\sX$-module by $T_{\sX/k}$: it is the tangent sheaf of $\sX$ over $k$. Descending a coherent sub-module $G\subset T_{U/k}$ to $\sX$ \emph{as a sub-module of $T_{\sX/k}$} is equivalent to the condition $(\partial s)^{-1}(s^*G)=(\partial t)^{-1}(t^*G)$ inside $T_{U\times_\sX U/k}$, which we abbreviate as $s^*G=t^*G$.

\begin{definition}
In the above notations, an \textbf{$1$-foliation} $\sF$ on $\sX$ is a $1$-foliation $\sF_U\subset T_{U/k}$ satisfying $s^*\sF_U=t^*\sF_U$ as sub-modules of $T_{U\times_\sX U/k}$.
\end{definition}

Using \autoref{lemma:pullback_under_etale} and standard arguments, one sees that the set of $1$-foliations on $\sX$ does not depend on the choice of the presentation $U\to \sX$. So we can intrinsically talk about $1$-foliations on $\sX$ as coherent sub-modules $\sF\subset T_{\sX/k}$ that are closed under Lie brackets and $p$-th powers.

\begin{example}\label{example:foliation_on_qt_stack}
Suppose that $\sX=[\Spec(A)/G]$, where $G$ is a finite \'{e}tale group acting $k$-linearly on $A$. Notice that $G$ acts on $T_{A/k}$ by $\sigma\cdot \partial=\sigma^{-1}\circ \partial \circ \sigma$ for $\sigma\in G$ and $\partial\in T_{A/k}$. Then a $1$-foliation $\sF\subset T_{A/k}$ defines an $1$-foliation on $\sX$ if and only if $\sF$ is a $G$-sub-module of $T_{A/k}$.
\end{example}

Recall that regularity of schemes and local freeness of coherent modules are properties that are local in the \'{e}tale topology. Therefore we can state:

\begin{definition}
Let $\sF$ be a $1$-foliation on $\sX$. The \textbf{singular locus} of $\sF$, written $\Sing(\sF)$, is the locus of $\sX$ where either $\sX$ is not regular, or where $\sX$ is regular but $\sF$ is not a sub-bundle of $T_{\sX/k}$. We say that $\sF$ is \textbf{regular} if $\Sing(\sF)=\emptyset$.
\end{definition}

It is easily seen that $\Sing(\sF)$ is closed in $\sX$, and that if $f\colon U\to \sX$ is an \'{e}tale morphism from a scheme, then $f^{-1}\Sing(\sF)=\Sing(f^*\sF)$. %See also \autoref{lemma:sing_locus_1_foliation} below.

\subsection{Infinitesimal quotients}
Next we explain how to generalize the construction of infinitesimal quotients. The essential observation is given by the following lemma.

\begin{lemma}\label{lemma:infinitesimal_quotient_stacks}
Let $f\colon V\to U$ be an \'{e}tale morphism of normal $k$-schemes of finite type, and $\sF\subset T_{U/k}$ be a $1$-foliation. Then we have a Cartesian diagram
		\begin{equation*}\label{eqn:infinitesimal_quotient_stacks}
		\begin{tikzcd}
		V\arrow[d, "q_V"] \arrow[r, "f"] & U\arrow[d, "q_U"] \\
		V/f^*\sF \arrow[r, "g"] & U/\sF
		\end{tikzcd}
		\end{equation*}
where $q_U,q_V$ are the quotient morphisms, and $g$ is \'{e}tale.
\end{lemma}
\begin{proof}
Since $q_U$ is a universal homeomorphism and $f$ is \'{e}tale, it follows from the topological invariance of the small \'{e}tale site \cite[04DZ]{Stacks_Project} that there exists an \'{e}tale morphism $h\colon W\to U/\sF$ such that $V=W\times_{h,U/\sF,q_U}U$. Since $q$ is a universal homeomorphism of height one, so is the projection $q\colon V\to W$. Since $U/\sF$ is normal and $h$ is \'{e}tale, by the Jacobson correspondence \cite[Theorem 2.5.12]{Posva_Singularities_of_quotients_by_foliations} it follows that $q$ is the quotient by a $1$-foliation $\sG\subset T_{V/k}$. We have to show that $\sG=f^*\sF$: this can be checked at after localizing over the generic points of $U$. So we reduce to the case where: $U=\Spec(L)$ is the spectrum of a field and $\sF\subset \Der_k(L)$, the quotient $U/\sF$ is the spectrum of the field $K=L^\sF$, and $W$ is the spectrum of an \'{e}tale $K$-algebra
		$$A=\bigoplus_{i=1}^n K[x]/(h_i), \quad h_i'\neq 0.$$
Therefore $V=W\times_{U/\sF}U$ is the spectrum of $B=\bigoplus_{i=1}^n L[x]/(h_i)$. Write $K_i=K[x]/(h_i)$ and $L_i=L[x]/(h_i)$. Since $L_i$ is unramified over $L$ we have $\Der_k(L_i)=\Der_k(L)\otimes_L L_i$, and the pullback of $\sF$ is the sub-module $\sF_i=\sF\otimes_L L_i$. To conclude the proof, it suffices to show that $K_i=L_i^{\sF_i}$ for each $i$.

I claim that $K_i\subseteq L_i^{\sF_i}$. To see this, we unravel the action of $\sF_i$. Take $\partial\in \sF$, and regard it as an element of $\Hom_L(\Omega_{L/k},L)$. Then $\partial\otimes 1\in \sF_i$, regarded as a derivation of $L_i$, is the composition:
		$$\begin{tikzcd}
		L_i\arrow[r, "d_{L_i}"] &
		\Omega_{L_i/k} \cong \Omega_{L/k}\otimes_L L_i
		\arrow[r, "\partial\otimes \id"] &
		L\otimes_L L_i
		\cong L_i.
		\end{tikzcd}$$
Let $ax^j\in L$ with $a\in L$. The arrow $L_i\to \Omega_{L/k}\otimes L_i$ sends $ax^j$ to $d_L(a)\otimes x^j$. So if $a\in K$ we have $(\partial\otimes 1)(ax^j)=0$. This shows that $K_i\subseteq L_i^{\sF_i}$.

To conclude, consider the string of equalities
	$$[L_i:K_i]=[L:K]=p^{\rk_L\sF}=p^{\rk_{L_i}\sF_i}=[L_i:L_i^{\sF_i}]$$
where the second one and the last one holds by the Jacobson correspondence, and the first one holds by definition of $L_i$. Since $K_i\subseteq L_i^{\sF_i}$ we deduce that $K_i= L_i^{\sF_i}$ as desired.
\end{proof}

\begin{construction}\label{construction:infinitesimal_qt_of_stacks}
Let $\sX$ be a normal separated Deligne--Mumford stack of finite type over $k$, and $\sF\subset T_{\sX/k}$ be a $1$-foliation. Choose an \'{e}tale presentation $U\to \sX$, and let $V=U\times_\sX U$ with \'{e}tale projections $s,t\colon V\rightrightarrows U$. Since $\sF_U$ descends to $\sX$, we have $s^*\sF= t^*\sF$ as $1$-foliations on $V$; we let $\sF_V$ stand for it. Applying \autoref{lemma:infinitesimal_quotient_stacks} to $V\to U\to U/\sF_U$ and to $V/\sF_V\to U/\sF_U\to U^{(1)}$, we get two commutative diagrams
		$$\begin{tikzcd}
		V\arrow[r, "s"] \arrow[d, "q_V" left] & U\arrow[d, "q_U"] \\
		V/\sF_V\arrow[r, "\bar{s}"] \arrow[d, "q_{V/\sF_V}" left] & U/\sF_U\arrow[d, "q_{U/\sF_U}"] \\
		V^{(1)}\arrow[r, "s^{(1)}"]& U^{(1)}
		\end{tikzcd}
		\quad
		\text{and}
		\quad
		\begin{tikzcd}
		V\arrow[r, "t"] \arrow[d, "q_V" left] & U\arrow[d, "q_U"] \\
		V/\sF_V\arrow[r, "\bar{t}"] \arrow[d, "q_{V/\sF_V}" left] & U/\sF_U\arrow[d, "q_{U/\sF_U}"] \\
		V^{(1)}\arrow[r, "t^{(1)}"]& U^{(1)}
		\end{tikzcd}
		$$
where $\bar{s},\bar{t}$ are \'{e}tale and all squares are cartesian. As $s^{(1)},t^{(1)}\colon V^{(1)}\rightrightarrows U^{(1)}$ is an \'{e}tale groupoid, and since 
		$$V^{(1)}\times_{s^{(1)},U^{(1)}}U/\sF_U\cong V/\sF_V\cong V^{(1)}\times_{t^{(1)},U^{(1)}}U/\sF_U$$
over $V^{(1)}$, we deduce from \autoref{lemma:pullback_groupoid} below that
		\begin{equation}\label{eqn:infinitesimal_qt_groupoid}
		\bar{s},\bar{t}\colon V/\sF_V\rightrightarrows U/\sF_U
		\end{equation}
is an \'{e}tale groupoid.
\end{construction}

\begin{lemma}\label{lemma:pullback_groupoid}
Let $s,t\colon V\rightrightarrows U$ be an \'{e}tale groupoid of schemes. Suppose that $U'\to U$ is a morphism with the property that there exists an isomorphism 
		$$V\times_{s,U}U'\cong V\times_{t,U}U'\quad \text{over }V.$$
Then the base-change $V\times_UU'\rightrightarrows U'$ is an \'{e}tale groupoid.
\end{lemma}
\begin{proof}
The symbol $V'=V\times_UU'$ is any of $V\times_{s,U}U'$ or $V\times_{t,U}U'$: without loss of generality, say it denotes the first one. Denote the second projection by $s'\colon V'\to U'$. The isomorphism $V'\cong V\times_{t,U}U'$ followed by the second projection defines $t'\colon V'\to U'$.

We check the existence of compositions for $(U',V',s',t')$: the other groupoid properties are verified in a similar way. Denote by $h\colon U'\to U$ the morphism given in the statement. 

For a scheme $T$, suppose we have $\alpha,\beta,\gamma\in U'(T)$ and $f,g\in V'(T)$ such that
		$$s'(f)=\alpha, \quad 
		t'(f)=\beta=s'(g), \quad
		t'(g)=\gamma.$$
We suggestively denote this by $f\colon \alpha\to\beta$, etc. Under $V'\cong V\times_{t,U}U'$ we have $f=(f^\flat\colon h(\alpha)\to h(\beta), \beta)$, and under $V'=V\times_{s,U}U'$ we have $g=(g^\flat\colon h(\beta)\to h(\gamma), \beta)$. The composition $g^\flat\circ f^\flat\colon h(\alpha)\to h(\gamma)$ exists as $V\rightrightarrows U$ is a groupoid, and we define the desired composition $g\circ f$ as the pair $(g^\flat\circ f^\flat, \alpha)$ in $V(T)\times_{s,U(T)}U'(T)=V'(T)$.
\end{proof}

We check that the output of our construction does not depend on the choice of presentation $U\to \sX$.

\begin{lemma}
The quotient of \autoref{eqn:infinitesimal_qt_groupoid} does not depend on the choice of the \'{e}tale presentation $U\to \sX$.
\end{lemma}
\begin{proof}
Given another presentation $U'\to \sX$, we can form the fiber product $\sU=U\times_{\sX}U'$. Then $\sU$ is \'{e}tale surjective over both $U$ and $U'$, and $\sU\to \sX$ is a presentation. Hence it suffices to consider the case where $U'\to \sX$ factors as 
		$$U'\overset{h}{\longrightarrow} U\longrightarrow\sX$$
where $h$ is \'{e}tale surjective. Let $V=U\times_\sX U$ and $V'=U'\times_{\sX}U'$. We have a morphism of groupoids $f_\bullet=(f_1,f_0)\colon (V'\rightrightarrows U')\to (V\rightrightarrows U)$ where $f_0\colon U'\to U$ is $h$ and $f_1\colon V'\to V$ is given by the left vertical arrow is the cartesian square
		$$\begin{tikzcd}
		U'\times_\sX U' \arrow[r]\arrow[d] & U'\times U' \arrow[d, "h\times h"] \\
		U\times_\sX U \arrow[r] & U\times U.
		\end{tikzcd}$$
In particular, notice that $f_1$ is \'{e}tale. The morphism of groupoids $f_\bullet$ descends to a morphism of associated quotient stacks $\mathfrak{f}\colon [V'\rightrightarrows U']\to [V\rightrightarrows U]$ which is an isomorphism by \cite[04T5]{Stacks_Project}.

To compare the quotients of these groupoids by the foliation $\sF$, it is more convenient to use a bottom-up approach and start from the Frobenius twists $V^{(1)}\rightrightarrows U^{(1)}$, etc. By \autoref{lemma:infinitesimal_quotient_stacks} and \autoref{construction:infinitesimal_qt_of_stacks}, we have a commutative diagram
	$$\begin{tikzcd}
	& V'/\sF_{V'} \arrow[rrr, shift left]\arrow[rrr, shift right] \arrow[dd] \arrow[dl, "\bar{f}_1" above left]
	 &&& U'/\sF_{U'}\arrow[dd]\arrow[dl, "\bar{f}_0"] \\
	V/\sF_V \arrow[crossing over, rrr, shift left]\arrow[crossing over, rrr, shift right] \arrow[dd]
	 &&& U/\sF_U &
	\\
	& V'^{(1)} \arrow[dl, "f^{(1)}_1" above left] \arrow[rrr, shift left]\arrow[rrr, shift right]
	&&& U'^{(1)}\arrow[dl, "f^{(1)}_0"] \\
	V^{(1)}\arrow[rrr, shift left]\arrow[rrr, shift right]
	&&& U^{(1)}\arrow[from=uu, crossing over] &
	\end{tikzcd}$$	
where every vertical face is cartesian. So we get a morphism of groupoids 
	$$\bar{f}_\bullet=(\bar{f}_1,\bar{f}_0)\colon
		(V'/\sF_{V'}\rightrightarrows U'/\sF_{U'}) \to
		(V/\sF_V \rightrightarrows U/\sF_U)$$
which descends to a morphism
	$$\bar{\mathfrak{f}}\colon [V'/\sF_{V'}\rightrightarrows U'/\sF_{U'}]\to
		[V/\sF_V \rightrightarrows U/\sF_U].$$
By construction $\bar{\mathfrak{f}}$ is the pullback of $\mathfrak{f}^{(1)}$, and so it is an isomorphism. This completes the proof.
\end{proof}

%We used the following formal property of cartesian squares in the previous proof.

%\begin{lemma}\label{lemma:cartesian_squares}
%Consider a commutative diagram
%		$$\begin{tikzcd}
%		C\arrow[r] \arrow[d] & B \arrow[d] \\
%		C'\arrow[r]\arrow[d] & B'\arrow[d] \\
%		A\arrow[r] & S
%		\end{tikzcd}$$
%where the bottom and the big squares are cartesian. Then the upper square is cartesian.
%\end{lemma}
%\begin{proof}
%Consider the diagram of solid arrows
%		$$\begin{tikzcd}
%		T\arrow[drr, bend left, "\psi"]\arrow[ddr, bend right, "\varphi"]\arrow[dr, dotted, "\eta"] && \\
%		&C\arrow[r, "u"] \arrow[d, "m"] & B \arrow[d, "g"] \\
%		&C'\arrow[r, "s"]\arrow[d, "h"] & B'\arrow[d, "f"] \\
%		&A\arrow[r, "t"] & S
%		\end{tikzcd}$$
%where $s\varphi=g\psi$, and we wish to prove the existence of an arrow $\eta$ that is unique with the property that $m\eta=\varphi$ and $u\eta=\psi$. We have $fg\psi=fs\varphi=th\varphi$ by commutativity of the bottom square, so as the big square is cartesian we obtain an arrow $\eta\colon T\to C$ that is unique with the property that $hm\eta=h\varphi$ and $u\eta=\psi$. So it remains to show that $m\eta=\varphi$. As the bottom square is cartesian, it suffices to show that $sm\eta=g\psi=s\varphi$. The second equality holds by assumption, while the first one follows from $g\psi=gu\eta=sm\eta$ (using the definition of $\psi$ and the commutativity of the upper square). So we are done.
%\end{proof}

\begin{lemma}
The quotient stack defined by the groupoid \autoref{eqn:infinitesimal_qt_groupoid} is a separated Deligne--Mumford stack of finite type over $k$, with diagonal representable by schemes.
\end{lemma}
\begin{proof}
The quotient stack, which we temporarily denote by $\sQ$, is a Deligne--Mumford stack by \cite[Proposition 4.3.1]{Laumon_Moret_Bailly_Champs_algebriques}. Since $U/\sF_U$ is of finite type over $k$, so is $\sQ$ (see \cite[036O]{Stacks_Project}). To study its diagonal, consider the 2-cartesian diagram 
		$$\begin{tikzcd}
		V/\sF_V\arrow[rr, "\bar{s}\times\bar{t}"]\arrow[d] && U/\sF_U\times_k U/\sF_U\arrow[d] \\
		\sQ\arrow[rr,"\Delta"] && \sQ\times_k\sQ
		\end{tikzcd}$$
where the vertical arrows are \'{e}tale (cf.\ the diagram in the proof of \cite[Proposition 4.3.2]{Laumon_Moret_Bailly_Champs_algebriques}). Then $\Delta$ is representable by schemes because $V/\sF_V$ and $U/\sF_U$ are schemes. The morphism $\bar{s}\times\bar{t}$ is proper because it is the base-change of the morphism $s^{(1)}\times t^{(1)}$ which is proper by assumption. As properness is local on the target for the \'{e}tale topology  (see \cite[02YJ]{Stacks_Project}), we deduce that $\Delta$ is proper, which by definition means that $\sQ$ is separated.
\end{proof}

\begin{definition}\label{def:infinitesimal_qt_stack}
Let $\sF$ be a $1$-foliation on a normal separated Deligne--Mumford stack $\sX$ of finite type over $k$. The quotient $\sX/\sF$ is the normal separated Deligne--Mumford stack defined as the quotient of the groupoid \autoref{eqn:infinitesimal_qt_groupoid}. It comes with a finite purely inseparable morphism $q\colon \sX\to \sX/\sF$ of height one.
\end{definition}

The essential example to keep in mind is the following one.
\begin{example}[Infinitesimal quotients of affine quotient stacks]\label{example:local_structure_of_infinitesimal_qt}
Consider $\sX=[\Spec(A)/G]$, where $G$ is a finite \'{e}tale group over $k$. Let $\sF\subset T_{\sX/k}$ be a $1$-foliation, corresponding to the $G$-$A$-module $\sF_A\subset T_{A/k}$. Then I claim that
		$$\sX/\sF=\left[\Spec\left(A^{\sF_A}\right)/G\right],
		\quad \text{with coarse moduli }
		\Spec\left( \left( A^{\sF_A}\right)^G\right)=\Spec\left( \left( A^G\right)^{\sF_A}\right).$$
Indeed, in this case the groupoid scheme defining $\sX$ is given by $a,p_2\colon G\times\Spec(A)\rightrightarrows \Spec(A)$ where $a$ is the action and $p_2$ is the second projection. In particular $p_2^*\sF_A=a^*\sF_A$ by \autoref{example:foliation_on_qt_stack}. So in this case, \autoref{construction:infinitesimal_qt_of_stacks} produces the cartesian diagram
		$$\begin{tikzcd}
		G\times \Spec(A)\arrow[r, "a\times p_2"]\arrow[d, "\id\times q"] & \Spec(A)\times\Spec(A) \arrow[d,"q\times q"] \\
		G\times \Spec(A^{\sF_A})\arrow[r, "u"] & \Spec(A^{\sF_A})\times \Spec(A^{\sF_A})
		\end{tikzcd}$$
where $q\colon \Spec(A)\to \Spec(A^{\sF_A})$ is the quotient morphism. It is easy to observe that the $G$-action on $A$ restricts to a $G$-action on $A^{\sF_A}$, and the map $u^*\colon A^{\sF_A}\otimes A^{\sF_A}\to \prod_G A^{\sF_A}$ is simply the (co-)restriction of $(a\times p_2)^*$. In particular, we obtain that $\sX/\sF=\left[\Spec\left(A^{\sF_A}\right)/G\right]$. Its coarse moduli space is the spectrum of the sub-ring of $G$-invariants of $A^{\sF_A}$, which also can be described as the  sub-ring of $\sF_A$-constants of $A^G$, for both rings are equal to the intersection $A^{\sF_A}\cap A^G$.
\end{example}

\subsection{Properties of infinitesimal quotients and the Jacobson correspondence}
Next we generalize \cite[Lemma 2.5.10]{Posva_Singularities_of_quotients_by_foliations}.

\begin{proposition}\label{lemma:regularity_of_qt}
Assume that $\sX$ is regular, and let $\sF\subset T_{\sX/k}$ be a $1$-foliation. Then $\sF$ is regular if and only if $\sX/\sF$ is regular. 
\end{proposition}
\begin{proof}
If $\sX=[V\rightrightarrows U]$ as above, then $U$ and $V$ are regular. Thus by \cite[Lemma 2.5.10]{Posva_Singularities_of_quotients_by_foliations}, the quotients $U/\sF_U$ and $V/\sF_V$ are regular if and only if $\sF_U$ and $\sF_V$ are regular. Since the two morphisms $V\rightrightarrows U$ are \'{e}tale, one sees that $\sF_V$ is regular if and only if $\sF_U$ is regular. Therefore we obtain: $X/\sF$ is regular if and only if $U/\sF_U$ and $V/\sF_V$ are regular, if and only if $\sF_U$ is regular, if and only if $\sF$ is regular.
\end{proof}

\begin{lemma}\label{lemma:qt_preserves_tameness}
Suppose that $\sX$ is tame, and let $\sF\subset T_{\sX/k}$ be a $1$-foliation. Then $\sX/\sF$ is tame.
\end{lemma}
\begin{proof}
Let $c_\sX\colon \sX\to X$ and $c_{\sY}\colon \sY=\sX/\sF\to Y$ be the (Keel--Mori) coarse moduli morphisms. We have a commutative diagram
		$$\begin{tikzcd}
		\sX\arrow[d, "c_\sX"]\arrow[r, "q"] & \sY\arrow[d, "c_{\sY}"]\\
		X\arrow[r, "q_c"] & Y.
		\end{tikzcd}$$
Since $q,c_\sX,c_{\sY}$ induce bijections on geometric points, so does $q_c$. So $q_c$ is universally injective and surjective. Since $q_*\sO_{\sX}$ is a finite $\sO_{\sY}$-module, we see that $c_{\sX,*}\sO_{\sX}=\sO_X$ is a finite $c_{\sY,*}\sO_{\sY}=\sO_Y$-module. Therefore $q_c$ is a universal homeomorphism \cite[04DF]{Stacks_Project}. 

Now it is well-known that we can find an \'{e}tale covering $\{X_i\to X\}$ inducing cartesian diagrams
		$$\begin{tikzcd}
		\left[\Spec(A_i)/G_i\right] \arrow[d]\arrow[r] & \sX\arrow[d] \\
		\Spec\left( A_i^{G_i}\right)=X_i\arrow[r] & X
		\end{tikzcd}$$
where the $G_i$'s are finite \'{e}tale $k$-groups whose orders are invertible in $k$ (see e.g.\ \cite[Theorem 4.3.1]{Alper_Stacks_and_moduli}). By \cite[04DZ]{Stacks_Project} we may assume that $X_i=Y_i\times_Y Y$ for some \'{e}tale morphism $Y_i\to Y$. Then $\{Y_i\to Y\}$ is an \'{e}tale covering, and by \cite[Theorem 3.2]{Abramovich_Olsson_Vistoli_Tame_stacks} it is sufficient to check that $\sY\times_YY_i$ is tame. Thus we reduce to the case where $\sX=[\Spec(A)/G]$ with $G$ finite discrete of order invertible in $k$. Then $\sX/\sF=[\Spec(A^{\sF_A})/G]$ by \autoref{example:local_structure_of_infinitesimal_qt}, and this quotient stack is tame by the assumptions on $G$.
\end{proof}

To conclude, let us extend the Jacobson correspondence (see \cite[Theorem 2.5.12]{Posva_Singularities_of_quotients_by_foliations}) to the stack-theoretic setting. We require the following lemma.

\begin{lemma}\label{lemma:descent_of_model_foliations}
Let $f\colon U\to U'$ be a purely inseparable morphism of height one, where $U$ and $U'$ are normal connected schemes of finite type over $k$. Let $s',t'\colon V'\rightrightarrows U'$ be an \'{e}tale groupoid scheme and let $s,t\colon V\rightrightarrows U$ be its pullback along $f$. Then $\Der_{f^{-1}\sO_{U'}}(\sO_U)$ descends to a sub-module of $T_{[V\rightrightarrows U]/k}$.
\end{lemma}
%Consider a cartesian square of \'{e}tale groupoid schemes
%		$$\begin{tikzcd}
%		V\arrow[r, shift right, "t" below]\arrow[r, shift left, "s" above]\arrow[d] & U \arrow[d, "f"] \\
%		V'\arrow[r, shift right, "t'" below]\arrow[r, shift left, "s'" above] & U'
%		\end{tikzcd}$$
%where $U,U'$ are normal connected $k$-schemes of finite type, and $f$ factors the Frobenius $F_{U/k}$. 
\begin{proof}
Write $G=\Der_{f^{-1}\sO_{U'}}(\sO_U)=T_{U/U'}\subset T_{U/k}$. We have to check that the natural descent datum of $T_{U/k}$ for $V\rightrightarrows U$, namely
		$$\begin{tikzcd}
		(\partial t)\circ (\partial s)^{-1}\colon 	
		s^*T_{U/k} \arrow[r, "\sim"] &
		T_{U\times_\sX U/k} \arrow[r, "\sim"]
		& t^*T_{U/k}
		\end{tikzcd},$$
restricts to a descent datum for $G$. Since $(\partial s)^{-1}(s^*G)=T_{V/V'}=(\partial t)^{-1}(t^*G)$, this is indeed the case.
\end{proof}

\begin{theorem}[Jacobson's correspondence]\label{thm:stacky_Jacobson}
Let $\sX$ be a Deligne--Mumford (resp. tame Deligne--Mumford) separated normal connected stack of finite type over $k$, with diagonal representable by schemes. Then there is a bijection between
	\begin{enumerate}
		\item $1$-foliations of corank $r$ on $\sX$, and
		\item factorizations $F_{\sX/k}\colon \sX\to \sY\overset{g}{\to} \sX^{(1)}$, where $\sY$ is a Deligne--Mumford (resp. tame Deligne--Mumford) separated normal connected stack of finite type over $k$ with diagonal representable by schemes, $g$ is representable by schemes and has degree $p^r$.
	\end{enumerate}
The bijection is explicitly given by $\sF\mapsto \sX/\sF$ and $(\sX\to\sY)\mapsto \Der_{\sO_\sY}(\sO_\sX)$. 
\end{theorem}

\begin{remark}
With the notations of \autoref{thm:stacky_Jacobson}, if $\sY\to \sX^{(1)}$ is representable by schemes then $\sX\to \sY$ is automatically representable by schemes (notice that $F_{\sX/k}$ is representable and apply \cite[Lemme 3.12]{Laumon_Moret_Bailly_Champs_algebriques}).
\end{remark}

\begin{proof}
Choose an \'{e}tale presentation $h^{(1)}\colon U^{(1)}\to \sX^{(1)}$. Then for every representable $\sY\to \sX^{(-1)}$, the pullback morphism $U^{(1)}\times_{\sX^{(1)}}\sY\to\sY$ is an \'{e}tale presentation. In particular, the second set in the statement is in bijection with the set of commutative diagrams
		$$\begin{tikzcd}
		U\arrow[r, "h"] \arrow[d] & \sX\arrow[d] \\
		W\arrow[r]\arrow[d] & \sY \arrow[d, "g"] \\
		U^{(1)}\arrow[r, "h^{(1)}"] & \sX^{(1)}
		\end{tikzcd}$$
where all squares are cartesian, all horizontal arrows are \'{e}tale and all vertical arrows are purely inseparable. By the usual Jacobson correspondence and \autoref{lemma:descent_of_model_foliations}, this set is in bijection with the set of $1$-foliations $\sF_U\subset T_{U/k}$ of corank $r$ with a descent datum for $U\times_\sX U\rightrightarrows U$ that is compatible with the natural one of $T_{U/k}$. By definition, the latter set is in bijection with the $1$-foliations of corank $r$ on $\sX$. The specialization to tame stacks is ensured by \autoref{lemma:qt_preserves_tameness} (notice that if $\sY$ is given, then $\sX^{(1)}$ is the quotient of $\sY$ by a $1$-foliation), and the description of the bijection is easily verified.
\end{proof}

\section{Resolution on surfaces}\label{section:resolution_surfaces}
In this section we prove our resolution results for $1$-foliations on surfaces. 

\begin{theorem}\label{thm:surface_foliations_resolution}
Let $S$ be a normal surface over an algebraically closed field $k$ of characteristic $p>0$, and let $\sF$ be a $1$-foliation of rank one on $S$. Then there exists a projective birational morphism $f\colon S'\to S$ such that:
	\begin{enumerate}
		\item $S'$ is regular,
		\item $f$ is an isomorphism over the locus where $\sF$ has at worst multiplicative singularities \cite[Definition 2.4.7]{Posva_Singularities_of_quotients_by_foliations}, and
		\item the induced $1$-foliation $f^*\sF$ has at worst multiplicative singularities.
	\end{enumerate}
\end{theorem}
Let us recall that $\sF$ has multiplicative singularities at $s\in \Sing(\sF)$ if there are formal coordinates $x,y$ at $s$ such that $\sF\otimes \widehat{\sO}_{s}$ is generated by $x\partial_x+\lambda y\partial_y$ for some $\lambda\in\bF_p^\times$.
\begin{proof}
By resolving the singularities of $S$ (which by definition are contained in the singular locus of $\sF$), we may assume that $S$ is regular. As $\Sing(\sF)$ is a finite set the result is local on $S$, so we may also assume that $S$ is affine, say $S=\Spec(A)$. The function field $K(S/\sF)$ is purely inseparable of degree $p$ over $K(S^{(1)})$, so we can write $K(S/\sF)=K(S^{(1)})[u^{1/p}]$ for some $u\in \Gamma(S^{(-1)},\sO_{S^{-1)}})$. 

According to a theorem of Giraud \cite{Giraud_Forme_normale_sur_une_surface_en_car_pos}, there is a regular proper model $f\colon V\to S^{(1)}$, depending on $u$ and given by a sequence of smooth blow-ups, such that every point $v\in V$ has a finite \'{e}tale neighbourhood $(v\in \fV)$ on which the function $f^*u\in \sO_{V}(V)$ can be written 
		\begin{equation}\label{eqn:etale_neighbourhoods}
		(f^*u)|_\fV=t^p+\bold{x}^\bold{a}
		\end{equation}
where $t\in \sO_\fV(\fV)$, $\bold{x}=(x_1,\dots,x_n)$ is a subset of a $p$-basis \cite[\S 2.2]{Posva_Singularities_of_quotients_by_foliations} of $\sO_{\fV,v}$ over $\sO_{\fV^{(1)},v}$, so $1\leq n\leq 2$, and $\bold{a}\neq \bold{0}$ modulo $p$ (\footnote{
		Here and latter, we write $\bold{x}^\bold{a}$ to denote the product $\prod_{i=1}^n x_i^{a_i}$.
}). Since $\sO_{\fV,v}$ has a perfect residue field, we can in fact assume that $\bold{x}$ is a partial system of local parameters (see \cite[Paragraph before Remarque 1.6]{Giraud_Forme_normale_sur_une_surface_en_car_pos}). The morphism $f$ induces a commutative diagram
	$$\begin{tikzcd}
	V^{(-1)}\arrow[r]\arrow[d, "f^{(-1)}"] & V\left[(f^*u)^{1/p}\right]^\nu\arrow[d]\arrow[r] & V\arrow[d, "f"] \\
	S\arrow[r] & S/\sF\arrow[r] & S^{(1)}
	\end{tikzcd}$$
where the vertical arrows are proper and birational. Observe that 
	$$V^{(-1)}/f^{(-1),*}\sF=V\left[(f^*u)^{1/p}\right]^\nu.$$ Indeed, both varieties are normal and finite over $V$, with function field $K(S/\sF)$; therefore both of them are given by the normalization of $V$ insider $K(S/\sF)$.

We claim that $f^{(-1)}\colon V^{(-1)}\to S$ is the morphism we are looking for. The fact that $f^{(-1)}$ can be taken to be an isomorphism over $\Reg(\sF)$ will be shown at the very end. In the meantime, to simplify the notations, we may replace $S$ by $V^{(-1)}$, etc., and assume that \'{e}tale neighbourhoods as in \autoref{eqn:etale_neighbourhoods} already exist over $S^{(1)}$.

We are going to study the normalization of $B=\widehat{\sO}_{S^{(1)},s}[T]/(T^p-u)$, when \emph{$s$ is a closed point}. By \autoref{lemma:quotient_is_normalization}, this normalization is the completed local ring of $S/\sF$ at $s$. Notice that $B$ is local and complete \cite[07N9]{Stacks_Project}, hence Henselian. Let us write $\widehat{\sO}_{S^{(1)},s}=\sO_s$ with maximal ideal $\fm_s$. As $\fV\to S^{(1)}$ is \'{e}tale it induces an isomorphism on completions, so the above decomposition of $u|_\fV$ descends to $\sO_{s}$. Replacing $T$ by $T-t$, we may assume that $t=0$. We distinguish several cases.
	\begin{enumerate}
		\item $\bold{x}=(x)$ and $\bold{a}=(a)$ where $(a,p)=1$. By \autoref{prop:Euclid_algorithm} there exists $w$ in the normalization of $B$ such that $w^p=x$ and $w^a=T$. Thus we have $B[w]\subseteq B^\nu$ with $B[w]\cong \sO_{s}[W]/(W^p-x)$. 
		If $y\in \sO_s$ is such that $\fm_s=(x,y)$ then we obtain
				$$B^\nu=B[w]= k\llbracket x^{1/p},y\rrbracket.$$
		Therefore $B^\nu$ is regular.
		\item $\bold{x}=(x_1,x_2)$ and $\bold{a}=(a_1,pa_2)$. The element $w=T/x_2^{a_2}$ satisfies the monic equation $w^p=x_1^{a_1}$, and so $w$ belongs to $B^\nu$. This gives the partial normalization
				$$B[w]\cong \sO_s[W]/(W^p-x_1^{a_1}).$$
		Now case (a) above applies, and we find that $B^\nu$ is regular.
		\item $\bold{x}=(x_1,x_2)$ and $\bold{a}=(a_1,a_2)$ where $(a_i,p)=1$ for $i=1,2$. Then $(x_1,x_2)=\fm_s$ and:
				\begin{equation}\label{eqn:only_singularities}
				\begin{array}{l}
				\text{there is a system of local parameters }x,y \text{ for }\sO_s\text{ such that} \\
				B^\nu \cong \left(
				k\llbracket x,y\rrbracket [T]/(T^p-x^ay^b)\right)^\nu \text{ for }a,b\text{ coprime with }p.
				\end{array}
				\end{equation}
	\end{enumerate}
So we have obtained that: the surface $S/\sF$ is regular outside finitely many closed points, the completed stalks of which can be described as in \autoref{eqn:only_singularities}. 

We now reformulate our singularity statements into properties of $\sF$. Let $s\in S$ be a closed point. If $s$ is a regular point of $S/\sF$, we apply \cite[Lemma 2.5.10]{Posva_Singularities_of_quotients_by_foliations} to get that $\sF$ is regular at $s\in S$. Suppose that $s\in S/\sF$ is not regular, so it is given with its $\widehat{\sO}_{S^{(1)},s}$-algebra structure as in \autoref{eqn:only_singularities}. Then applying \autoref{prop:inv_subring_normaliz_hypersurface} to $\sO_{S^{(1)},s}\hookrightarrow \sO_{S/\sF,s}\hookrightarrow \sO_{S,s}$ we see that there are formal local parameters $u,v$ of $\widehat{\sO}_{S,s}$ and $\alpha,\beta\in\bF_p^\times$ such that, if we write
		$$D_s=\alpha u\frac{\partial}{\partial u}+\beta v\frac{\partial}{\partial v}\in\Der_k^\text{cont}\left(\widehat{\sO}_{S,s}\right),$$
then in $\widehat{\sO}_{S,s}$ the subrings of $\sF\otimes\widehat{\sO}_{S,s}$-constants and of $D_s$-constants coincide. Since $D_s$ is $p$-closed, by the Jacobson correspondence we obtain that $D_s$ generates formally $\sF\otimes\widehat{\sO}_{S,s}$ at $s$. Since $D_s$ is multiplicative, we have obtained that: at every closed point of $S$ the foliation $\sF$ is either regular or generated formally by a multiplicative derivation. Non-closed points of $S$ have codimension at most $1$, so $\sF$ is necessarily regular there.

Returning to our original non-simplified notations, it remains to prove that $f^{(-1)}\colon V^{(-1)}\to S$ (where $S$ is regular) can be taken to be an isomorphism over the locus where $\sF$ has at worst multiplicative singularities. The above calculations show that
	\begin{equation}\label{eqn:u_good_then_F_good}
	\begin{matrix}
	u\in \sO_{S^{(1)},s}\text{ is \'{e}tale-locally} \\
	\text{of the form }\autoref{eqn:etale_neighbourhoods}
	\end{matrix}		
	\quad \Longrightarrow \quad
	\begin{matrix}
	\sF\text{ has at worst multiplicative}\\
	\text{singularities at }s\in S.
	\end{matrix}
	\end{equation}
It follows that we can refine our application of Giraud's result as follows. We can decompose $f^{(-1)}$ as
		$$f^{(-1)}=\pi_1\circ\cdots\circ\pi_N\colon 
		V^{(-1)}=S_N\to \dots \to S_0=S$$
where each $\pi_{i}$ a sequence of blow-ups of infinitely closed points, and if $b_i\in S$ is the image of the $\pi_i$-exceptional divisors then $b_i\neq b_j$ for $i\neq j$. By \autoref{eqn:u_good_then_F_good} we see that the locus of non-multiplicative singularities of $\sF$ is contained in $\{b_i\}_{i=1}^N$. If $\sF$ happens to have multiplicative singularities at $b_i$ (which can happen, see \autoref{rmk:u_vs_F} below), then we suppress the blow-ups of $\pi_i$ from the definition of $f^{(-1)}$ (which is possible by our definition of $\pi_i$). Thus we may assume that the locus of non-multiplicative singularities of $\sF$ is precisely $\{b_i\}$, and the theorem is proved.
\end{proof}

\begin{remark}\label{rmk:u_vs_F}
Let us use the notations of the proof of \autoref{thm:surface_foliations_resolution}. Then the morphism $f$ is very sensitive to the choice of $u\in \sO_{S^{(1)}}$ inducing $K(S)=K(S^{(1)})(u^{1/p})$, while the singularities of $S/\sF$ are the those of the \emph{normalized} $p$-cyclic cover induced by $u$. In particular, the converse of \autoref{eqn:u_good_then_F_good} is not true. Consider for example the case $S=\bA^2_{x,y}$ and $\sF$ generated by $\partial_x$. If $w=w(x^{p},y)\in k[x,y]$ is not a $p$-th power, observe that
		$$\sF=\{\partial\in T_{S/k}\mid \partial(w)=0\}.$$
By Jacobson's correspondence, this implies that $K(S/F)=K(S^{(1)})[w]$ and thus $S/F=S^{(1)}[w]^\nu$. However the element $w^p\in H^0(S^{(1)},\sO_{S^{(1)}})$ usually does not possess the property \autoref{eqn:etale_neighbourhoods}.
\end{remark}

%\begin{remark}
%Given a normal surface $S$ and a $1$-foliation $\sF$ on it, then amongst the \emph{regular} models $S'\to S$ satisfying the conclusions of \autoref{thm:surface_foliations_resolution} there is a minimal one. 
%\end{remark}

Next we use weighted blow-ups to completely resolve the foliation singularities. 

\begin{theorem}\label{thm:stacky_resolution_surfaces}
Let $S$ be a normal surface over an algebraically closed field $k$ of characteristic $p>0$, and let $\sF$ be a $1$-foliation on $S$. Then there exists a proper birational morphism $f\colon \sS'\to S$ such that 
	\begin{enumerate}
		\item $\sS'$ is a regular tame Deligne--Mumford stack on which the induced $1$-foliation $f^*\sF$ is everywhere regular, and
		\item $f$ is an isomorphism over $\Reg(\sF)$.
	\end{enumerate}
Moreover, if $p=2$ we can take $\sS'$ to be a scheme and $f\colon \sS'\to S$ to be projective.
\end{theorem}
\begin{proof}
By \autoref{thm:surface_foliations_resolution}, we can blow-up $S$ along $\Sing(\sF)$ until $S$ is regular and if $\sF$ is not regular at a point $s\in S$, then we can find formal coordinates $x,y$ of $\widehat{\sO}_{S,s}$ such that $\sF\otimes\widehat{\sO}_{S,s}$ is generated by $\partial=x\partial_x+\lambda y\partial_y$ with $\lambda\in \bF_p^\times$.

Let us begin with the case $p=2$. Then $\lambda=1$ is the only possibility at any singular point. Let $S'\to S$ be the blow-up of the finite reduced closed set $\Sing(\sF)$: the computations of \autoref{example:smth_blowup_of_foliations} show that the induced $1$-foliation $\sF'$ on $S'$ is regular. 

Now let us come back to the general case $p>0$. Let $\Lambda\in \{1,\dots,p-1\}$ be such that its image modulo $p$ is $\lambda$. Consider the Rees algebra $\fI_\bullet=(x,1)+(y,\Lambda)$ on $\widehat{\sO}_{S,s}$ (cf.\ \autoref{eqn:weighted_ideals}). As $(x,y)=\widehat{\fm}_{S,s}$ there exists an increasing function $\phi\colon \bZ_+\to\bZ_+$ such that $\fI_n\supset \fm_{S,s}^{\phi(n)}$. This gives a sequence of closed immersions
		$$\Spec(\widehat{\sO}_{S,s}/\fI_n)\hookrightarrow
		\Spec\left(\widehat{\sO}_{S,s}/\fm_{S,s}^{\phi(n)}\right)
		\hookrightarrow S.$$
Thus there is a sheaf of ideals $I_n$ on $S$ that is co-supported at $s$, and such that $I_n\widehat{\sO}_{S,s}=\fI_n$. Then $I_\bullet=\bigoplus_n I_nt^n\subset \sO_S[t]$ is a Rees algebra on $\sO_S$: the defining properties are trivially verified away from $s$, and at $s$ they can be checked on the completion, where they hold by definition.

Let $f_s\colon \sS=\Bl_{I_\bullet}S\to S$ be the blow-up of $I_\bullet$. By compatibility with flat base-change, we have $\sS\otimes \widehat{\sO}_{S,s}\cong \Bl_{\fI_\bullet}\widehat{\sO}_{S,s}$. By fppf descent, we get that $\sS$ is regular tame Deligne--Mumford in a neighbourhood of $f_s^{-1}(s)$. Moreover, $f_s$ is an isomorphism above the complement of $s$.

\begin{claim}\label{claim:weighted_blow_up_surface}
$f_s^*\sF$ is regular in a neighbourhood of $f^{-1}_s(s)$.
\end{claim}
\begin{proof}\renewcommand{\qedsymbol}{$\lozenge$}
By fppf descent, it suffices to check this on $\Bl_{\fI_\bullet}\widehat{\sO}_{S,s}$. By fppf descent again, we may assume that $\sO_{S}=k[x,y]$ and $\fm_s=(x,y)$. Using \autoref{example:weighted_blow_ups}, one computes that on the $D_+(x)$-chart $\bA^2_{u,v}$ it holds that
		$$\partial= u\partial_u +(\lambda -\lambda)v\partial_v=u\partial_u.$$
So on $D_+(x)$, the saturated pullback of $\sF$ is generated by $\partial_u$; in particular, it is regular. On the $D_+(y)$-chart $\bA^2_{u',v'}$, one finds 
		$$\partial=(1-\lambda \lambda^{-1})u'\partial_{u'}+\lambda \lambda^{-1}v'\partial_{v'}=v'\partial_{v'}.$$
So on $D_+(y)$, the saturated pullback of $\sF$ is generated by $\partial_{v'}$: hence it is regular.
\end{proof}

So we see that after a weighted blow-up, we can resolve the singularities of $\sF$ at $s$. We can repeat this procedure at every other point of the finite discrete closed set $\Sing(\sF)$. In the end we obtain a proper morphism $f\colon \sS'\to S$ which is an isomorphism above the complement of $\Sing(\sF)$, and such that $\sS'$ is a regular tame Deligne--Mumford stack on which $f^*\sF$ is regular. This completes the proof.
\end{proof}

\subsection{Functorial resolution for two-dimensional foliated Deligne--Mumford stacks}
In this sub-section, we extend the methods of \autoref{thm:surface_foliations_resolution} and \autoref{thm:stacky_resolution_surfaces} to obtain a functorial resolution theorem for $1$-foliations on two-dimensional Deligne--Mumford stacks. We let $k$ denote an algebraically closed field of characteristic $p>0$.

At first glance, it seems that one can obtain a resolution in three steps: first by resolving the underlying stack, then by blowing-up the non-multiplicative singularities of the $1$-foliation, and finally by performing a weighted blow-up along the multiplicative singularities. However, the last step is not completely functorial, as shown by the following example.

\begin{example}[Non-uniqueness of weighted blow-ups]
Consider the $1$-foliation $\sF$ on $S=\Spec(k\llbracket x,y\rrbracket)$ generated by $x\partial+\lambda y\partial_y$ where $\lambda\in \bF_p$. Another generator of $\sF$ is given by $\lambda^{-1}x\partial_x+y\partial_y$. So according to the computations of \autoref{claim:weighted_blow_up_surface}, we can resolve $\sF$ by blowing-up either $\fI=(x,1)+(y,\Lambda)$ or $\fI'=(x,\Lambda')+(y,1)$ where $\Lambda,\Lambda'\in \{1,\dots,p-1\}$ respectively lift $\lambda$ and $\lambda^{-1}$. But in general $\fI\neq \fI'$, and consequently $\Bl_\fI S\not\cong \Bl_{\fI'} S$. 

For example, suppose that $p=5$ and that $\lambda=3$. Then $\Lambda=3$ and $\Lambda'=2$. The element $x^{3-3}y^{\lceil 3/\Lambda\rceil}=y$ belongs to $\fI_3$. However it does not belong to $\fI'_3=(x^3,x^2y,xy,y^2)$. So $\fI\not\subseteq \fI'$.
\end{example}

In this example, the issue is that the two different generators of $\sF$ induce different endomorphisms of the cotangent space $\fm/\fm^2$. In fact, \autoref{prop:Rees_algebra_is_functorial} essentially shows that if two generators of $\sF$ induces the same diagonalizable endomorphism of $\fm/\fm^2$, with eigenvalues in $\bF_p^\times$, then they give rise to the same Rees algebras. Understanding this point is the key in turning the methods introduced above into a functorial resolution process for $1$-foliations on two-dimensional Deligne--Mumford stacks. 

\medskip
Let us now prove our enhancement of \autoref{thm:stacky_resolution_surfaces}. We need some preliminary definitions and results.

\begin{definition}
Let $\sF$ be a $1$-foliation on a normal variety $X$. We let $\NMSing(\sF)$ be the set of those $x\in \Sing(\sF)$ such that $\sF$ does not have multiplicative singularities at $x$. 
\end{definition}

\begin{lemma}\label{lemma:pullback_non_lc_locus}
Let $f\colon V\to U$ be an \'{e}tale morphism of normal varieties over $k$, and $\sF$ be a $1$-foliation on $U$. For a closed point $u\in U$, we have $u\in \NMSing(\sF)$ if and only if $f^{-1}(u)\subset \NMSing(f^*\sF)$.
\end{lemma}
\begin{proof}
Let $v\in f^{-1}(u)$ be a (closed) point. Then $f$ induces an isomorphism of completions $\widehat{\sO}_{V,v}\cong \widehat{\sO}_{U,u}$. Since having non-multiplicative singularities is checked at the completion, the result follows.
\end{proof}

\begin{notation}
We define the following total order on $\bF_p$: $\lambda < \lambda'$ if and only if $\Lambda< \Lambda'$, where $\Lambda,\Lambda'\in \{0,\dots,p-1\}$ are the unique lifts of $\lambda$ and $\lambda'$ respectively.
\end{notation}

\begin{definition}\label{def:smallest_eigenvalue}
Let $(s\in S)$ be a germ of regular surface over $k$, and $\sF$ be a singular $1$-foliation of rank one on $S$ with multiplicative singularities. Consider the natural representation
		$$\rho_{\sF,u}\colon \sF/\fm_u\sF\longrightarrow
		\End_k(\fm_u/\fm_u^2).$$ 
We let
		\begin{equation}\label{eqn:smallest_eigenvalue}
		\lambda(\sF,u)=\min\left\{\lambda\in \bF_p \ \big| \ \ 
		\begin{matrix}
		\text{there is }\xi\in \sF/\fm_u\sF \text{ such that}\\
		\text{the eigenvalues of }\rho_{\sF,u}(\xi)
		\text{ are }\{1,\lambda\}.
		\end{matrix}
		\right\}
		\end{equation}
\end{definition}

\begin{remark}\label{rmk:unique_endomorphism}
Notations as in \autoref{def:smallest_eigenvalue}. Then if $p=2$, or if $\lambda(\sF,u)\neq -1$, there is a unique endomorphism in $\im(\rho_{\sF,u})$ with eigenvalues $1$ and $\lambda(\sF,u)$. If $p>2$ and $\lambda(\sF,u)=-1$ then there are exactly two (multiply one by $-1$ to get the other).
\end{remark}

\begin{lemma}\label{lemma:normal_form_smallest_eigenvalue}
Notations as in \autoref{def:smallest_eigenvalue}. Then $\lambda(\sF,u)$ is well-defined, and there are formal coordinates $x,y\in \widehat{\sO}_{U,u}$ such that $\sF\otimes \widehat{\sO}_{U,u}$ is generated by $x\partial_x+\lambda(\sF,u) y\partial_y$.
\end{lemma}
\begin{proof}
Since $\sF$ has a multiplicative singularity at $u$, we can find formal coordinates $x,y\in \widehat{\sO}_{U,u}$ such that $\sF\otimes \widehat{\sO}_{U,u}$ is generated by $x\partial_x+\lambda y\partial_y$ for some $\lambda\in\bF_p^\times$. Since $\sF/\fm_u\sF$ is one-dimensional over $k$, we deduce that the elements in the image of $\rho_{\sF,u}$ are simultaneously diagonalizable with eigenvalues $\{\alpha,\lambda\alpha\}$ where $\alpha\in k$. Thus $\lambda(\sF,u)=\min\{\lambda,\lambda^{-1}\}$. If $\lambda(\sF,u)=\lambda$ then the ordered pair of coordinates $x,y$ works. If $\lambda(\sF,u)=\lambda^{-1}$ then the ordered pair $y,x$ works, since $y\partial_y+\lambda^{-1}x\partial_x$ also generates $\sF\otimes \widehat{\sO}_{U,u}$.
\end{proof}

\begin{lemma}\label{lemma:smallest_eigenvalue_invariant}
Let $f\colon V\to U$ be an \'{e}tale morphism of regular surfaces, and let $\sF$ be a $1$-foliation of rank one on $U$. If $u\in \Sing(\sF)$ is a multiplicative singularity and $v\in f^{-1}(u)$, then $\lambda(\sF,u)=\lambda(f^*\sF,v)$.
\end{lemma}
\begin{proof}
Notice that $\lambda(f^*\sF,v)$ is well-defined, because $v$ is a multiplicative singularity of $f^*\sF$ by \autoref{lemma:pullback_non_lc_locus}. Let $x,y\in \widehat{\sO}_{U,u}$ be formal coordinates such that $\sF\otimes\widehat{\sO}_{U,u}$ is generated by $x\partial_x +\lambda(\sF,u) y\partial_y$, and let $w,z\in \widehat{\sO}_{V,v}$ be formal coordinates such that $f^*\sF\otimes \widehat{\sO}_{V,v}$ is generated by $w\partial_w+\lambda(f^*\sF,v) z\partial_z$. The morphism $f$ induces a $k$-linear isomorphism $\widehat{\sO}_{V,v}\cong \widehat{\sO}_{U,u}$, so $f^*x,f^*y$ is a system of parameters of $\widehat{\sO}_{V,v}$. Moreover, $f$ also induces a commutative diagram
		\begin{equation}\label{eqn:representation_invariant}
		\begin{tikzcd}
		\sF/\fm_u\sF \arrow[d, "\cong"] \arrow[rr, "\rho_{\sF,u}"] && \End_k(\fm_u/\fm_u^2)\arrow[d, "\cong"] \\
		f^*\sF/\fm_v f^*\sF  \arrow[rr, "\rho_{f^*\sF,v}"] && \End_k(\fm_v/\fm_v^2)
		\end{tikzcd}
		\end{equation}
where the vertical arrows are $k$-linear isomorphisms. Since eigenvalues of $k$-linear endomorphisms are preserved by $k$-linear conjugations, it follows from the definition \autoref{eqn:smallest_eigenvalue} that $\lambda(\sF,u)=\lambda(f^*\sF,v)$.
\end{proof}

\begin{corollary}\label{cor:Rees_algebra_invariant}
Notations as in \autoref{lemma:smallest_eigenvalue_invariant}. Let $x,y\in \widehat{\sO}_{U,u}$ and $z,w\in \widehat{\sO}_{V,v}$ be formal parameters such that $\sF\otimes\widehat{\sO}_{U,u}$ is generated by $x\partial_x +\lambda(\sF,u) y\partial_y$, and such that $f^*\sF\otimes\widehat{\sO}_{V,v}$ is generated by $w\partial_w+\lambda(f^*\sF,v) z\partial_z$. Assume that $p=2$, or that $\lambda(\sF,u)\neq -1$. Then
		$$f^*\left[(x,1)+(y,\Lambda(\sF,u))\right]
		= (w,1)+(z,\Lambda(f^*\sF,v))$$
as Rees $\widehat{\sO}_{V,v}$-algebras, where $\Lambda(\bullet,\bullet)\in \{1,\dots,p-1\}$ denotes the unique lift of $\lambda(\bullet,\bullet)$.
\end{corollary}
\begin{proof}
The left-hand side is equal to $(f^*x,1)+(f^*y,\Lambda(\sF,u))$. By \autoref{lemma:smallest_eigenvalue_invariant} we have $\Lambda(\sF,u)=\Lambda(f^*\sF,v)$. Denote this common value by $\Lambda$, and its reduction modulo $p$ by $\lambda$. \autoref{prop:Rees_algebra_is_functorial} shows that $(w,1)+(z,\Lambda)=(f^*x,1)+(f^*y,\Lambda)$ as soon as $w\partial_w+\lambda z\partial_z$ and $(f^*x) \partial_{f^*x}+\lambda (f^*y)\partial_{f^*y}$ reduce to the same endomorphism of $\fm_v/\fm_v^2$. This holds by the commutativity of \autoref{eqn:representation_invariant}, the fact that $\sF/\fm_v\sF$ is one-dimensional, the minimality of $\lambda$ and \autoref{rmk:unique_endomorphism}.
\end{proof}

\begin{theorem}\label{thm:functorial_resolution}
Let $\sS$ be any normal separated two-dimensional Deligne--Mumford stack of finite type over $k$, and $\sF\subset T_{\sS/k}$ be a $1$-foliation of rank one. Then there exists a proper birational morphism
		$$R(\sS,\sF)\colon \sS'\longrightarrow \sS$$
such that:
	\begin{enumerate}
		\item $\sS'$ is a regular two-dimensional Deligne--Mumford stack of finite type over $k$,
		\item $\sS'$ is tame if $\sS$ is so,
		\item $R(\sS,\sF)$ is an isomorphism over $\Reg(\sF)$,
		\item $R(\sS,\sF)^*\sF$ is regular.
	\end{enumerate}
Moreover, $R(\bullet,\bullet)$ is functorial with respect to \'{e}tale morphisms: if $f\colon \sT\to \sS$ is \'{e}tale, then 
		$$\begin{tikzcd}
		\sT' \arrow[r]\arrow[d, "R(\sT{,}f^*\sF)" left] & \sS'\arrow[d, "R(\sS{,}\sF)"] \\
		\sT\arrow[r, "f"] & \sS
		\end{tikzcd}$$
is cartesian.
\end{theorem}
\begin{proof}
We divide the proof in two steps: the construction of $R(\bullet,\bullet)$, and the proof of functoriality.

\medskip
\textsc{Step 1: Construction of $R$.} We divide this step in two parts.

\textsc{Step 1.1.} \emph{First, assume that $\sS$ is regular.} Choose an \'{e}tale groupoid presentation $\sS=[s,t\colon V\rightrightarrows U]$. Then $U$ and $V$ are regular surfaces (possibly disconnected). Let $\sF_U$ and $\sF_V$ be the $1$-foliations on $U$ and $V$ corresponding to $\sF$. Recall that $s^*\sF_U=\sF_V=t^*\sF_U$. By \autoref{lemma:pullback_non_lc_locus} we have
		$$s^{-1}\left(\NMSing(\sF_U)\right)=\NMSing(\sF_V)=t^{-1}\left(\NMSing(\sF_U)\right).$$
In particular, we have isomorphisms
		$$V\times_{s,U}\Bl_{\NMSing(\sF_U)}(U)\cong \Bl_{\NMSing(\sF_V)}(V)\cong V\times_{t,U}\Bl_{\NMSing(\sF_U)}(U)$$
over $V$. By \autoref{lemma:pullback_groupoid} we obtain an \'{e}tale groupoid $\Bl_{\NMSing(\sF_V)}(V)\rightrightarrows \Bl_{\NMSing(\sF_U)}(U)$ whose associated quotient stack $\sS_1$ is projective and birational over $\sS$. In other words, $\sS_1\to \sS$ is the blow-up of $\sS$ along the image of $\NMSing(\sF_U)$. The stack $\sS_1$ is still regular separated and Deligne--Mumford, and tame if $\sS$ is so.

We continue blowing-up the finite set of non-multiplicative singularities of $\sF$. By \autoref{thm:surface_foliations_resolution}, after finitely many steps we obtain a projective birational morphism
		$$R_{+,\pre}(\sS,\sF)\colon \sS_{+,\pre}\longrightarrow \sS$$
such that $R_{+,\pre}(\sS,\sF)^*\sF=\sF_{+,\pre}$ has only multiplicative singularities. 

The set of (multiplicative) singularities of $\sF_{+,\pre}$ is finite, and the quantity $\lambda(\sF_{+,\pre},\bullet)$ is well-defined at each of them by \autoref{lemma:smallest_eigenvalue_invariant}. Let 
		$$Q_{\pre}(\sS_{+,\pre},\sF_{+,\pre})\colon
		\sS_+\longrightarrow \sS_{+,\pre}$$
be the (usual) blow-up of every point $u\in \Sing(\sF_{+,\pre})$ for which $\lambda(\sF_{+,\pre},u)=-1$. Let $\sF_+=Q_{\pre}(\sS_{+,\pre},\sF_{+,\pre})^*\sF_{+,\pre}$: using \autoref{eqn:pullback_blowup_surface}, an \'{e}tale-local computation shows that $\sF_+$ is regular if $p=2$, and that $\lambda(\sF_+,u)\neq -1$ for every $u\in \Sing(\sF_+)$ if $p>2$.

If $p=2$ we go to Step 1.2. Otherwise, assuming that $\Sing(\sF_+)$ is non-empty, we want to perform a weighted blow-up as in the proof of \autoref{thm:stacky_resolution_surfaces}. The fact that $\sS_+$ is a stack is not an issue, cf.\ \cite{Quek_Rydh_Weighted_blowups}, and we only have to produce the Rees $\sO_{\sS_+}$-algebra that we want to blow-up. Let $[s_+,t_+\colon V_+\rightrightarrows U_+]$ be an \'{e}tale groupoid presentation of $\sS_+$.  By descent, we have to produce a Rees algebra $\fI_{U_+}$ on $U_+$ with the property that $s_+^*\fI_{U_+}=t_+^*\fI_{U_+}$ as Rees $\sO_{V_+}$-algebras.

We denote by $\sF_{U_+}$ and $\sF_{V_+}$ the pullbacks of $\sF_+$ to $U_+$ and $V_+$. The Rees algebra $\fI_{U_+}$ will be co-supported on $\Sing(\sF_{U_+})$: as in the proof of \autoref{thm:stacky_resolution_surfaces}, it suffices to specify its value after completion at every point of $\Sing(\sF_{U_+})$. So if $u\in \Sing(\sF_{U_+})$, let $x,y$ be the formal coordinates given by \autoref{lemma:normal_form_smallest_eigenvalue} applied to $\sF_{U_+}$, and let 
	$$\fI_{U_+}\otimes \widehat{\sO}_{U_+,u}=(x,1)+(y,\Lambda(\sF_{U_+},u))$$ 
where $\Lambda(\sF_{U_+},u)\in \{1,\dots,p-1\}$ is the unique lift of $\lambda(\sF_{U_+},u)$.

The fact that $s_+^*\fI_{U_+}=t_+^*\fI_{U_+}$ follows at once from \autoref{cor:Rees_algebra_invariant}. By descent, we obtain a Rees $\sO_{\sS_+}$-algebra $\fI$.

Finally, let $Q(\sS_+,\sF_+)\colon \sS'=\Bl_\fI(\sS_+)\to \sS_+$ be the blow-up of $\fI$. Then $\sS'$ is regular and separated, and tame if $\sS$ is so. Etale-local computations, as in \autoref{claim:weighted_blow_up_surface}, show that $Q(\sS_+,\sF_+)^*\sF_+$ is regular. Hence in case $\sS$ is regular, the composition 
	$$R(\sS,\sF)=Q(\sS_+,\sF_+)\circ Q_{\pre}(\sS_{+,\pre},\sF_{+,\pre})\circ R_+(\sS,\sF)\colon \sS'\longrightarrow \sS$$ 
provides a resolution. For the rest of the proof, we will write $Q_+(\sS,\sF)=Q(\sS_+,\sF_+)\circ Q_{\pre}(\sS_{+,\pre},\sF_{+,\pre})$.

\medskip
\textsc{Step 1.2.} \emph{Construction of $R$, general case.} If $\sS$ is only normal, let $\pi_{\sS}\colon \sS_0\to \sS$ be the resolution of singularities of $\sS$ provided by \cite[Theorem 17.5]{Cossart_Jannsen_Saito_Desingularization_inv_and_strategy}. Then we let 
	$$R(\sS,\sF)=R(\sS_0,\pi_{\sS}^*\sF)\circ \pi_{\sS}.$$ 
By construction, this morphism satisfies the conditions (a)-(d) of the statement (\footnote{
		The fact that $\sS_0$ is tame if $\sS$ is so is not stated in \cite{Cossart_Jannsen_Saito_Desingularization_inv_and_strategy}, but follows from the fact that the desingularization morphism obtained there is a sequence of usual blow-ups along the singular locus, see \cite[Theorem 1.2]{Cossart_Jannsen_Saito_Desingularization_inv_and_strategy}.
}).

\medskip
\textsc{Step 2: Functoriality of $R$.} Let $f\colon \sT\to \sS$ be an \'{e}tale morphism. We can decompose $R(\sS,\sF)$ as a sequence
		$$\begin{tikzcd}
		\sS' \arrow[rr, "Q_+(\sS_+{,}\sF_+)"] &&
		\sS_+ \arrow[rr, "R_+(\sS_0{,} \sF_0)"] &&
		\sS_0 \arrow[rr, "\pi_\sS"] &&
		\sS
		\end{tikzcd}$$
where the morphisms are as in Step 1, and where for simplicity we write $\sF_0=\pi_\sS^*\sF$ and $\sF_+=R_+(\sS_0,\sF_0)^*\sF_0$. To prove the functoriality statement, it suffices to show that we have a diagram
		$$\begin{tikzcd}
		\sT' \arrow[r, "f'"]\arrow[d, "Q_+(\sT_+{,}f_+^*\sF_+)" left] &
		\sS' \arrow[d, "Q_+(\sS_+{,}\sF_+)"] \\
		\sT_+ \arrow[r, "f_+"] \arrow[d, "R_+(\sT_0{,}f_0^*\sF_0)" left] &
		\sS_+  \arrow[d, "R_+(\sS_0{,}\sF_0)"] \\
		\sT_0 \arrow[r, "f_0"] \arrow[d, "\pi_{\sT}" left] &
		\sS_0 \arrow[d, "\pi_\sS"] \\
		\sT \arrow[r, "f"] & \sS
		\end{tikzcd}$$
where each one of the three small squares is cartesian.

The bottom square, involving $\pi_\sS$ and $\pi_\sT$, is cartesian since the resolution morphism constructed in \cite[\S 17]{Cossart_Jannsen_Saito_Desingularization_inv_and_strategy} is functorial with respect to \'{e}tale morphisms (see in particular \cite[Theorem 17.2]{Cossart_Jannsen_Saito_Desingularization_inv_and_strategy}).

The middle square is cartesian because the morphisms $R_+$ are repeated blow-ups of the non-multiplicative loci, which are \'{e}tale-invariant by \autoref{lemma:pullback_non_lc_locus}.

Finally, the top square is cartesian by construction of $Q$ and by \autoref{lemma:smallest_eigenvalue_invariant} and \autoref{cor:Rees_algebra_invariant}. This completes the proof of the theorem.
\end{proof}

\begin{remark}[Non-uniqueness of $R$]\label{rmk:several_algorithms}
The function $R=Q\circ R_+\circ \pi$ constructed in the proof of \autoref{thm:functorial_resolution} depends on several choices. For example, it depends on the functorial resolution of singularities $\pi$ of two-dimension Deligne--Mumford stacks elaborated in \cite{Cossart_Jannsen_Saito_Desingularization_inv_and_strategy}, and on the algorithm $R_+$ of Giraud \cite{Giraud_Forme_normale_sur_une_surface_en_car_pos}. Both are based on usual scheme-theoretic blow-ups, and it is not unlikely that \emph{weighted} blow-ups can provide short-cuts in both cases (as they do for resolution of singularities in characteristic zero, see e.g.\ \cite{Abramovich_Birational_geometry_using_weighted_blowups}). Furthermore, the last piece $Q$ depends on the ordering of $\bF_p$ that we have chosen: any other \'{e}tale-invariant way of selecting a weighted blow-up leads to a potentially different transformation.
\end{remark}

\section{Resolution on threefolds}\label{section:resolution_3folds}
Next we consider the case of $1$-foliations on threefolds. We begin with the corank $1$ case, from which we will deduce the corank $2$ case. The case $p=2$ will mostly be treated separately.

\subsection{Corank $1$ case}
The strategy in the corank $1$ case is similar to the surface case.

\begin{theorem}\label{thm:resolution_3fold_rank_2}
Let $X$ be a normal threefold over an algebraically closed field $k$ of characteristic $p>0$, and let $\sF$ be a $1$-foliation of corank $1$ on $X$. Then there exists a birational projective morphism $f\colon X'\to X$ such that $X'$ is regular, and such that the induced $1$-foliation $f^*\sF$ has everywhere at worst multiplicative singularities.
\end{theorem}
\begin{proof}
The argument is in essence similar to that of \autoref{thm:surface_foliations_resolution}. By \cite{Cossart_Piltant_Resolution_of_singularities_for_3_folds_in_pos_char_I} we may assume that $X$ is regular. Since $\sF$ has corank $1$, by \autoref{lemma:quotient_is_normalization} the field extension $K(X/\sF)\supset K(X^{(1)})$ has degree $p$. So there exists $u\in K(X^{(1)})$ such that $K(X/\sF)=K(X^{(1)})[u^{1/p}]$. 

Let $D$ be the divisor of poles of $u$, and let $\{U_i\}$ be an affine open cover of $X^{(1)}$, sufficiently refined so that we can find $v_i\in \sO_{X^{(1)}}(U_i)$ with $\text{div}(v_i)=D\cap U_i$ (this is possible since $X^{(1)}$ is regular). Then $uv_i^p\in \sO_{X^{(1)}}(U_i)$ for each $i$, and on an intersection $U_{ij}$ the two sections $uv_i^p$ and $uv_j^p$ differ by the unit $(v_i/v_j)^p\in \sO_{X^{(1)}}(U_{ij})^*$. 

\begin{claim}\label{claim:Cossart_method_applies}
There is a projective birational morphism $f\colon V\to X^{(1)}$ such that: $V$ is regular, and \'{e}tale-locally at every $s\in V$ we can write
		\begin{equation}\label{eqn:etale_neighbourhoods_II}
		f^*(uv_i^p)=t^p+\bold{x}^{\bold{a}}
		\end{equation}
where $i$ is such that $f(s)\in U_i$, the collection $\bold{x}$ is part of a system of local parameters and $\bold{a}\neq \bold{0}$ modulo $p$.
\end{claim}

This is a slight generalization of the main result of Cossart's thesis \cite{Cossart_Polyhedre_caracteristique}. In order to prove the claim, a glimpse into the methods of \cite{Cossart_Polyhedre_caracteristique} is necessary: for the sake of readability, the necessary details will be given in the \autoref{appendix}.

We claim that $f^{(-1)}\colon V^{(-1)}\to X$ is the morphism we are looking for. As in the proof of \autoref{thm:surface_foliations_resolution}, we observe that 
		$$V^{(-1)}/f^{(-1),*}\sF=V\left[(f^*u)^{1/p}\right]^\nu.$$
Therefore, to simplify the notations, we replace $X$ by $V^{(-1)}$, etc. Let us now compute the normalization of $B=\widehat{\sO}_{X^{(1)},s}[T]/(T^p-uv_i^p)$ when $s\in f^{-1}(U_i)$ is a closed point: by \autoref{lemma:quotient_is_normalization} this normalization is the local ring of $X/\sF$ at $s$. We may assume that $t=0$.
	\begin{enumerate}
		\item Assume that for some $j$ we have $a_j=0$ modulo $p$. Say $j=3$: then we can partially normalize the germ by introducing the element $w=T/x_3^{a_3/p}\in B^\nu$, namely
			$$B[w]\cong \widehat{\sO}_{X^{(1)},s}[W]/(W^p-x_1^{a_1}x_2^{a_2}).$$
		As only two local parameters are involved in this cyclic cover, the normalization of $B[w]$ now reduces to one of the surface cases done in the proof of \autoref{thm:surface_foliations_resolution}. Similarly if $\bold{x}$ is a partial system of parameters. In particular we get that $B^\nu$ is regular, or that 
				\begin{equation}\label{eqn:surface_type_cases}
				\begin{array}{l}
				\text{there is a system of local parameters }x,y,z \text{ for }\widehat{\sO}_{X^{(-1)},s}\text{ such that} \\
				B^\nu \cong \left(
				k\llbracket x,y,w \rrbracket [T]/(T^p-x^{a}y^{b})\right)^\nu
				\text{ with }a,b\text{ coprime to }p.
				\end{array}
				\end{equation}
		\item The remaining case reads $\bold{x}=(x_1,x_2,x_3)$ and $a_i\neq 0$ modulo $p$ for every $i$. In that case,
				\begin{equation}\label{eqn:threefold_type_case}
				\begin{array}{l}
				\text{there is a system of local parameters }x,y,z \text{ for }\widehat{\sO}_{X^{(1)},s}\text{ such that} \\
				B^\nu \cong \left(
				k\llbracket x,y,z\rrbracket [T]/(T^p-x^ay^bz^c)\right)^\nu \text{ where }a,b,c\text{ are coprime to }p.
				\end{array}
				\end{equation}
	\end{enumerate}
Recall that $B^\nu$ is the completed local ring of $X/\sF$ at $s$. So if $s\in X/\sF$ is regular, then $\sF$ is regular at $x$ by \cite[Lemma 2.5.10]{Posva_Singularities_of_quotients_by_foliations}. If $s\in X/\sF$ is of type \autoref{eqn:surface_type_cases} or \autoref{eqn:threefold_type_case}, then we apply \autoref{prop:inv_subring_normaliz_hypersurface} to $\sO_{X^{(1)},s}\hookrightarrow \sO_{X/\sF,s}\hookrightarrow \sO_{X,s}$ and obtain that the foliation $\sF$ is formally generated at $s$, up to saturation, by multiplicative derivations. This completes the proof.
\end{proof}

\begin{remark}\label{rmk:blow_up_over_reg_locus}
Similarly to what we observed in \autoref{rmk:u_vs_F}, the morphism $f$ we obtain in the proof of \autoref{thm:resolution_3fold_rank_2} is highly sensitive to the choice of the functions $v_i^ps$, and they may not have the \'{e}tale-local form \autoref{eqn:etale_neighbourhoods_II} at some point $x\in X^{(1)}$ while $\sF$ has multiplicative singularities at $x\in X$. Since the blowing-up procedure of \cite{Cossart_Polyhedre_caracteristique} is more intricate than the two-dimensional case of \cite{Giraud_Forme_normale_sur_une_surface_en_car_pos}, we leave open whether we can ensure $f^{(-1)}$ is an isomorphism over the locus where $\sF$ has at worst multiplicative singularities.
\end{remark}

As in the surface case, the next step is to use weighted blow-ups to improve the singularities of the $1$-foliation. The case $p=2$ can be dealt with using usual blow-ups, uniformly for the rank $1$ and $2$ cases, and we will do so in \autoref{thm:resolution_in_char_2} below.

\begin{theorem}\label{thm:stacky_resolution_rk_2}
Let $X$ be a normal threefold over an algebraically closed field of characteristic $p>2$, and let $\sF$ be a $1$-foliation of corank one on $X$. Then there exists a proper birational morphism $f\colon \sX'\to X$ such that $\sX'$ is a regular tame Deligne--Mumford stack and $f^*\sF$ is regular.
\end{theorem}
\begin{proof}
By \autoref{thm:resolution_3fold_rank_2} we may assume that $X$ is regular and $\sF$ has at worst multiplicative singularities. We will perform two rounds of weighted blow-ups: first we simplify $\Sing(\sF)$, then we resolve the singularities of $\sF$. Accordingly, we divide our proof in two steps.

\medskip
\textsc{Step 1: Resolving $\Sing(\sF)$.} For an arbitrary closed point $s\in \Sing(\sF)$, we can find \'{e}tale-local coordinates $x,y,z\in\widehat{\sO}_{X,s}$ such that $\sF\otimes \widehat{\sO}_{X,s}$ is generated
		$$\begin{cases}
		\text{either by }&
		 \partial_z,\ x\partial_x+\lambda y\partial_y,
		\\
		\text{or by }& \partial=x\partial_x+\lambda y\partial_y,\ \psi=x\partial_x+\mu z\partial_z,
		\end{cases}$$
for appropriate $\lambda,\mu\in\bF_p^\times$. So we see that $\Sing(\sF)$ is of pure dimension one, and that $s$ is a singular point of $\Sing(\sF)$ if and only if we are in the second case. 

So the closed set $\Sing(\Sing(\sF))$ is finite. Let $s$ be any point in it. We will perform a weighted blow-up centered at $s$, defined in terms of the coordinates $x,y,z$: it will be algebraizable by the argument given in the proof of \autoref{thm:stacky_resolution_surfaces}, so let us only explain which weighted blow-ups we shall perform.

Choose integers $b,c\geq 1$ coprime to $p$ such that $\lambda=2b$ and $\mu=2c$ hold in $\bF_p$ (recall that $p\neq 2$), and blow-up the ideal $(x,1)+(y,b)+(z,c)$. Let us compute the pullbacks of $\partial$ and $\psi$, using the computations of \autoref{example:weighted_blow_ups}:
	\begin{itemize}
		\item On the $D_+(x)$-patch with coordinates $u,v,w$, one finds
				$$\partial=u\partial_u+(\lambda -b)v\partial_v-cw\partial_w=u\partial_u+bv\partial_v-cw\partial_w$$
		and
				$$\psi=u\partial_u-bv\partial_v+(\mu-c)w\partial_w=u\partial_u-bv\partial_v+cw\partial_w.$$
		Therefore the saturated pullback of $\sF$ on $D_+(x)$ is generated by $\partial_u$ and by $y\partial_y+\lambda z\partial_z$ for an appropriate $\lambda\in\bF_p$. 
		\item On the $D_+(y)$-patch with coordinates $u,v,w$, one finds similarly that
				$$\partial=-u\partial_u-2v\partial_v-\mu w\partial_w, \quad 
				\psi=u\partial_u+\mu w\partial_w$$
		and so the saturated pullback of $\sF$ is generated by $\partial_v$ and $u\partial_u+\mu w\partial_w$.
		\item The $D_+(z)$-patch is similar to the $D_+(y)$ one.
	\end{itemize}
Therefore, if we perform such weighted blow-ups at every singular point of $\Sing(\sF)$, we obtain a proper model $g\colon\sX\to X$ where $\sX$ is a regular tame Deligne--Mumford stack, such that the singular locus of $g^*\sF$ is regular.

\medskip
\textsc{Step 2: Resolving $\sF$.}
Refreshing our notations, we have a regular tame Deligne--Mumford $X$ together with a $1$-foliation $\sF\subset T_{X/k}$ or corank one with the following properties:
	\begin{itemize}
		\item $X$ is schematic at the generic points of $\Sing(\sF)$;
		\item $\Sing(\sF)$ is regular one-dimensional;
		\item we can find \'{e}tale morphisms $f_i\colon U_i\to X$, where the $U_i$'s are schemes, such that $\Sing(\sF)$ is contained in the union of the images of the $U_i$'s, and such that for each closed $s\in f_i^{-1}\Sing(\sF)=\Sing(f_i^*\sF)$, we have preferred coordinates $x,y,z\in \sO_{U_i,s}$ such that 
	\begin{equation}\label{eqn:preferred_parameters}
	\sF\otimes \sO_{U_i,s}=\sO_{U_i,s}\cdot \partial_z +\sO_{U_i,s}\cdot (x\partial_x+\lambda_s y\partial_y)
	\end{equation}			
for some $\lambda_s\in\bF_p^\times$. (These $U_i$ are the coordinate charts of the weighted blow-ups performed in Step 1.)
	\end{itemize}	
If $x,y,z\in \sO_{U_i,s}$ are as in \autoref{eqn:preferred_parameters}, then the blow-up of $(x,y)$ with weights $(1,\Lambda_s)$ will resolve the singularities of $\sF\otimes \sO_{U_i,s}$, where $\Lambda_s\in\{1,\dots,p-1\}$ lifts $\lambda_s\in\bF_p^\times$: the computation is just as in the proof of \autoref{thm:stacky_resolution_surfaces}. The problem is that these coordinates are not unique (\footnote{
	Consider for example the new coordinates $u=x+y^{d},v=y$, where $d\in\{1,\dots,p-1\}$ denotes the lift of $\lambda^{-1}\in\bF_p^\times$. Then $x\partial_x+\lambda y\partial_y=u\partial_u+\lambda v\partial_v$. More generally, see the proof of \autoref{prop:Rees_algebra_is_functorial}.
		}) 
and $\lambda_s$ is not uniquely defined (\footnote{
		For example, switching $x$ and $y$ replaces $\lambda_s$ by $\lambda_s^{-1}$. 
		}).
Thus we need to be careful if we want to define a global weighted blow-up of $X$ along $\Sing(\sF)$. 

To begin with, we prove that $\lambda_s$ is uniquely defined up to inversion.

\begin{claim}\label{claim:generic_value_of_lambda}
The set-valued function $s\mapsto \{\lambda_s,\lambda_s^{-1}\}$, where $s$ runs through the closed points of $\Sing(f_i^*\sF)$, is constant along the preimage of any irreducible component of $\Sing(\sF)$.
\end{claim}
\begin{proof}\renewcommand{\qedsymbol}{$\lozenge$}
Let $E\subset \Sing(f_i^*\sF)$ be an irreducible component. If $\eta$ is the generic point of $E$, consider the sub-$\sO_{U_i,\eta}$-module
		$$\fI(\sF,E)=(\sF\otimes \sO_{U_i,\eta})\cap \fm_\eta T_{U_i/k,\eta}.$$
By definition, for $\partial\in \fI(\sF,E)$ we have $\partial(\fm_\eta)\subseteq \fm_\eta$. Thus we obtain a $k(\eta)$-linear representation
		\begin{equation*}
		\rho_{\sF,E}\colon \fI(\sF,E)/\fm_\eta \fI(\sF,E)\longrightarrow \End_{k(\eta)}(\fm_\eta/\fm_\eta^2).
		\end{equation*}
If $s\in E$ is a closed point, we have a preferred set of coordinates $x,y,z$ of $\sO_{U_i,s}$ for $\sF\otimes\sO_{U_i,s}$ as described in \autoref{eqn:preferred_parameters}, and the ideal $I_E\otimes\sO_{U_i,s}$ is $(x,y)$. So localizing at this ideal and writing down the representation $\rho_{\sF,E}$ in terms of these coordinates, one sees that 
	$$\dim_{k(\eta)}\fI(\sF,E)/\fm_\eta\fI(\sF,E)=1,$$ 
and that the image $\xi$ of $x\partial_x+\lambda_s y\partial_y$ generates this vector space. The eigenvalues of $\rho_{\sF,E}(\xi)$ are precisely $1$ and $\lambda_s$. As the eigenvalues of $\rho_{\sF,E}(a\cdot \xi)$ are $\{a,a\lambda_s\}$ for any $a\in k(\eta)$, we see that up to inversion the value of $\lambda_s$ is determined by finding an element of $\im(\rho_{\sF,E})$ which has $1$ as eigenvalue. As $\im(\rho_{\sF,E})$ does not depend on $s$, this shows that the function $s\mapsto \{\lambda_s,\lambda_s^{-1}\}$ is locally constant on the components of $\Sing(f_i^*\sF)$. Finally, observe that the eigenvalues of $\im(\rho_{\sF,E})$ can be found after passing to the completion $\widehat{\sO}_{U_i,\eta}$. As $f_i$ is \'{e}tale, this computation can be carried out directly at the completion of $\sO_{X,f_i(\eta)}$ (recall that $f_i(\eta)\in X$ is a schematic point). Thus the function $s\mapsto \{\lambda_s,\lambda_s^{-1}\}$ depends only on $f_i(\eta)$. This proves our claim.
\end{proof}

Thanks to the claim, we can assign to each irreducible component $\sE\subset \Sing(\sF)$ a value $\lambda_\sE\in\bF_p^\times$, such that for any $f_i\colon U_i\to X$ the constant value of $s\mapsto \{\lambda_s,\lambda_s^{-1}\}$ along $f_i^{-1}\sE$ is equal to $\{\lambda_\sE,\lambda_\sE^{-1}\}$. (We have to choose between an element and its inverse, but this choice will not matter.) We let $\Lambda_\sE\in\{1,\dots,p-1\}$ be the lift of $\lambda_\sE\in\bF_p^\times$.

\medskip
Next, we want to show that the Rees $\sO_{U_i,s}$-algebras $(x,1)+(y,\Lambda_\sE)$ can be glued together as $s$ moves along $\Sing(f_i^*\sF)$. This can be be checked at the generic points of $\Sing(f_i^*\sF)$. It is useful to introduce the following terminology:

Let $E$ be an irreducible component of $\Sing(f_i^*\sF)$ which dominates the component $\sE$ of $\Sing(\sF)$. At the generic point $\eta\in E$, we say that an ordered system of coordinates $x,y\in \sO_{U_i,\eta}$ is \emph{adapted to $(\sF,\lambda_\bullet)$} if $\fI(\sF,E)$ is generated by $x\partial_x +\lambda_\sE y\partial_y$. 

\begin{claim}
Every generic point of $\Sing(f_i^*\sF)$ has an ordered system of coordinates that is adapted to $(\sF,\lambda_\bullet)$.
\end{claim}
\begin{proof}\renewcommand{\qedsymbol}{$\lozenge$}
By \autoref{eqn:preferred_parameters} and \autoref{claim:generic_value_of_lambda} there is an ordered system of parameters $x,y,z$ of $\sO_{U_i,s}$ such that either
		$$x\partial_x+\lambda_\sE y\partial_y \quad \text{or}
		\quad x\partial_x+\lambda_\sE^{-1} y\partial_y$$
generates $\fI(\sF,E)$. In the first case $x,y$ is adapted. In the second case $y,x$ is adapted.
\end{proof}

\begin{claim}\label{claim:weighted_ideal_indep_of_choices}
Let $\eta$ be a generic point of $\Sing(f_i^*\sF)$. Assume that $u,v\in \sO_{U_i,\eta}$ and $x,y\in\sO_{U_i,\eta}$ are adapted to $(\sF,\lambda_\bullet)$. Then $(u,1)+(v,\Lambda_\sE)=(x,1)+(y,\Lambda_\sE)$ as Rees $\sO_{U_i,\eta}$-algebras.
\end{claim}
\begin{proof}\renewcommand{\qedsymbol}{$\lozenge$}
By extension-retraction for ideals along faithfully flat extensions \cite[Theorem 7.5.(ii)]{Matsumura_Commutative_Ring_Theory}, it suffices to prove equality at the completion $\widehat{\sO}_{U_i,\eta}$. Write $\xi=x\partial_x+\lambda y\partial_y$ and $\xi'=u\partial_u+\lambda v\partial_v$. As in \autoref{claim:generic_value_of_lambda}, consider the $k(\eta)$-linear representation
		$$\rho_{\sF,E}\colon \fI(\sF,E)/\fm_\eta \fI(\sF,E)\longrightarrow \End_{k(\eta)}(\fm_\eta/\fm_\eta^2).$$
Since $\fI(\sF,E)/\fm_\eta\cdot \fI(\sF,E)$ is one-dimensional, and since $\rho_{\sF,E}(\xi)$ and $\rho_{\sF,E}(\xi')$ have the same two distinct eigenvalues (namely $1$ and $\lambda$), we deduce that $\rho_{\sF,E}(\xi)=\rho_{\sF,E}(\xi')$. So we are in position to apply \autoref{prop:Rees_algebra_is_functorial} and conclude.
\end{proof}

Now fix an irreducible component $\sE$ of $\Sing(\sF)$ and let $E$ be an irreducible component of $\Sing(f_i^*\sF)$ which dominates $\sE$. Let $s,s'\in E$ be two closed points with preferred set of local parameters $x,y,z\in\sO_{U_i,s}$ and $x',y',z'\in \sO_{U_i,s'}$ as in \autoref{eqn:preferred_parameters}: changing the order of the parameters if necessary, we may assume that $\lambda_s=\lambda_\sE=\lambda_{s'}$. By \autoref{claim:weighted_ideal_indep_of_choices} the Rees algebras
		$$(x,1)+(y,\lambda_\sE)\subset \sO_{U_i,s}[t]\quad 
		\text{and}\quad
		(x',1)+(y',\lambda_\sE)\subset \sO_{U_i,s'}[t]$$
induce, by extension of scalars to $\sO_{U_i,E}$, the \emph{same} Rees $\sO_{U_i,E}$-algebra. It follows that there exists a Rees algebra $\fJ_{\sE,i}\subset \sO_{U_i}[t]$ satisfying
		$$\fJ_{\sE,i}\otimes\sO_{U_i,s}=(x,1)+(y,\lambda_\sE)\quad \text{for all }s\in f_i^{-1}(\sE)$$
where $x,y\in \sO_{U_i,s}$ are as in \autoref{eqn:preferred_parameters}. Furthermore, it also follows from \autoref{claim:weighted_ideal_indep_of_choices} that on an \'{e}tale intersection $U_{ij}=U_i\times_X U_j$, we have $\fJ_{\sE,i}|_{U_{ij}}=\fJ_{\sE,j}|_{U_{ij}}$ (because all the computations are performed at the completed stalks). Thus $\{\fJ_{\sE,i}\}_i$ descends to a Rees $\sO_X$-algebra $\fJ_\sE\subset \sO_X[t]$. 

Let $\sX=\Bl_{\fJ_\sE}(X)$. It is a tame regular Deligne--Mumford stack. For every closed point $s\in \Sing(f_i^*\sF)$ we have
	$$\sX\otimes_X\sO_{U_i,s}=\Bl_{(x,1)+(y,\lambda_\sE)}\Spec(\sO_{U_i,s})$$ 
where $x,y\in\sO_{U_i,s}$ are as in \autoref{eqn:preferred_parameters}. So as observed before \autoref{claim:generic_value_of_lambda} the pullback of $\sF$ to $\sX$ is regular along the preimage of $\sE$.

Repeating this operation for every component $\sE$ of $\Sing(\sF)$, we obtain the model $\sX'\to X$ with the desired properties. This completes the proof.
\end{proof}

\begin{remark}\label{rmk:functorial_resolution_dim_3?}
I do not know whether the proofs of \autoref{thm:resolution_3fold_rank_2} and \autoref{thm:stacky_resolution_rk_2} can be extended to $1$-foliations on three-dimensional Deligne--Mumford stacks. A first issue is that we do not currently have a resolution theorem for three-dimensional stacks, as the method of \cite{Cossart_Piltant_Resolution_of_singularities_for_3_folds_in_pos_char_I} is not functorial with respect to, say, \'{e}tale morphisms. But even if we restrict to regular stacks, a more serious difficulty is that the methods of \cite{Cossart_Polyhedre_caracteristique} are not functorial either (see \cite[II.C.5.4]{Cossart_Polyhedre_caracteristique}), and so it is unclear if we can define the analogue of the transformation $R_+\colon \sS_+\to\sS$ from the proof of \autoref{thm:functorial_resolution}.
\end{remark}

\subsection{Corank $2$ case}
The corank $2$ case is an easy application of the stacky Jacobson's correspondence and of the corank $1$ case.

\begin{theorem}\label{thm:stacky_resolution_rk_1}
Let $Y$ be a normal threefold over an algebraically closed field of characteristic $p>2$, and let $\sG$ be a $1$-foliation of corank two on $Y$. Then there exists a proper birational morphism $f\colon \sY'\to Y$ such that $\sY'$ is a regular tame Deligne--Mumford stack and $f^*\sG$ is everywhere regular.
\end{theorem}
\begin{proof}
Let $X^{(1)}=Y/\sG$. Then the induced morphism $X\to Y$ has degree $p^2$ and so corresponds to a $1$-foliation $\sF$ of corank one on $X$ by the schematic Jacobson correspondence \cite[Theorem 2.5.12]{Posva_Singularities_of_quotients_by_foliations}. We apply \autoref{thm:stacky_resolution_rk_2} to $(X,\sF)$: there is a proper birational morphism $\sX'\to X$ such that $\sX'$ is a tame regular Deligne--Mumford stack, such that the pullback $\sF'$ of $\sF$ to $\sX'$ is everywhere regular. Let $\sY'=\sX'/\sF'$ (see \autoref{def:infinitesimal_qt_stack}). Then we have an induced birational proper morphism $f\colon \sY'\to Y$ and a commutative diagram
		$$\begin{tikzcd}
		\sX'\arrow[r]\arrow[d] & \sY'\arrow[d, "f"]\arrow[r] & \sX'^{(1)}\arrow[d] \\
		X\arrow[r] & Y\arrow[r] & X^{(1)}.
		\end{tikzcd}$$
Since $\sF'$ is regular, it follows from \autoref{lemma:regularity_of_qt} that $\sY'$ is regular. It is tame Deligne--Mumford by \autoref{lemma:qt_preserves_tameness}. By \autoref{thm:stacky_Jacobson} we have $\sY'/f^*\sF=\sX'^{(1)}$, and another application of \autoref{lemma:regularity_of_qt} shows that $f^*\sG$ is regular everywhere. This concludes the proof. 
\end{proof}

\begin{remark}[Schematic counterpart]\label{rmk:schematic_rk_1_result}
One can wonder whether we can obtain an analogue of \autoref{thm:resolution_3fold_rank_2} for $1$-foliations $\sG$ of corank $2$ on threefolds (i.e.\ a resolution statement that does not involve stacks as underlying spaces). 

By either taking the coarse moduli of $\sY'$ in \autoref{thm:stacky_resolution_rk_1}, or the infinitesimal quotient $X'/f^*\sF$ in \autoref{thm:resolution_3fold_rank_2} (where $\sF$ is as in the proof of \autoref{thm:stacky_resolution_rk_1}), one obtains the following: given a corank two $1$-foliation $\sG$ on a normal threefold $Y$, there exists a projective birational $g\colon Y'\to Y$ such that $Y'$ has only (possibly non-isolated) linearly reductive quotient singularities, and $g^*\sG$ is locally and up-to-saturation generated by multiplicative derivations. 

If we insist on having a regular underlying variety, then the problem is to lift of the local multiplicative derivations generating $g^*\sG$ to a resolution of $Y'$, aiming to use \cite[Remark 2.4.9]{Posva_Singularities_of_quotients_by_foliations} on the resolution. From this perspective, it might be better to use $Y'=X'/f^*\sF$, for then we have very explicit local descriptions of $\Sing(Y')$, given by \autoref{eqn:surface_type_cases} and \autoref{eqn:threefold_type_case}. The local description of $g^*\sG$ is then afforded through \autoref{prop:inv_subring_normaliz_hypersurface}. This is very handy in case $p=2$, as we shall see in \autoref{thm:resolution_in_char_2} below.

Take $Y'=X'/f^*\sF$ and assume that $\Sing(f^*\sF)$ is a regular curve (on $X'$). Then we can show that the local multiplicative derivations defining $g^*\sG$ indeed lift to a resolution of $Y'$. Indeed, in this case the singularities of $Y'$ are of type \autoref{eqn:surface_type_cases}, which by \autoref{prop:inv_subring_normaliz_hypersurface} can be written
		$$k\llbracket x^iy^j, z^p\mid i+\lambda j=0 \ (p)\rrbracket, \quad \lambda\in\bF_p^\times,$$
and $g^*\sG$ is generated up-to-saturation by $x\partial_x$. Such a singularity is formal-locally the product of a smooth curve germ $k\llbracket z^p\rrbracket$ with a germ of surface linearly reductive quotient singularity $k\llbracket x^iy^j\mid i+\lambda j=0 \ (p)\rrbracket$. So a canonical resolution is afforded by the minimal resolution of the surface germ, and by \autoref{prop:T-equivariance_toric_surfaces} the multiplicative derivation $x\partial_x$ lifts to the resolution.

In case $\Sing(f^*\sF)$ is not regular, I do not know if we can lift the local generators of $g^*\sG$ to a resolution of $Y'$. A generalization of \autoref{prop:T-equivariance_toric_surfaces} to dimension $3$ would give a positive answer to this question.
\end{remark}

\subsection{$p=2$ case}
If $p=2$, then we do not need weighted blow-ups to resolve $1$-foliation singularities.

\begin{theorem}\label{thm:resolution_in_char_2}
Let $X$ be a normal threefold over a perfect field $k$ of characteristic $p=2$, and let $\sF$ be a $1$-foliation on $X$. Then there exists a birational projective morphism $f\colon X'\to X$ such that $X'$ is a regular threefold, and such that the induced $1$-foliation $f^*\sF$ is regular.
\end{theorem}
\begin{proof}
We start with the rank $1$ case, using the method sketched in \autoref{rmk:schematic_rk_1_result} above. The rank $2$ case will follow using the same "Jacobson trick" as in \autoref{thm:stacky_resolution_rk_1}.

\medskip
\textsc{Case 1: $\sF$ has rank $1$.} If $Y=X/\sF$, then $X=Y^{(-1)}/\sG$ for a $1$-foliation $\sG$ of corank $1$. The proof of \autoref{thm:resolution_3fold_rank_2}, applied to $(Y^{(-1)},\sG)$, shows that after replacing $X$ by a birational projective model, we may assume that at those points $s\in X$ where $X$ is not regular, we are in one of the following situations:
	\begin{enumerate}
		\item \autoref{eqn:surface_type_cases}: we have $\widehat{\sO}_{X,s}\cong k\llbracket x^2,xy,y^2,z\rrbracket\cong k\llbracket u,v,w,t\rrbracket /(t^2-uv)$, and under this isomorphism $\widehat{\sF}_s$ is generated up to saturation by $t\frac{\partial}{\partial t}$;
		\item \autoref{eqn:threefold_type_case}: we have $\widehat{\sO}_{X,x}\cong k\llbracket x^2,y^2,z^2,xyz\rrbracket\cong k\llbracket u,v,w,t\rrbracket /(t^2-uvw)$, and under this isomorphism $\widehat{\sF}_s$ is generated up to saturation by $t\frac{\partial}{\partial t}$.
	\end{enumerate}
Now I claim the following:
	
\begin{claim}
If we repeatedly blow-up the irreducible components of $\Sing(X)$ (one at a time, with its reduced structure), then after finitely many steps we obtain a resolution of $X$ and the formal local up-to-saturation generators of $\sF$ lift along every blow-up.
\end{claim}	
\begin{proof}\renewcommand{\qedsymbol}{$\lozenge$}
Since we have a global procedure to follow, it suffices to compute what happens at the completions. We consider again the two cases from above.
	\begin{enumerate}
		\item In the first case, the singular locus is given by $u=v=t=0$. So we have to compute what happens when we blow it up: we do so on the affine space model $\Spec(k[u,v,w,t]/(t^2-uv))\subset \bA^4$.
			\begin{enumerate}
				\item One chart is given by
			$$\frac{k[u,v,w,t]}{(t^2-uv)}\longrightarrow
			\frac{k[u,v',w,t']}{(t'^2-v')},\quad
			(u,v,w,t)\mapsto (u,uv',w,ut').$$
	This blow-up chart is regular. By \autoref{eqn:pullback_blowup_threefold}, one sees that $t\frac{\partial}{\partial t}$ lifts to $t'\frac{\partial}{\partial t'}$. 
				\item A second chart is given by
						$$\frac{k[u,v,w,t]}{(t^2-uv)}\longrightarrow
			\frac{k[u',v',w,t]}{(1-u'v')},\quad
			(u,v,w,t)\mapsto (tu',tv',w,t).$$
			This blow-up chart is regular. By \autoref{eqn:pullback_blowup_threefold}, one sees that $t\frac{\partial}{\partial t}$ lifts to $t\frac{\partial}{\partial t}-u'\frac{\partial}{\partial u'}-v'\frac{\partial}{\partial v'}$. 
				\item The remaining chart is similar to the first one.
			\end{enumerate}
		\item In the second case, the singular locus is the union of the three lines $(u=v=t=0),(u=w=t=0)$ and $(v=w=t=0)$. Our procedure will start by blowing-up any one of them: by symmetry, we may assume that we blow-up $(u=v=t=0)$. We use the affine model $\Spec(k[u,v,w,t]/(t^2-uvw))$ for our computations just as above.
			\begin{enumerate}
				\item One chart is given by
						$$\frac{k[u,v,w,t]}{(t^2-uvw)}\longrightarrow
			\frac{k[u,v',w,t']}{(t'^2-v'w)},\quad
			(u,v,w,t)\mapsto (u,uv',w,ut').$$
			The derivation $t\frac{\partial}{\partial t}$ lifts to $t'\frac{\partial}{\partial t'}$.
				\item A second chart is given by
						$$\frac{k[u,v,w,t]}{(t^2-uv)}\longrightarrow
			\frac{k[u',v',w,t]}{(1-u'v'w)},\quad
			(u,v,w,t)\mapsto (tu',tv',w,t).$$
			This chart is regular. The derivation $t\frac{\partial}{\partial t}$ lifts to $t\frac{\partial}{\partial t}-u'\frac{\partial}{\partial u'}-v'\frac{\partial}{\partial v'}$. 
				\item The remaining chart is similar to the first one.
			\end{enumerate}
	\end{enumerate}
In particular, we see that one blow-up transforms the singularities \autoref{eqn:threefold_type_case} into singularities of type \autoref{eqn:surface_type_cases}, which in turn are resolved by a single blow-up. Moreover, the derivations lift accordingly. This proves the claim.
\end{proof}

Refreshing our notations, we can now assume the following: $X$ is regular, and $\sF$ is everywhere formally generated up to saturation by a continuous multiplicative derivation. But $\sF$ might not be regular yet. I claim that $Z=\Sing(\sF)$ is regular, and that $\sF$ becomes regular on $\Bl_{Z}X$. %let me blow-up the ideal $\sI$ generated by the image of the evaluation pairing $\sO_X\times \sF\to \sO_X$. I claim that $\sI$ defines a regular subscheme of $X$, and that $\sF$ becomes regular on the blow-up (notice that any element of $\sF$ lifts to $\Bl_{\sI}X$ by \autoref{prop:blowup_and_derivations}). 
It suffices to check this after base-changing to $\Spec(\widehat{\sO}_{X,s})$ for a closed point $s\in X$. By \cite[Remark 2.4.9]{Posva_Singularities_of_quotients_by_foliations}, there exists formal coordinates $x,y,z$ for $\widehat{\sO}_{X,s}$ such that $\sF\otimes \widehat{\sO}_{X,s}$ is generated by either
		$$\frac{\partial}{\partial x}
		\quad \text{or}\quad
		x\frac{\partial}{\partial x}+y\frac{\partial}{\partial y}
		\quad \text{or}\quad
		x\frac{\partial}{\partial x}+y\frac{\partial}{\partial y}+z\frac{\partial}{\partial z}.$$
In the first case, $\sF$ is actually regular at $s$. In the second case the singular locus $Z$ is cut out by $(x,y)$ and in the third case it is cut out by $(x,y,z)$; so $Z$ is indeed a regular subscheme of $X$. In particular $\Bl_{Z}X$ is regular. It follows at once from the computations in \autoref{example:smth_blowup_of_foliations} that the formal local generators of $\sF$ lift to derivations which generates, after saturation, a regular foliation. This concludes the rank $1$ case.
		
\medskip
\textsc{Case 2: $\sF$ has rank $2$.} Consider the quotient $Y^{(1)}=X/\sF$ and the factorization
		$$F_{Y/k}\colon Y\overset{h}{\longrightarrow} X\longrightarrow Y^{(1)}.$$
Then $h$ is defined by a $1$-foliation $\sG$ of rank $1$. By the first case, we may assume that $Y$ is regular and that $\sG$ is everywhere regular. Then $X$ is regular by one implication of \cite[Lemma 2.5.10]{Posva_Singularities_of_quotients_by_foliations}, and as $Y^{(1)}$ is also regular we deduce that $\sF$ is regular by the other implication of \cite[Lemma 2.5.10]{Posva_Singularities_of_quotients_by_foliations}.
\end{proof}

\appendix
\section{Globalization of Cossart's algorithm}\label{appendix}
In this appendix we prove \autoref{claim:Cossart_method_applies}, which is necessary for the proof of \autoref{thm:resolution_3fold_rank_2}. This claim is a slight generalization of the main result of \cite{Cossart_Polyhedre_caracteristique}. Following the referees' suggestion, we give a general overview of the main definitions, tools and ideas of \emph{op.cit.}\ before proving our claim.

\subsection{Definitions}
Let $k$ be an algebraically closed field of characteristic $p>0$, and $X$ be a regular variety of dimension three over $k$. We fix a global function $f\in H^0(X,\sO_X)$ which is not a $p$-th power.

\subsubsection{Some ideals}\label{section:appendix_def}
Given a snc reduced divisor $E$ on $X$ and $Y\subset X$ a regular closed sub-scheme having normal crossings with $E$, we define the sheaf
		$$\sD(X,E,Y)\colon U\mapsto 
		\{\partial\in \Der(\sO_X)(U)\mid 
		\partial(I_{E\cap U})\subseteq I_{E\cap U} 
		\ \text{and} \
		\partial(I_{Y\cap U})\subseteq I_{Y\cap U}\}.$$
Evaluating such a sheaf of differential operators at $f$, we obtain the ideal sheaf
		$$\sJ(X,f,E,Y)=\im\left[
		\sD(X,E,Y)\overset{\ev_f}{\longrightarrow} \sO_X\right].$$
In case $Y$ is empty, we simply suppress it from the notations.

By \cite[Proposition I.A.2]{Cossart_Polyhedre_caracteristique} there is a unique way of writing
		$$\sJ(X,f,E,Y)=
		\sH(X,f,E)\cdot J(X,f,E,Y)$$
where $\sH(E,f,E)$ and $J(X,f,E,Y)$ are ideal sheaves uniquely with the following properties: $\sH(X,f,E)$ is the ideal of a snc divisor whose support is contained in $E$, and for every generic point $\eta$ of $E$ we have $\sH(X,f,E)_\eta=\sJ(X,f,E,Y)_\eta$. The ideal $\sH(X,f,E)$ depends on $(X,f,E)$, but not on $Y$ \cite[Proposition I.A.3]{Cossart_Polyhedre_caracteristique}. 

For $x\in X$ we let:
		$$\alpha(Y,x)=\ord_x J(X,f,E,Y), \
		\alpha(x)= \ord_x J(X,f,E,\overline{\{x\}}), \
		\nu(x)=\ord_x J(X,f,E).$$
We always have $\nu(x)\leq \alpha(Y,x)\leq \alpha(x)\leq \nu(x)+1$ \cite[Proposition I.A.3]{Cossart_Polyhedre_caracteristique}.
If $Y$ is irreducible with generic point $\eta$, we write $\alpha(Y)=\alpha(Y,\eta)$.
%Note that if $Y$ is a stratum of $E$, then $J(X,f,E,Y)=J(X,f,E)$ and $\alpha(Y,x)=\nu(x)$ for $x\in Y$ \cite[I.A.6.(6)]{Cossart_Polyhedre_caracteristique}.

For any $n\geq 0$ we set
		$$\Sing_n(X,f,E)=\{x\in X\mid \nu(x)\geq n\}.$$
We write $\Sing(X,f,E)=\Sing_1(X,f,E)$, and we let $\nu(X)=\sup\{\nu(x)\mid x\in X\}$.

%Given a point $x\in X$, Cossart also defines the ideals $\sJ(X,f,(u,\lambda))$ and $I(X,f,(u,\lambda))$ in $\sO_{X,x}$ where $(u,\lambda)$ is a $p$-basis at $x$ satisfying some conditions \cite[I.A.6]{Cossart_Polyhedre_caracteristique}. We shall only say that these conditions depend only on $E$ and $Y$ \cite[Definition I.A.1.1]{Cossart_Polyhedre_caracteristique}, that $\sJ(X,f,(u,\lambda))$ is generated by some partial derivatives of $f$ induced by the $p$-basis \cite[I.A.1.(15-17), A.6.(1)]{Cossart_Polyhedre_caracteristique} and that $I(X,f,(u,\lambda))=\sH(X,f,E)^{-1}\cdot \sJ(X,f,(u,\lambda))$.

\subsubsection{Blow-ups}\label{section:appendix_blow_up}
We allow ourselves to blow-up $(X,E)$ in a very restrictive way. Let us call a sequence
		$$ (X(n), E(n))\to (X(n-1),E(n-1))\to \dots \to (X(0),E(0))=(X,E)$$
a \emph{modification} of $(X,E)$ if the following conditions, defined inductively on the sequence, are satisfied:
	\begin{enumerate}
		\item we denote by $f\in H^0(X(i),\sO_{X(i)})$ the pullback of $f$;
		\item each $\pi(i)\colon X(i+1)\to X(i)$ is the blow-up of a regular closed sub-scheme $Y(i)\subset \Sing(X(i),f,E(i))$ having normal crossings with $E(i)$;
		\item $E(i+1)= \red\left[\pi(i)^{-1}(E(i)\cup Y(i))\right]$ is snc.
	\end{enumerate}
We write $\Sing(X(i))=\Sing(X(i),f,E(i))$ for simplicity. If a point $x(i+t)\in X(i+t)$ maps to $x(i)\in X(i)$ along this sequence of blow-ups, we say that $x(i+t)$ is \emph{$\nu$-close} to $x(i)$ if $\nu(x(i+t))=\nu(x(i))$.

Most of the time, we will consider blow-ups $\pi(i)$ satisfying the additional condition
	\begin{enumerate}\setcounter{enumi}{3}
		\item for every $x\in Y(i)$, if $\eta$ is the generic point of the irreducible component of $Y(i)$ passing through $x$, then $\alpha(Y(i),x)=\alpha(Y(i),\eta)$.
	\end{enumerate}
When this additional condition is satisfied, $\pi(i)$ will be called a \emph{permissible blow-up}.

\subsubsection{Goal}\label{section:appendix_goal}
Cossart proves the following: if $\dim X=3$, there exists a modification $X(n)$ of $(X,f,\emptyset)$ such that $J(X(n),f,E(n))=\sO_{X(n)}$. Equivalently, we eventually have $\nu(X(n))=0$. This implies, by \cite{Giraud_Forme_normale_sur_une_surface_en_car_pos}, that \'{e}tale-locally the function $f$ has the form $t^p+\bold{x}^{\bold{a}}$ as in \autoref{eqn:etale_neighbourhoods_II}. 

%Our goal is to show that Cossart's procedure depends only on the divisors $E(i)$ and on the ideals $\sJ(X(i),f,E(i))$ and $\sJ(X(i),f,E(i),Y(i))$ encountered along the modification.
%is insensitive to the change $f\rightsquigarrow u^p f$ where $u\in H^0(X,\sO_X^\times)$.

\subsubsection{Transformation laws}
Let us explain how the ideal sheaves defined above behave under modifications and permissible blow-ups:
	\begin{enumerate}
		\item If $\pi(i)\colon X(i+1)\to X(i)$ is a modification, with blow-up center $Y(i)$ irreducible, then \cite[Proposition I.E.1]{Cossart_Polyhedre_caracteristique}:
			$$\begin{cases}
			\sJ(X(i+1),f,E(i+1))=& \sJ(X(i),f,E(i),Y(i))\cdot \sO_{X(i+1)},\\
			\sH(X(i+1),f,E(i+1))=& \sH(X(i),f,E(i),Y(i))\cdot I_{Y(i)}^{\alpha(Y(i))}\cdot \sO_{X(i+1)},\\
			J(X(i+1),f,E(i+1))=& I_{Y(i)}^{-\alpha(Y(i))}\cdot J(X(i),f,E(i),Y(i))\cdot \sO_{X(i+1)}.\\
			\end{cases} $$
		\item Suppose that $\pi(i)$ is a permissible blow-up. Then for any $x\in X(i+1)$ we have $\nu(x)\leq \nu(\pi(i)(x))$ \cite[Proposition I.E.2]{Cossart_Polyhedre_caracteristique}.
	\end{enumerate}
	
\subsubsection{Stability}
We record the following result:

\begin{claim}\label{claim:gluing_together_the_ideals}
Let $\phi\in H^0(X,\sO_X^\times)$. Then the ideals $\sJ$, $\sH$, $J$, the invariants $\alpha$ and $\nu$, the singular loci, the notions of modifications and permissible blow-ups, and the transformations laws, are left unchanged by the substitution $f \rightsquigarrow \phi^pf$.
\end{claim}
\begin{proof}\renewcommand{\qedsymbol}{$\lozenge$}
It suffices to show that $\sJ(X,f,E,Y)=\sJ(X,\phi^pf,E,Y)$. Let $\partial\in \sD(X,E,Y)$: then $\partial(\phi^pf)=\phi^p\partial(f)$. As $\phi^p$ is invertible, this implies that $\partial(\phi^pf)\in \sJ(X,f,E,Y)$ and $\partial(f)\in \sJ(X,\phi^pf,E,Y)$, so the desired equality follows.
\end{proof}

\subsection{Reduction to $\dim \Sing_{\nu(X(n))}X(n)\leq 1$}\label{section:appendix_reduction_to_dim_1}
We continue with the notations of \autoref{section:appendix_def}.
Suppose that a modification $X(i)$ of $X$ is given. Let $S(2,i)$ be the union of the codimension $2$ irreducible components of $\Sing_{\nu(X(i))}X(i)$. Then \cite[Propositions II.B.2-3]{Cossart_Polyhedre_caracteristique}:
	\begin{enumerate}
		\item $S(2,i)$ is regular, and $\Sing X(i)=S(2,i)$ along $S(2,i)$;
		\item if $E(i)$ has normal crossings with $S(2,i)$, and if $\pi(i)\colon X(i+1)\to X(i)$ is the blow-up of $S(2,i)$, then $\nu(X(i+1))\leq \nu(X(i))$, and if equality occurs then $S(2,i+1)$ is empty. (While $\pi(i)$ is an isomorphism, the divisors $E(i)$ and $E(i+1)$ are not the same!)
	\end{enumerate}
The blow-up $\pi(i)$ is not necessarily permissible, but this is the only situation where we will consider non-permissible blow-ups \cite[Remarque II.B.3.1]{Cossart_Polyhedre_caracteristique}.

%Assume we can find a modification $X(i+t)$ of $X(i)$ where $E(i+t)$ is transverse to $S(2,i+t)$. As we want to decrease $\nu(X(\bullet))$, the above shows that we can reduce to the study of the case $\dim \Sing_{\nu(X(i))}X(i)\leq 1$.

The transversality of $E(i)$ and $S(2,i)$ can be achieved through repeated permissible blow-ups along the intersection of the two divisors \cite[II.B.3-10]{Cossart_Polyhedre_caracteristique}. The sequence of blow-ups is dictated by the algorithm \cite[Algorithme II.B.10.8]{Cossart_Polyhedre_caracteristique} where the relevant ingredients are: transversality, the invariants $m$ and $\sigma$ \cite[II.B.8.2-3]{Cossart_Polyhedre_caracteristique}, and the conditions $\star$ and $\star\star$ \cite[II.B.5.1, II.B.10.2]{Cossart_Polyhedre_caracteristique}. 

We have the following stability result.

\begin{claim}\label{claim:m_sigma_S(2,i)}
Let $\phi\in H^0(X,\sO_X^\times)$ and fix a modification $X(i)$. Then over $X(i)$ the invariants $m$ and $\sigma$, and the conditions $\star$ and $\star\star$, are left unchanged by the substitution $f\rightsquigarrow \phi^pf$.
%, assuming $X(i)\to X$ is also invariant by this substitution.
\end{claim}
\begin{proof}\renewcommand{\qedsymbol}{$\lozenge$}
Indeed, the two latter conditions are formulated in terms of the divisors $E(i),S(2,i)$ and of the notion of permissibility. The quantities $m$ and $\sigma$ are invariants of irreducible embedded curves $C\subset X(i)$, and by definition they depend only on $E(i)$ and on the history of the modification. 
\end{proof}

\subsection{The algorithm}
We keep the notations of \autoref{section:appendix_def}.
Recall that we wish to find a modification $X(n)$ with $\nu(X(n))=0$. The previous sub-sections show that, by working inductively, it suffices to solve the following problem: given a modification $X(i)$ with $\dim \Sing_{\nu(X(i))}X(i)\leq 1$, find a sequence of permissible blow-ups $X(i+n)\to X(i)$ such that $\nu(X(i+n))<\nu(X(i))$. 

Note that if a permissible blow-up $X(i+1)\to X(i)$ satisfies $\nu(X(i+1))=\nu(X(i))$, then we still have $\Sing_{\nu(X(i+1))}X(i+1)$ \cite[II.C.0.3]{Cossart_Polyhedre_caracteristique}.

Cossart designs an intricate inductive proof in order to show that this problem can be solved. The main additional invariant is the function $\kappa$, defined from every point of any modification $X(i)$ and taking values in $\{0,1,\dots,7\}$ (and being sometimes decorated by an additional letter). Following Cossart, we first define characterize the cases $\kappa\in \{0,1\}$ before succinctly describing the algorithm of the proof.

\subsubsection{$\kappa=0$}
For a closed point $x\in \Sing X(i)$, we set $\kappa(x)=0$ if \cite[Definition II.C.1]{Cossart_Polyhedre_caracteristique}: $x$ is isolated in $\Sing_{\nu(x)}X(i)$, and there is a sequence of permissible blow-ups
		$$X(i+n)\to \dots \to X(i)$$
such that:
	\begin{enumerate}
		\item the set $Y(i+t)\subset X(i+t)$ of points that are $\nu$-closed to $x$ is finite,
		\item $X(i+t+1)\to X(i+t)$ is the blow-up of $Y(i+t)$,
		\item $Y(i+n)$ is empty.
	\end{enumerate}
These are, of course, the easiest singularities to resolve. 
It is clear that:

\begin{claim}\label{claim:kappa=0}
Let $\phi\in H^0(X,\sO_X^\times)$. Then the condition $\kappa=0$ is left unchanged by the substitution $f \rightsquigarrow \phi^pf$.
\end{claim}

\subsubsection{$\kappa=1$}
For a given point $x\in \Sing X(i)$, let $\sC(x)$ be the set of closed irreducible curves $C$ through $x$ with $\alpha(C)\geq \nu(x)$ and $(\nu(C),m(C),\alpha(C),\sigma(C))$ maximal for the lexicographic order. Then the condition $\kappa(x)=1$ indicates the success of an algorithm defined using $\nu$, permissibility and cardinality of $\sC(\bullet)$ \cite[II.C.2]{Cossart_Polyhedre_caracteristique}. 

This is the generic case: we have $\kappa=1$ on a open dense open subset of the dimension one components of $\Sing_{\nu(X(i))}X(i)$ \cite[Proposition C.2.3]{Cossart_Polyhedre_caracteristique}.

\begin{claim}\label{claim:kappa=1}
Let $\phi\in H^0(X,\sO_X^\times)$. Then the condition $\kappa=1$ is left unchanged by the substitution $f \rightsquigarrow \phi^pf$.
\end{claim}
\begin{proof}\renewcommand{\qedsymbol}{$\lozenge$}
In view of what we have explained, this follows from \autoref{claim:gluing_together_the_ideals} and \autoref{claim:m_sigma_S(2,i)}.
\end{proof}

\subsubsection{Good and bad points}
Assume that we have defined the conditions $\kappa\in \{0,\dots,q\}$ with $q\leq 7$. Let $x\in X(i)$ with $\kappa(x)=q$. We say that $x$ is a \emph{good} point if the algorithm \cite[II.C.3,``Algorithme du point bon"]{Cossart_Polyhedre_caracteristique} is a success. This algorithm is defined in terms of permissibility, cardinality of $\sC(\bullet)$, $\nu$ and $\kappa\leq q$; its success indicates, roughly speaking, that we have a sequence of blow-ups of curves such that eventually no point is $(\nu,\kappa)$-close to $x$. If this algorithm fails, then $x$ is a \emph{bad} point. 
Obviously,

\begin{claim}\label{claim:good_point}
Let $\phi\in H^0(X,\sO_X^\times)$. Then the condition ``$x$ is a good point" is left unchanged by the substitution $f \rightsquigarrow \phi^pf$, assuming the conditions $\kappa\in \{0,\dots,\kappa(x)\}$ are.
\end{claim}

Cossart's strategy is to prove that bad points become good upon several rounds of blow-ups, and that further blow-ups make $(\nu,\kappa)$ drop (for the lexicographic order). In view of the definition of $\kappa=0$, this leads to the desired drop of $\nu$.

\subsubsection{The blow-up procedure}\label{section:appendix_Cossart_algorithm}
The sequence of permissible blow-ups used to make $\nu(X(i))$ drop is specified in \cite[Algorithme II.C.5.3]{Cossart_Polyhedre_caracteristique}. The blow-ups centers are selected using the notions of dimension, permissibility, good and bad points, $\nu$, $m$, $\alpha$, $\sigma$, $\kappa$ and $\sC$.

After reducing to $\dim \Sing_{\nu(X(i))}X(i)\leq 1$, the algorithm blow-ups points and curves. Let us note that we blow-up curves not contained in $\Sing_{\nu(X(i))}X(i)$ only once we have achieved $\dim \Sing_{\nu(X(i))}X(i)=0$. We need to ensure that if we blow-up a curve $C$ to improve the situation at a point $x$, then this will not mess up the situation at another point of $C\cap \Sing_{\nu(X(i))}X(i)$: this is shown in \cite[Proposition I.C.5.1]{Cossart_Polyhedre_caracteristique} (using only the tools introduced above).

\medskip
We claim now that the conditions $\kappa\in \{2,\dots,7\}$ are invariant under $f\rightsquigarrow \phi^p f$ for $\phi\in H^0(X,\sO_X^\times)$. We illustrate this in several important cases (including $\kappa=2$ and the ``joyful case"), and leave the remaining ones to the reader.

\subsection{The case $\kappa=2$}
\subsubsection{More notations}
We follow \cite[I.A.5]{Cossart_Polyhedre_caracteristique}. 
Let $x\in X(n)$ be an arbitrary point. We say that a $p$-basis $(\lambda_i)$ of $\sO_{X(n),x}$ is \emph{adapted} at $x$ if $\fm_{X,x}=(\lambda_1,\dots,\lambda_r)$ with $r\leq s$ and $I_{E(n),x}=(\lambda_1,\dots,\lambda_t)$ with $t\leq r$. Usually we write $u_i=\lambda_i$ for $i\leq r$, to distinguish the local parameters inside the $p$-basis, and let $(u,\lambda)=(\lambda_i)$. 
Since a $p$-basis of $\sO_{X(n),n}$ is also a differential basis, we can write
		$$df=\sum_{i=1}^s \D^{u,\lambda}_i(f) \cdot d\lambda_i
		\quad \text{in }
		\Omega_{\sO_{X(n),x}/k}^1.$$
Then we let $\DM^{u,\lambda}_i(f)=\lambda_i\D^{u,\lambda}_i(f)$ for $i=1,\dots,s$. See \cite[I.A.1]{Cossart_Polyhedre_caracteristique}.

By definition of a $p$-basis, we can write
		$$f=\sum_{\bold{a}<\bold{p}}f_{\bold{a}}^p\bold{\lambda}^{\bold{a}}, \quad
		f_{\bold{a}}\in \sO_{X(n),x}, \ \text{and} \
		\bold{\lambda}^{\bold{a}}=
		\prod_{i=1}^s\lambda_i^{\bold{a}(i)}.$$
We write $R(f,u,\lambda)=f_{\bold{0}}^p$.
We let $\ord_{u_i}$ indicate the order at the generic point of the divisor $\text{div}(u_i)$, and we set
		$$A(i)=\ord_{u_i}(f-f_{\bold{0}}^p), \quad
		h(u,\lambda)=\prod_{1\leq i\leq t}u_i^{A(i)}.$$
We also set 
	$$\sJ(X(n),f,(u,\lambda))=(\DM^{u,\lambda}_i(f)\mid 1\leq i\leq s)$$
and 
	$$I(X(n),f,(u,\lambda))=\sH(X(n),f,E(n))^{-1}\cdot \sJ(X(n),f,(u,\lambda));$$ 
these are ideals of $\sO_{X(n),x}$.

\begin{claim}\label{claim:A(i)}
Let $\phi\in H^0(X,\sO_X^\times)$. Then $R(\phi^pf,u,\lambda)=\phi^p R(f,u,\lambda)$. Moreover, the $A(i)$s, $h(u,\lambda)$, $\sJ(X(n),f,(u,\lambda))$ and $I(X(n),f,(u,\lambda))$ are invariant under the substitution $f\rightsquigarrow \phi^pf$.
\end{claim}
\begin{proof}\renewcommand{\qedsymbol}{$\lozenge$}
We have $\phi^pf=\sum_{\bold{a}<\bold{p}}(\phi f_{\bold{a}})^p\bold{\lambda}^{\bold{a}}$ so $R(\phi^pf,u,\lambda)=\phi^p R(f,u,\lambda)$. Moreover
		$$\ord_{u_i}(\phi^pf-\phi^pf_{\bold{0}}^p)=\ord_{u_i}(f-f_{\bold{0}}^p)$$
so $A(i)$ is invariant. Invariance of $h(u,\lambda)$ follows immediately, while the invariance of the ideals $\sJ$ and $I$ is clear.
\end{proof}

Now let $S=\bigoplus_{j=1}^N kU_j$ be a finitely generated $k$-vector space, and $\mathfrak{a}$ be an homogeneous ideal of $k[S]=k[U_1,\dots, U_N]$. There is a smallest sub-vector space $S_1\subset S$ such that 
		$$\mathfrak{a}=(k[S_1]\cap \mathfrak{a})k[S].$$
We let $\VDir(\mathfrak{a})=S_1$.
		
Let $(R,\fm,k)$ be a regular local ring of dimension $3$ with residue field $k$. Given any $g\in R$, let $v=\ord_\fm(g)$ and write $\cl_\fm^v(g)$ the image of $g$ in $\fm^v/\fm^{v+1}$. Given a non-trivial ideal $I$ of $R$, let $\cl^\bullet_\fm I$ be the ideal of $\gr_\fm(R)$ generated by the collection $\{\cl^v_\fm(g)\mid g\in I\}$. Choose a regular set of parameters $(u_1,u_2,u_3)$ for $R$, so that we can write $\gr_\fm(R)=k[U_1,U_2,U_3]$. Then we set $\VDir(I)=\VDir(\cl^\bullet_\fm I)$. See \cite[I.E.1.5]{Cossart_Polyhedre_caracteristique}.

We will apply this construction for $R=\sO_{X(n),x}$, where $x$ is a closed point, and $(u_i)$ a $p$-basis that is adapted at $x$. In this situation, we let \cite[I.E.2.1]{Cossart_Polyhedre_caracteristique}:
		$$\VDir(x)=
		\begin{cases}
			\VDir\left( \cl_\fm^\bullet 
						J(X(n),f,E(n),\overline{ \{x\} }) 
					\right) 
				& \text{if }\alpha(x)=\nu(x), \\
			\VDir\big( \cl_\fm^\bullet 
						J(X(n),f,E(n)) 
					\big) 
				& \text{if }\alpha(x)=\nu(x)+1.
		\end{cases}$$
It follows at once from \autoref{claim:gluing_together_the_ideals} that:
\begin{claim}\label{claim:VDir}
Let $\phi\in H^0(X,\sO_X^\times)$. Then $\VDir(x)$ is invariant under the substitution $f\rightsquigarrow \phi^pf$.
\end{claim}

\subsubsection{Definition of $\kappa=2$}
Let $X(n)$ be a modification such that $\dim \Sing_{\nu(X(n))} X(n)\leq 1$. Let $x\in \Sing X(n)$ be a closed point such that $\kappa(x)\notin \{0,1\}$. We say that $\kappa(x)=2$ \cite[IV.A.1]{Cossart_Polyhedre_caracteristique} if there exists a $p$-basis $(u,\lambda)$ adapted at $x$ such that:
	\begin{enumerate}
		\item if $\alpha(x)=\nu(x)$, then $E(n)\subset \text{div}(u_2u_3)$ and $\VDir(x)\not\subset kU_1\oplus kU_2$;
		\item if $\alpha(x)=\nu(x)+1$, then $E(n)=\text{div}(u_3)$ and:
			\begin{enumerate}
				%\item letting $G_2=\cl^v\left(
						%\frac{\D^{u,\lambda}_2(f)}{h(u,\lambda)}\right)$ we have
						%$$G_2(U_1,U_2,0)\in k[U_1^p,U_2];$$
				\item letting $G_2=\cl^v\left(
						\frac{\D^{u,\lambda}_2(f)}{h(u,\lambda)}\right)$, we have $U_1\in \VDir(G_2(U_1,U_2,0))$;
				\item $\frac{\D^{u,\lambda}_1(f)}{h(u,\lambda)}\in u_3\sO_{X(n),x}$.
			\end{enumerate}
	\end{enumerate}
(Here we have made use of the remark immediately following \cite[Definition IV.A.1]{Cossart_Polyhedre_caracteristique} to shorten the definition.)

\begin{claim}\label{claim:kappa=2}
Let $\phi\in H^0(X,\sO_X^\times)$. Then the condition $\kappa=2$ is left unchanged by the substitution $f \rightsquigarrow \phi^pf$.
\end{claim}
\begin{proof}\renewcommand{\qedsymbol}{$\lozenge$}
The condition $\kappa\notin \{0,1\}$ is invariant, and so are the quantities $\alpha$ and $\nu$. We take the same $p$-bais $(u,\lambda)$. By \autoref{claim:A(i)} and \autoref{claim:VDir}, $\VDir(x)$ and $h(u,\lambda)$ are invariant. So it suffices to check that the invariance of the conditions (i) and (ii) above. We have $D^{u,\lambda}_i(\phi^pf)=\phi^pD^{u,\lambda}_i(f)$, so invariance of (ii) is clear. For (i), simply observe that
		$$\cl^v\left(
						\frac{\D^{u,\lambda}_2(\phi^pf)}{h(u,\lambda)}\right)
				=\gamma^p G_2$$
where $\gamma=\cl^0(\phi)\in k$, and that $\VDir(G_2(U_1,U_2,0))=\VDir(\gamma^pG_2(U_1,U_2,0))$.
\end{proof}

\subsection{The ``joyful case" and some $\kappa=3$ cases}
We fix a modification $X(i)$ with $\dim \Sing_{\nu(X(i))}X(i)\leq 1$. We review the definition of $\kappa$ in a few more cases.

\subsubsection{More notations (bis)}
Let $x\in X(i)$ be a closed point, and $(u,\lambda)$ be a $p$-basis adapted at $x$. We let \cite[I.E.3]{Cossart_Polyhedre_caracteristique}:
		$$v(x)=\dim_k \VDir(x).$$
Moreover \cite[I.E.4]{Cossart_Polyhedre_caracteristique}, if $I$ is an ideal of $\sO_{X(i),x}$, we let $e(I)$ denote the smallest $m\geq 0$ such that there exists a $p$-basis $(u,\lambda)$ adapted at $x$ with $\VDir(I)\subset \bigoplus_{j=1}^m kU_j$. Then we set:
		$$e(x)=
		\begin{cases}
			e(J(X(i),f,E(i), \{x\}) & \text{if } \alpha(x)=\nu(x), \\
			e(J(X(i),f,E(i)) & \text{if } \alpha(x)=\nu(x)+1.
		\end{cases}$$
It is clear that:
\begin{claim}\label{claim:v_and_e}
Let $\phi\in H^0(X,\sO_X^\times)$. Then $v(x)$ and $e(x)$ are invariant under the substitution $f\rightsquigarrow \phi^pf$.
\end{claim}

\subsubsection{The joyful case}
We say that a closed point $x\in X(i)$ is a \emph{joyful case} if \cite[VI, Introduction]{Cossart_Polyhedre_caracteristique}:
	\begin{enumerate}
		\item $v(x)=1$, $e(x)=2$, $\alpha(x)=\nu(x)$;
		\item there exists a $p$-basis $(u,\lambda)$ adapted at $x$ such that $\VDir(x)=kU_1\oplus kU_2$ and $\text{div}(u_1u_2)\subseteq E(i)$;
		\item let $x'$ be a point on the blow-up $\pi\colon \Bl_xX(i)\to X(i)$ that is $\nu$-closed to $x$ and not on the strict transform of $V(u_1,u_2)$. Then there is a $p$-basis $(v,\mu)$ adapted at $x'$ such that
				$$\begin{cases}
				f &= (v_2^av_3^b)^p (\gamma v_1^{1+\nu(x)}+v_2g)+ R(f,v,\mu),\\
				\text{div}(v_2)&=\pi^{-1}(x), \\
				v_1v_2&=u_1+u_2, \\
				1+\nu(x)&\neq 0 (p),
				\end{cases}$$
where $a,b\geq 0$ and $\gamma\in \sO_{\Bl_xX(i),x'}^\times$ and $g\in \sO_{\Bl_xX(i),x'}^\times$.
	\end{enumerate}
	
\begin{claim}\label{claim:joyful_case}
Let $\phi\in H^0(X,\sO_X^\times)$. Then the definition of the joyful case is invariant under the substitution $f\rightsquigarrow \phi^pf$.
\end{claim}
\begin{proof}\renewcommand{\qedsymbol}{$\lozenge$}
Thanks to \autoref{claim:gluing_together_the_ideals}, to \autoref{claim:VDir} and to \autoref{claim:v_and_e}, we only have to check the claim for the condition (c) above. We take the same $(v,\mu)$. Multiplying the first equality by $\phi^p$, we get
		$$\phi^p f = (v_2^av_3^b)^p ((\phi^p\gamma) v_1^{1+\nu(x)}+v_2(\phi^p g))+ \phi^p R(f,v,\mu).$$
As $\phi^p R(f,v,\mu)=R(\phi^p f, v, \mu)$ by \autoref{claim:A(i)}, we are done.
\end{proof}

\subsubsection{$\kappa=3(A)$}
We say that $x\in X(i)$ satisfies $\kappa(x)=3(A)$ if \cite[VI.A.3]{Cossart_Polyhedre_caracteristique}:
	\begin{enumerate}
		\item $\nu(x)\geq 1$ and $\kappa(x)>2$,
		\item there is a $p$-basis $(u,\lambda)$ adapted at $x$ such that $\text{div}(u_1)\subseteq E(n)$, $A(1)=0$ modulo $p$ and
			$$\ord_x\left( 
			\frac{\D^{u,\lambda}_1(f)}{h(u,\lambda)} \right) 
			= \nu(x)-1.$$
	\end{enumerate}
Clearly,
\begin{claim}\label{claim:kappa=3A}
Let $\phi\in H^0(X,\sO_X^\times)$. Then the definition of $\kappa=3(A)$ is invariant under the substitution $f\rightsquigarrow \phi^pf$.
\end{claim}

\subsubsection{$\kappa=3(B)$}
We say that $x\in X(i)$ satisfies $\kappa(x)=3(A)$ if \cite[VI.B.1]{Cossart_Polyhedre_caracteristique}:
	\begin{enumerate}
		\item $\nu(x)=\alpha(x)$ and $\kappa(x)>2$,
		\item there is a $p$-basis $(u,\lambda)$ adapted at $x$ such that there is $\partial\in \sD(X(i),E(i))$ with
				$$\frac{\partial(f)}{h(u,\lambda)}
				= u_2 u_3^{\nu(x)-1} \mod u_2(u_1,u_2),$$
		and moreover one of the following two conditions is satisfied:
			\begin{enumerate}
				\item $E(i)=\text{div}(u_1)$ and
				$\VDir\left( \cl^v[ J(X(i),f,E(i)) \mod (u_1)] 
				\right) = kU_2\oplus kU_3$, or
				\item $E(n)=\text{div}(u_1u_2)$, $\text{div}(u_1)>\text{div}(u_2)$ for the history of the modification (see \cite[I.B.8]{Cossart_Polyhedre_caracteristique}), and
						$$f=h(u,\lambda)(\gamma u_1^{\nu(x)}+u_2A)+ R(f,u,\lambda)$$
				where $\gamma\in \sO_{X(i),x}^\times$, $A\in \fm_x^{\nu(x)-1}$ and $\nu(x)\geq 2$.
			\end{enumerate}				
	\end{enumerate}

\begin{claim}\label{claim:kappa=3B}
Let $\phi\in H^0(X,\sO_X^\times)$. Then the definition of $\kappa=3(B)$ is invariant under the substitution $f\rightsquigarrow \phi^pf$.
\end{claim}
\begin{proof}\renewcommand{\qedsymbol}{$\lozenge$}
We take the same $(u,\lambda)$. In (b) we replace $\partial$ by $\phi^{-p}\partial$. For (b.ii) we substitute $\gamma \rightsquigarrow \phi^p\gamma$, $A\rightsquigarrow \phi^p A$ and use \autoref{claim:A(i)}.
\end{proof}

\subsection{The case $\kappa=4(A)$}
We treat one last case, as its definition requires the use of characteristic polyhedra---an important tool in \cite{Cossart_Polyhedre_caracteristique} which has not been featured yet in our overview.

\subsubsection{More notations (ter)}
We follow \cite[IV.C.1-2]{Cossart_Polyhedre_caracteristique}. 
Let $R$ be a regular local ring with and $u_1,\dots,u_d$ be a regular system of parameters. For a convex subset $\Delta\subset \bR^d_+$ such that $\Delta=\Delta+\bR^d_+$, let $I(\Delta)_u$ be the ideal of $R$ generated by the collection $\{\prod_i u_i^{A(i)}\mid A\in \Delta\}$.

Let $I$ be an ideal of $R$. Let $v=\ord_{u_1}(J\cdot R/(u_2,\dots,u_d))$. We denote by $\Delta[J;u_d,\dots,u_2;u_1]$ the smallest convex subset $\Delta$ of $\bR^{d-1}_1$ satisfying the following property: if $D_\Delta\subset \bR^d_+$ is the convex body defined by
		$$D_\Delta=\left\{ (x_1,\dots,x_d) \mid
			x_1\geq v \text{ or }
			\frac{1}{v-x_1}(x_2,\dots,x_d)\in \Delta\right\}$$
then $J\subseteq I(D_\Delta)_u$.

\subsubsection{Definition of $\kappa=4(A)$}
Let $x$ be a closed point of $\Sing X(i)$. We say that $\kappa(x)=4(A)$ if:
	\begin{enumerate}
		\item $\kappa(x)>3$;
		\item $\alpha(x)=1+\nu(x)$ and $1+\nu(x)\neq 0(p)$;
		\item $E(i)\neq \emptyset$ around $x$;
		\item there is a $p$-basis $(u,\lambda)$ adapted at $x$ such that:
			\begin{enumerate}
				\item $\ord_x [I(X(i),f,(u,\lambda))\mod (u_2,u_3)]=\nu(x)+1$;
				\item $\Sing_{\nu(x)}X(i)\subseteq \text{div}(u_1u_2)\cup \text{div}(u_2u_3)$ and $\text{div}(u_1)\not\subset E(i)$;
				\item $\Delta[I(X(i),f,(u,\lambda)):u_3,u_2:u_1]$ has only one vertex $(a,b)$ with $a+b<1$;
				\item $\Delta\left[ \left(
					\frac{\D^{u,\lambda}_1(f)}{h(u,\lambda)} \right);
					u_3,u_2;u_1 \right]$ is prepared \cite[IV.C.3.3]{Cossart_Polyhedre_caracteristique}.
			\end{enumerate}
	\end{enumerate}
	
\begin{claim}\label{claim:kappa=4A}
Let $\phi\in H^0(X,\sO_X^\times)$. Then the definition of $\kappa=4(A)$ is invariant under the substitution $f\rightsquigarrow \phi^pf$.
\end{claim}
\begin{proof}\renewcommand{\qedsymbol}{$\lozenge$}
We take the same $(u,\lambda)$. Then invariance is clear for items (a),(b),(c),(d.i-ii). For (d.iii-iv), note that the ideals involved are invariant, so the associated convex bodies are also invariant.
\end{proof}

\subsection{Proof of \autoref{claim:Cossart_method_applies}}
We change our setting slightly. Let $X$ be a regular variety of dimension three over $k$, and let $u\in K(X)$ be a rational function which does not belong to $K(X)^p$. Let $D$ be the divisor of poles of $u$, and let $\{U_i\}$ be an affine open cover of $X$, sufficiently refined so that we can find $v_i\in \sO_{X}(U_i)$ with $\text{div}_{U_i}(v_i)=D\cap U_i$: this is possible since $X$ is regular. Then for all $i,j$:
		$$uv_i^p\in \sO_X(U_i), \quad
		v_i/v_j\in \sO_X(U_{ij})^*$$
with the usual notation $U_{ij}=U_i\cap U_j$. 

We claim that there exists a projective birational morphism $f\colon V\to X$ where $V$ is regular such that for each $i$, \'{e}tale-locally on $f^{-1}(U_i)$ we have
		$$f^*(uv_i^p)=t^p+\bold{x}^{\bold{a}}$$
where $\bold{x}$ is a $p$-basis and $\bold{a}\neq\bold{0} (p)$. In fact, we claim that Cossart's algorithm can be applied to find such a model $V$.

Over $U_{ij}$ we have $(uv_i^p)=(v_j/v_i)^p(uv_j^p)$. Thus by \autoref{claim:gluing_together_the_ideals} we see that: if $\pi\colon X'\to X$ is any proper birational morphism, with $X'$ regular, $E'$ an snc divisor on $X'$ and $Y'$ a regular closed sub-scheme having normal crossing with $E'$, the ideals $\sJ(\pi^{-1}(U_i),uv_i^p,E',Y')$ glue together into an ideal sheaf $\sJ(X',\{uv_i^p\},E',Y')$. As noticed in \autoref{claim:gluing_together_the_ideals}, over $U_{ij}$ the invariant $\nu$ can be defined using either $uv_i^p$ or $uv_j^p$, thus it is well-defined globally. We define modifications of $(X,\{uv_i^p\},\emptyset)$, and the associated snc divisors, as in \autoref{section:appendix_blow_up}.

If after blowing-up we achieve $\nu=0$ everywhere, then as mentioned in \autoref{section:appendix_goal} we will have reached our goal.

We apply Cossart's algorithm, which we briefly described in \autoref{section:appendix_reduction_to_dim_1} and \autoref{section:appendix_Cossart_algorithm}. The point is that the selection of the blow-up centers over $U_{ij}$ are insensitive to the change $uv_i^p\rightsquigarrow uv_j^p$, because it depends on notions that are themselves insensitive to such substitutions as our successive claims have shown. Thus we can run Cossart's algorithm globally over $X$. It terminates by the argument given in \cite[II.C.5.4]{Cossart_Polyhedre_caracteristique}---which, again, only uses notions that are insensitive to the substitutions $uv_i^p\rightsquigarrow uv_j^p$. This completes our proof.

\bibliographystyle{alpha}
\bibliography{Bibliography}

\end{document}